\documentclass[journal,onecolumn]{IEEETran}
\let\labelindent\relax

\usepackage{etoolbox}
\usepackage[shortlabels]{enumitem}
\usepackage{cite}
\usepackage{amsmath,amssymb,amsfonts,mathtools}
\usepackage{algorithmic}
\usepackage{graphicx} 
\usepackage{graphics} 
\usepackage{textcomp}
\usepackage{dsfont}
\usepackage{color}

\usepackage{epstopdf}        
\usepackage{lscape}
\usepackage{multirow}
\usepackage{calligra}
\usepackage{booktabs}
\usepackage{framed} 
\usepackage{picins}
\usepackage{empheq}
\usepackage{epsfig} 
\usepackage{times}
\usepackage{dblfloatfix}
\usepackage{blindtext}
\usepackage{comment}

\usepackage{diffcoeff}
\diffdef { p }
{
op-symbol = \partial ,
left-delim = \left .,
right-delim = \right | ,
subscr-nudge = +4 mu
}

\newtheorem{thm}{Theorem}
\newtheorem{prop}{Proposition}
\newtheorem{lemma}{Lemma}

\newtheorem{definition}{Definition}
\newtheorem{assumption}{Assumption}
\newtheorem{remark}{Remark}
\newtheorem{example}{Example}

\newcommand{\tcb}{\textcolor{black}}

\newcommand{\set}[1]{\left\{#1\right\}}
\renewcommand{\vec}{\mathbf}
\newcommand{\R}{\mathbb{R}}
\newcommand{\N}{\mathbb{N}}

\newcommand{\abs}[1]{\left\lvert#1\right\rvert}
\renewcommand{\vec}{\mathbf}
\newcommand{\dom}[1]{\text{dom}\left(#1\right)}

\newcommand*\tageq{\refstepcounter{equation}\tag{\theequation}}
\allowdisplaybreaks

\newcommand*{\QEDB}{\hfill\ensuremath{\square}}%
\graphicspath{'./Figs'}
\date{December 3, 2022}

\begin{document}
\title{Momentum-Based Nash Set Seeking over Networks via  Multi-Time Scale Hybrid Dynamic Inclusions}

\author{Daniel E. Ochoa and Jorge I. Poveda\vspace{-0.86cm}
\thanks{\tcb{The authors are with the Department of Electrical and Computer Engineering at the University of California San Diego, CA, USA. E-mail: {\tt dochoatamayo@ucsd.edu}.}}}
\maketitle

\begin{abstract}
Multi-time scale techniques, such as singular perturbations and averaging theory, have played an important role in the development of distributed Nash equilibrium seeking algorithms for network systems. Such techniques intrinsically rely on the \emph{uniform asymptotic stability} properties of the dynamics that evolve in each of the time scales of the closed-loop system. When such properties are absent, the synthesis of multi-time scale Nash equilibrium seeking algorithms is more challenging and it requires additional regularization mechanisms. In this paper, we investigate the synthesis and analysis of these mechanisms in the context of \emph{accelerated pseudogradient flows} with time-varying damping in non-cooperative games. Specifically, we introduce a new class of \emph{distributed and hybrid} Nash set seeking (NSS) algorithms that synergistically combine dynamic momentum-based flows with \emph{coordinated} discrete-time resets. The reset mechanisms can be seen as restarting techniques that allow individual players to choose their own momentum restarting policy to achieve better transient performance. The resulting closed-loop system is modeled as a hybrid dynamic inclusion, which is analyzed using tools from hybrid dynamical system's theory. Our algorithms are developed for potential games, as well as for monotone games for which a potential function does not exist. They can be implemented in games where players have access to gradient Oracles with full or partial information, as well as in games where players have access only to \emph{measurements} of their costs. In the latter case, we use tools from  hybrid extremum seeking control to achieve accelerated model-free Nash set seeking.
\end{abstract}
\vspace{-0.1cm}

\vspace{-0.2cm}
\section{Introduction}
\label{sec:introduction}
%
\IEEEPARstart{A}{mong} 
%
the different notions of equilibria related to game-theoretic models, the notion of Nash equilibrium (NE), introduced in \cite{NashPaper}, has become ubiquitous in many engineering and socio-technical systems such as transportation systems and energy markets. To converge to this equilibrium, different deterministic and stochastic Nash equilibrium seeking (NES) algorithms have been developed during the last decades, see \cite{BasarDNGT,MardenSafeExperimentationJournal,Grammatico21,StankovicNashSeeking,rosen,PersisGrammatico,Altman02,GNE19}. For games defined over networks, multi-time scale techniques that integrate fast consensus mechanisms have also become ubiquitous in the literature \cite{PavelGames}. Similarly, 
payoff-based and model-free NES dynamics, suitable for applications where precise mathematical models of the costs  are not available, have been studied using averaging theory for smooth ordinary differential equations (ODEs) \cite{Frihauf12a,Dither_ReUse} and for games with nonsmooth and hybrid dynamics \cite{PovedaKrsticBasar2020Journal,PoTe16}. 

\vspace{0.1cm}
\tcb{In the context of game-theoretic control system design, many results in the literature are somehow inspired or related to the time-invariant pseudogradient (PSG) flows studied by Rosen} in \cite{rosen}, which take the form $\dot{q}=-\mathcal{G}(q)$, where $\mathcal{G}$ is the \emph{pseudogradient} vector of the game \cite[Eq. (3.9)]{rosen}, and $q\in\mathbb{R}^n$ is the vector of actions of the players. For example, it is well-known that in potential games PSG flows can robustly minimize the potential function at a rate of order $\mathcal{O}(1/t)$. Additionally, for strongly monotone games, pseudogradient flows can achieve NES with exponential rates of convergence of order $\mathcal{O}(e^{-\kappa})$, with $\kappa$ being the strong monotonicity coefficient of $\mathcal{G}$. These stability and convergence results have become instrumental for the design of extended NES algorithms that incorporate additional mechanisms based on fast consensus dynamics \cite{GNE19,PavelGames}, projections \cite{Frihauf12a}, inertia \cite{yi2019operator}, proportional feedback terms \cite{PersisGrammatico}, tracking terms \cite{Hu15}, adaptive dynamics \cite{Grammatico21a}, etc. See also the recent work \cite{GaoCDC20} and references therein. However, while these results have provided significant insight into the design of NES dynamics, existing results still suffer from fundamental transient limitations inherited from PSG flows, which can be further exacerbated in games where the cost functions have shallow monotonicity properties.
 
On the other hand, compared to PSG flows, time-varying momentum-based dynamics, which are common in the optimization and machine learning literature  \cite{ODE_Nesterov,Wilson18,HamiltonianDescent,HighResolution2018,Muehlebach1,attouch2018rate,zero_order_poveda_Lina,OchoaPoveda20LCSS}, have not received as much attention in the context of games. In particular, in this paper we are interested in studying the Nash equilibria learning capabilities of the second-order dynamics
\begin{equation}\label{ODEmomentum}
\dot{q}=\frac{2}{\tau}(p-q),~~\dot{p}=-2\tau \mathcal{G}(q),~~~\dot{\tau}=\eta,
\end{equation}
with $\tau(0)=T_0\geq0$ and $\eta>0$, which are related to Nesterov's accelerated optimization algorithm whenever $\mathcal{G}$ is a gradient operator and via the \tcb{transformation to the momentum variable $p=\frac{\tau}{2}\dot{q}+q,~p(0)\in\R^n$}; see  \cite{ODE_Nesterov,Wilson18}. Such systems are particularly useful for optimization and estimation problems with cost functions having vanishing curvature at the optimal points, since they exhibit a geometric property, termed \emph{acceleration}, able to minimize smooth convex functions at a rate of order $\mathcal{O}(1/t^2)$. Moreover, in strongly convex optimization problems, systems of the form \eqref{ODEmomentum}, combined with suitable ``restarting'' heuristics, can achieve exponential rates of convergence; see \cite{ODE_Nesterov,Candes_Restarting}. Indeed, dynamics of the form \eqref{ODEmomentum} have been recently shown to accelerate the convergence in adaptive estimation problems \cite{annaswamytunners}, extremum seeking control \cite{zero_order_poveda_Lina}, and concurrent learning techniques \cite{OchoaL4DC}. Therefore, in light of intriguing numerical results in the context of games, it is natural to ask whether systems of the form \eqref{ODEmomentum} are also suitable for the \emph{robust} and \emph{efficient} solution of NES problems in noncooperative games when $\mathcal{G}$ is a pseudogradient operator, and also whether these dynamics can be extended to \tcb{network games and model-free setting}. In fact, existing results in the literature have focused only on (non-uniform) convergence results in centralized potential games \cite{gadjov2022exact}, or in momentum-based dynamics with maximally monotone operators $\mathcal{G}$ via Yosida regularizations of the form $\frac{1}{\lambda}(I-(I+\lambda \mathcal{G})^{-1})$, which are no suitable for distributed implementations \cite{attouch18}.

%

\vspace{0.01cm}
\textsl{Main Results:} In this paper, we provide answers to the above questions by using tools from nonlinear control theory. 
In particular, we first show that the direct implementation of systems of the form \eqref{ODEmomentum} is,
%
%
in general, not suitable for the efficient distributed solution of Nash set seeking (NSS) problems in non-cooperative games, even when $\mathcal{G}$ is strongly monotone and there exists a potential function. The limitations arise from three main structural issues: First, the dependence on a ``centralized'' momentum coefficient $\tau$ that precludes distributed implementations and that can also lead to \emph{uncoordinated} algorithms where players implement individual coefficients $\tau_i$, leading to poor transient performance, i.e., slower than PSG flows. Second, the vanishing damping in system \eqref{ODEmomentum} leads to lack of uniformity (with respect to $\tau(0)$) in the convergence properties of the dynamics, which makes them prone to instability under arbitrarily small additive disturbances, unavoidable in feedback-based implementations where gradients or neighboring states are estimated on-the-fly via multi-time scale techniques. Third, in non-potential games, standard Lyapunov functions used in optimization are not applicable, and, in fact, solutions of \eqref{ODEmomentum} might diverge even when the game is strongly monotone and the damping coefficient is uniformly lower bounded by a positive number. 

While the above features might suggest that momentum-based dynamics are problematic for games, it turns out that systems of the form \eqref{ODEmomentum} \emph{can} be used to efficiently and robustly find NE in a decentralized way, whenever they are combined with suitable distributed discrete-time dynamics that persistently reset some of the states of the players in a \emph{coordinated} way. However, in contrast to optimization problems \cite{ODE_Nesterov,Candes_Restarting,zero_order_poveda_Lina}, for general (non-potential) noncooperative games the frequency of the resets must occur in a certain frequency band in order to simultaneously achieve stability and acceleration. We establish these results using tools from hybrid dynamical systems (HDS) theory \cite{HDS}, and we extend the algorithms to decentralized network games and model-free settings via multi-time scale hybrid control theoretic tools.
The following original contributions are then presented in the paper.

\vspace{0.01cm}
\textsl{i)} We propose the first NSS algorithms with continuous-time \emph{dynamic} momentum and \emph{robust} asymptotic stability properties in non-cooperative games with $n$ players. The algorithms incorporate three main elements: a) a class of distributed \emph{continuous-time} pseudogradient-based dynamics with \emph{time-varying} momentum coefficients inspired by \eqref{ODEmomentum}; b) distributed periodic \emph{discrete-time} resets implemented by the players, which incorporate \emph{heterogeneous} reset policies that allow players to decide whether or not to restart their own momentum; c) a robust \emph{set-valued} distributed coordination mechanism that synchronizes the reset times of the players to induce suitable system-wide acceleration properties. The combination of the above three elements necessarily leads to a hybrid dynamic inclusion with non-unique solutions, which we analyze using tools from HDS \cite{HDS}.   
%

\textsl{ii)} To accommodate situations where players do not have access to full-information Oracles that provide evaluations of their pseudogradients, we introduce a new \emph{distributed} momentum-based  hybrid NSS algorithm for games with partial information, where players leverage communication with neighbors to estimate their actions on-the-fly in order to obtain gradient evaluations from the Oracle. The design of these dynamics follows similar multi-time scale ideas used for ODEs in the literature \cite{PavelGames}, but which are not directly applicable to systems of the form \eqref{ODEmomentum}. Indeed, unlike existing results based on fast consensus dynamics and ``reduced'' pseudogradient flows, our reduced dynamics are hybrid and set-valued, which prevents the direct application of standard singular perturbation tools for ODEs. 

%

\textsl{iii}) We present payoff-based versions of all our hybrid NSS algorithms, suitable for \emph{model-free} learning in non-cooperative games where players have access only to measurements of their cost. These dynamics exploit recent tools developed in the context of averaging-based hybrid extremum seeking control \cite{PoTe16,zero_order_poveda_Lina}, and their analyses is fundamentally different from other model-free non-hybrid algorithms studied in the literature, e.g. \cite{Frihauf12a,Dither_ReUse,Poveda:15}. In particular, the dynamics considered in this paper have set-valued jump maps that lead to non-unique solutions with non-trivial hybrid time domains having multiple simultaneous jumps in the standard continuous time domain, a behavior that is unavoidable in decentralized multi-agent HDS. We also show that these adaptive dynamics can approximately recover the acceleration properties of the model-based algorithms.

To our best knowledge, the algorithms of this paper (model-based, partial-information, and model-free) are the first in the literature that implement dynamic momentum and distributed restarting techniques in $n$-player noncooperative games.

\tcb{The rest of this paper is organized as follows. Section II presents preliminary concepts of game theory and hybrid dynamical systems. Section III presents the problem statement. Section IV presents a distributed coordination restarting mechanism, and uses it to formulate momentum-based hybrid NSS dynamics for games where players have access to full-information Oracles. Section V relaxes this assumption and uses multi-time scale techniques to study the case in which players only have access to local information from their neighbors. Section VI presents model-free NSS dynamics to tackle setups in which players can only measure the value of their local cost functions. Finally, Section VII presents the analysis and proofs.}
\vspace{-0.1cm}
\section{PRELIMINARIES}
\label{sec_prelim}
\subsubsection{Notation} Given a compact set $\mathcal{A}\subset\R^n$ and a vector $z\in\R^n$, we use $|z|_{\mathcal{A}}\coloneqq \min_{s\in\mathcal{A}}\|z-s\|_2$ to denote the minimum distance of $z$ to $\mathcal{A}$. \tcb{We use $\mathbf{1}_n$ to represent an $n$-dimensional vector with $1$ in all its entries, and define $\mathbf{1}_n\cdot A \coloneqq \set{x\in \mathbb{R}^n~:~ x_1=x_2=\hdots=x_n=a,~a\in A}$, for any set $A\subset \mathbb{R}$.} We use $\mathbb{S}^1\coloneqq \{z\in\R^2:z^2_1+z_2^2=1\}$ to denote the unit circle in $\R^2$, and $\mathbb{T}^n=\mathbb{S}^1\times\cdots\times\mathbb{S}^1$ to denote the $n^{th}$ Cartesian product of $\mathbb{S}^1$. We also use $r\mathbb{B}$ to denote a closed ball in the Euclidean space, of radius $r>0$, and centered at the origin. We use $I_n\in\R^{n\times n}$ to denote the identity matrix, and $(x,y)$ for the concatenation of the vectors $x$ and $y$, i.e., $(x,y)\coloneqq [x^\top,y^\top]^\top$. Also, we use $\mathcal{D}(k)$ to represent a diagonal matrix of appropriate dimension with diagonal given by the entries of a vector $k$. We also use $\overline{k}$ (resp. $\underline{k}$) to denote the largest (resp. smallest) entry of $k$. A function $\beta:\R_{\geq0}\times\R_{\geq0}\to\R_{\geq0}$ is said to be of class $\mathcal{K}\mathcal{L}$ if it is non-decreasing in its first argument, non-increasing in its second argument, $\lim_{r\to0^+}\beta(r,s)=0$ for each $s\in\R_{\geq0}$, and  $\lim_{s\to\infty}\beta(r,s)=0$ for each $r\in\R_{\geq0}$.

\vspace{0.1cm}
\subsubsection{Games} In this paper, we consider noncooperative games with $n\in\mathbb{Z}_{\ge 2}$ players, where each player $i$ can control its own action $q_i$, and has access to the actions $q_j$ of neighboring players $j\in\mathcal{N}_i\coloneqq \{j\in\mathcal{V}:(i,j)\in\mathcal{E}\}$, who are characterized by an undirected, connected, and time-invariant graph $\mathbb{G}=\{\mathcal{E},\mathcal{V}\}$, where $\mathcal{V}=\{1,2,\ldots,n\}$ is the set of players and $\mathcal{E}$ is the set of edges between players. We use $\mathcal{L}$ to denote the Laplacian matrix of the graph $\mathbb{G}$. The main goal of each player $i$ is to minimize its own cost function $\phi_i:\R^{n}\to\R$ by controlling its own action $q_i$. We assume that the costs $\phi_i$ are \emph{twice continuously differentiable}, and we use $q=(q_1,q_2,\ldots,q_n)$ to denote the overall vector of actions of the game. We also use denote as $q_{-i}$ the vector of all actions with the action of player $i$ removed. To simplify our exposition, we assume that the actions $q_i$ are scalars.  However, all our results also hold for vectorial actions by using suitable Kronecker products. We use $\mathcal{G}$ to denote the pseudogradient of the game, where $q\mapsto \mathcal{G}(q)\coloneqq \left(\frac{\partial \phi_1(q)}{\partial q_{1}},\frac{\partial \phi_2(q)}{\partial q_{2}},\ldots,\frac{\partial \phi_n(q)}{\partial q_{n}}\right)\in\R^n$. Following standard assumptions in the literature of fast NES \cite{Grammatico21,PersisGrammatico} and accelerated optimization \cite{ODE_Nesterov,Wilson18,HamiltonianDescent,HighResolution2018,Muehlebach1,zero_order_poveda_Lina,OchoaPoveda20LCSS}, in this paper we will work with the following assumptions.
\vspace{0.05cm}
\begin{assumption}\label{lipschitzassumption}
The mapping $\mathcal{G}$ is $\ell$-globally Lipschitz, i.e., there exists a constant $\ell>0$ such that $|\mathcal{G}(q)-\mathcal{G}(q')|\leq \ell|q-q'|$, for all $q,q'\in\R^n$. \QEDB
\end{assumption}
\vspace{0.05cm}
%
\begin{assumption}\label{cocoerciveassumption}
The mapping $\mathcal{G}$ is $1/\ell$-cocoercive, i.e., there exists $\ell$ such that $\big(\mathcal{G}(q)-\mathcal{G}(q')\big)^\top (q-q')\geq \frac{1}{\ell}|\mathcal{G}(q)-\mathcal{G}(q')|^2$ for all $q,q'\in \R^n$. Moreover, the map $q\mapsto |\mathcal{G}(q)|^2$ is radially unbounded. \QEDB
\end{assumption}

\vspace{0.05cm}
The first property of Assumption \ref{cocoerciveassumption} implies Assumption \ref{lipschitzassumption} via direct application of the Cauchy–Schwarz inequality. However, the converse is not necessarily true in non-potential games \cite{Vandenbeerghe19}. We will also use the following definition to characterize the monotonicity properties of the games.
\vspace{0.07cm}

\begin{definition}\label{main_assumptions}
A game with pseudogradient $\mathcal{G}$ is said to be:
\begin{enumerate}
    \item \emph{monotone}  if it satisfies $\big(\mathcal{G}(q)-\mathcal{G}(q')\big)^\top (q-q')\ge 0$, for all $q,~ q'\in\R^n$.
    \item  \emph{strictly monotone} if it satisfies $\big(\mathcal{G}(q)-\mathcal{G}(q')\big)^\top (q-q')>0$, for all $q\neq q'\in\R^n$.
    \item \emph{$\kappa$-strongly monotone} with $\kappa>0$, if it satisfies $\big(\mathcal{G}(q)-\mathcal{G}(q')\big)^\top (q-q')\geq \kappa|q-q'|^2$, for all $q,q'\in\R^n$.
    \item \emph{$\kappa$-strongly monotone quadratic} if it is a \tcb{$\kappa$-strongly monotone game} with $\mathcal{G}(q)=Aq+b$ for some $A\in\mathbb{R}^{n\times n}$ and $b\in\mathbb{R}^n$. 
    \item \emph{\tcb{potential}} if there exists a continuously differentiable and radially unbounded function $P:\R^n\to\R$, such that \tcb{$\mathcal{G}(q)=\nabla P(q)$}, for all $q\in\R^n$. \QEDB
\end{enumerate}
\vspace{0.1cm}
Monotone, strictly and $\kappa$-strongly monotone games that are also potential games will be referred to as \emph{\tcb{monotone potential games}}, \emph{\tcb{strongly monotone potential games}} and \emph{\tcb{$\kappa$-strongly monotone potential games}}, respectively. Their monotonicity properties are \tcb{defined by $\mathcal{G}$}. \QEDB
\end{definition}
%
%
%
%
%
%
%
%
\begin{remark}
Cocoercive mappings are monotone but not necessarily strongly monotone. However, games that are \tcb{$\kappa$-strongly monotone} and $\ell$-Lipschitz are also $\kappa/\ell^2$-cocoercive \cite[Prop. 2.1]{Vandenbeerghe19}. \QEDB
\end{remark}%
Strict monotonicity of  $\mathcal{G}$ implies that there is exactly one Nash equilibirium, if it exists. For \tcb{$\kappa$-strongly monotone game}s and \tcb{monotone potential games}, existence is always guaranteed \cite[Thm. 2.3.3]{facchinei2007finite}. In some cases, we will also work with the following assumption.
\begin{assumption}\label{regularassumption}
The function $\phi_i:\R\to \R$ is radially unbounded in $q_i$ for ever $q_{-i}\in \R^{n-1}$ and all $i\in \mathcal{V}$. \QEDB\end{assumption}

\vspace{0.1cm}
\subsubsection{Hybrid Dynamical Systems} To study our algorithms, in this paper we consider HDS with state $x\in\mathbb{R}^n$, and dynamics

\vspace{-0.5cm}
\begin{align}\label{HDS}
x\in C,~\dot{x}=F(x),~~~\text{and}~~~x\in D,~~x^+\in G(x),
\end{align}
where $x\in\R^n$ is the state of the system, $F:\R^n\to\R^n$ is called the flow map, $G:\R^n\rightrightarrows\R^n$ is a set-valued mapping called the jump map, and $C\subset\R^n$ and $D\subset\R^n$ are closed sets, called the flow set and the jump set, respectively \cite{HDS}. We use $\mathcal{H}=(C,F,D,G)$ to denote the data of the HDS $\mathcal{H}$. Solutions $x:\text{dom}(x)\to\R^n$ to system \eqref{HDS} are indexed by a continuous time parameter $t$, which increases continuously during flows, and a discrete-time index $j$, which increases by one during jumps. Therefore, solutions $x:\text{dom}(x)\to\mathbb{R}^n$ to system \eqref{HDS} are defined on \emph{hybrid time domains}. Solutions that have an unbounded time domain are said to be \emph{complete}. For a precise definition of hybrid time domains and solutions to HDS of the form \eqref{HDS}, we refer the reader to \cite[Ch.2]{HDS}. The following definitions will be instrumental to study the stability and convergence properties of systems of the form \eqref{HDS}.

\vspace{0.1cm}
\begin{definition}
The compact set $\mathcal{A}\subset C\cup D$ is said to be \emph{uniformly globally asymptotically stable} (UGAS) for system \eqref{HDS} if $\exists$ $\beta\in\mathcal{K}\mathcal{L}$ such that every solution $x$ satisfies:  
\begin{equation}\label{KLbound}
|x(t,j)|_{\mathcal{A}}\leq \beta(|x(0,0)|_{\mathcal{A}},t+j),~ \forall~(t,j)\in\text{dom}(x). 
\end{equation}
When $\beta(r,s)=c_1re^{-c_2s}$ for some $c_1,c_2>0$, the set $\mathcal{A}$ is \emph{uniformly globally exponentially stable} (UGES). When $\exists$ $T^*>0$ such that $\beta(r,s)=0$, $\forall$ $s\geq T^*,r>0$, the set $\mathcal{A}$ is said to be \emph{uniformly globally fixed-time stable} (UGFxS). \QEDB
\end{definition}

\vspace{0.1cm}
We will also consider $\varepsilon$-parameterized HDS of the form:
\begin{align}\label{HDS2}
x\in C_{\varepsilon},~~\dot{x}=F_{\varepsilon}(x),~~~\text{and}~~~x\in D_{\varepsilon},~x^+\in G_{\varepsilon}(x),
\end{align}
where $\varepsilon>0$. For these systems we will study \emph{semi-global practical stability} properties as $\varepsilon\to0^+$.

\vspace{0.1cm}
\begin{definition}\label{definitionSGPAS}
The compact set $\mathcal{A}\subset C\cup D$ is said to be \emph{Semi-Globally Practically Asymptotically Stable} (SGP-AS) as $\varepsilon\to0^+$ for system \eqref{HDS2} if $\exists$ $\beta\in\mathcal{K}\mathcal{L}$ such that for each pair $\delta{>}\nu{>}0$ there exists $\varepsilon^*>0$ such that for all $\varepsilon\in(0,\varepsilon^*)$ every solution of \eqref{HDS2} with $|x(0,0)|_{\mathcal{A}}\leq \delta$ satisfies  
\begin{equation}\label{SGPASbound}
|x(t,j)|_{\mathcal{A}}\leq \beta(|x(0,0)|_{\mathcal{A}},t+j)+\nu,
\end{equation}
$\forall~(t,j)\in\text{dom}(x)$. When the function $\beta$ has an exponential form, we say that $\mathcal{A}$ is \emph{semi-globally practically exponentially stable} (SGP-ES). ~\hfill\QEDB
\end{definition}

\vspace{0.1cm}
The notions of SGP-AS (-ES) can be extended to systems that depend on multiple parameters $\varepsilon=(\varepsilon_1,\varepsilon_2,\ldots,\varepsilon_{\ell})$. In this case, we say that $\mathcal{A}$ is SGP-AS as $(\varepsilon_{\ell},\ldots,\varepsilon_2,\varepsilon_{1})\to0^+$ where the parameters are tuned in order starting from $\varepsilon_1$. 
%

\vspace{0.05cm}
\begin{definition}
Consider the perturbed HDS
\begin{subequations}\label{perturbedHDS}
\begin{align}
&x+e\in C,~~~~~\dot{x}=F(x+e)+e,\\
&x+e\in D,~~~x^+\in G(x+e)+e,
\end{align}
\end{subequations}
where $\text{dom}(e)=\text{dom}(x)$, and $e$ is a measurable function satisfying $\sup_{(t,j)\in\text{dom}(e)}|e(t,j)|\leq\varepsilon$ with $\varepsilon>0$. System \eqref{perturbedHDS} is said to be R-UGAS (resp. R-UGES) if: 1) it is UGAS (resp. ES) when $\varepsilon=0$; and 2) it is SGP-AS (resp. SGP-ES) as $\varepsilon\to0^+$. \QEDB
\end{definition}
\vspace{0.05cm}

Finally, the following definition will be instrumental for the analysis of some of our algorithms.

\vspace{0.05cm}
\begin{definition}\label{closeness}
Two hybrid signals $x_1:\text{dom}(x_1)\to\mathbb{R}^n$ and $x_2:\text{dom}(x_2)\to\mathbb{R}^n$ are said to be $(T,J,\varepsilon)$-close if: (1) for each $(t,j)\in\text{dom}(x_1)$ with $t\leq T$ and $j\leq J$ there exists $s$ such that $(s,j)\in\text{dom}(x_2)$, with $|t-s|\leq \varepsilon$ and $|x_1(t,j)-x_2(t,j)|\leq \varepsilon$; (2) for each $(t,j)\in\text{dom}(x_2)$ with $t\leq T$ and $j\leq J$ there exists $s$ such that $(s,j)\in\text{dom}(x_1)$, with $|t-s|\leq\varepsilon$ and $|x_2(t,j)-x_1(t,j)|\leq \varepsilon$. \QEDB
\end{definition}
%

\vspace{-0.1cm}
\section{PROBLEM STATEMENT AND MOTIVATION}
\label{probl_statement}
%
%
A NE is defined as an action profile $q^*\in\R^{n}$ that satisfies   
\begin{equation}\label{NEcondition}
\phi_i(q_i^*,q^*_{-i})= \inf_{q_i\in\R}\phi_i(q_i,q_{-i}^*),~~\forall~i\in\mathcal{V}.
\end{equation}
Particularly, given a \tcb{monotone game} with pseudogradient $\mathcal{G}$, it follows that $q^*$ is a NE if and only if $\mathcal{G}(q^*)=0$ \cite[Prop. 2.1]{mertikopoulos2019learning}. 
%
%
%
%
%
%
Our goal is to \emph{efficiently} and \emph{robustly} find the set of points $q^*$ that satisfy \eqref{NEcondition}, denoted $\mathcal{A}_{\text{NE}}$, using algorithms with dynamic momentum. However, as the following example shows, this task is not trivial, even for potential games. 

\begin{figure}[t!]
    \centering
    \includegraphics[width=0.6\linewidth]{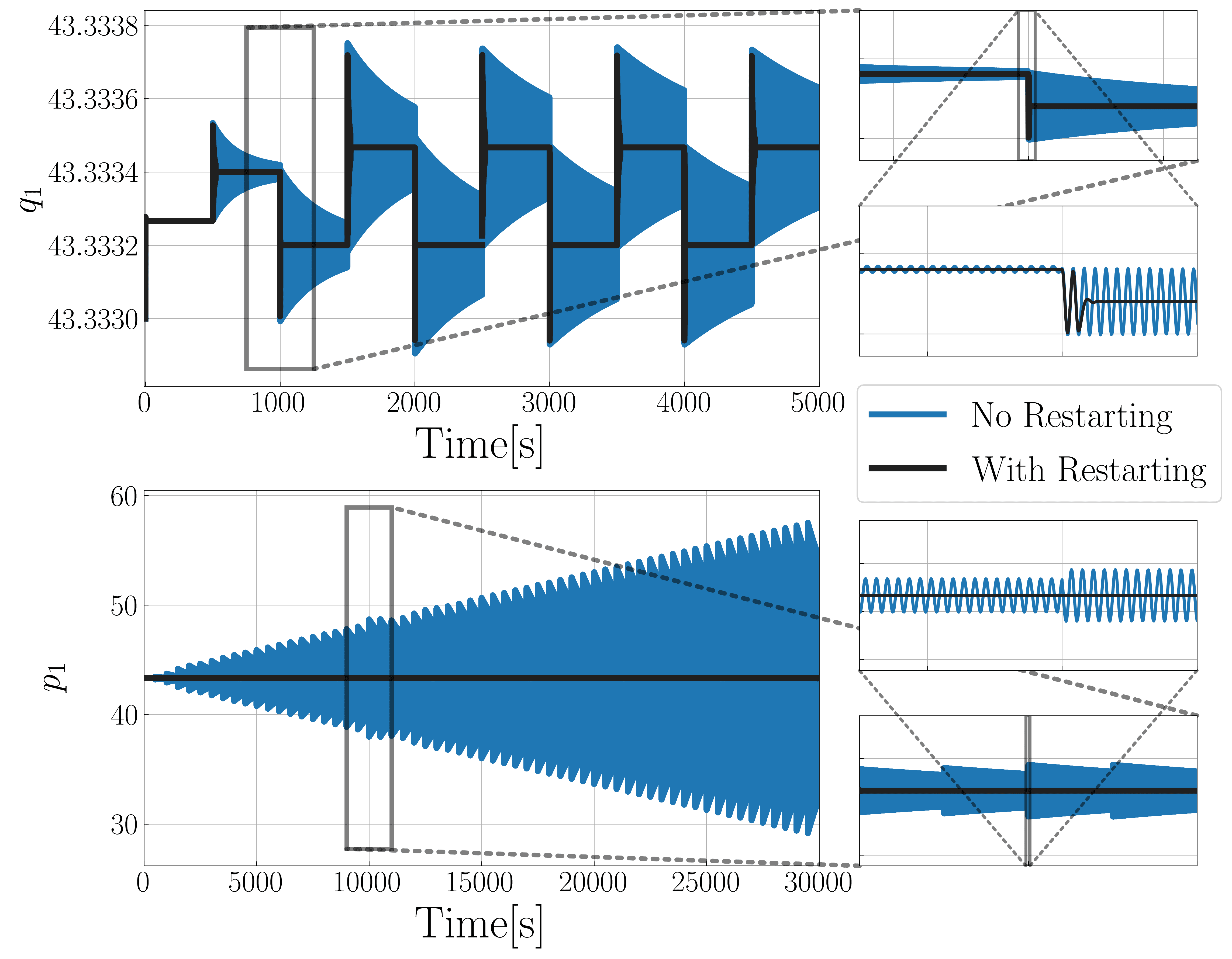}
    \caption{Instability of \eqref{ODEmomentum} in a duopoly game with perturbed gradients and $T_0= 2\sqrt{2}\times10^{-3}$. The instability can be removed by incorporating resets, which generate the stable trajectories shown in black.}
    \label{fig:example1}
    \vspace{-0.3cm}
\end{figure}

\vspace{0.1cm}
\begin{example}(\textsl{Instability Under Small Disturbances})\label{negative_example1}
Consider a duopoly game with pseudogradient $\mathcal{G}(q)=Aq+b$, where $A=[10,-5;-5,10]$, and $b=[-250,-150]$. This is a \tcb{$\kappa$-strongly monotone potential-game} studied in \cite[Sec. II]{Frihauf12a} using PSG flows. The unique NE is $q^*=[130/3,101/3]^\top$, and since $A$ is symmetric, the game has a (quadratic) potential function, which permits the direct application of \cite[Thm. 3]{ODE_Nesterov} to conclude convergence of all functions $q$ generated by \eqref{ODEmomentum} towards the NE $q^*$.
Nevertheless, if \eqref{ODEmomentum} is implemented with a perturbed gradient $\mathcal{G}(q)+e(t)$, where $t\mapsto e(t)$ is an arbitrarily small periodic disturbance satisfying $|e(t)|\leq \delta$ for all $t\geq0$, the highly oscillatory unstable behaviour shown in blue in Figure \ref{fig:example1} emerges. Indeed, by \cite[Thm.1]{PovedaTeelACC20}, in this case there is no $\beta\in \mathcal{K}\mathcal{L}$ such that the bound \eqref{SGPASbound} holds for the trajectories of system \eqref{ODEmomentum}. However, we will show that this bound actually exists when resets are used to restart $(\tau,\dot{q})$, generating the stable behavior shown in black color in Figure \ref{fig:example1}. \QEDB



\end{example}

\vspace{0.1cm}
The robustness issues illustrated in Example \ref{negative_example1} prevent the direct implementation of the momentum-based dynamics \eqref{ODEmomentum} in noisy environments, or in settings where some of the states or gradients are computed on-the-fly using multi-time scale techniques, such as in singular perturbation and averaging theory. In fact, such techniques usually require ``reduced'' or ``average'' systems with stability properties characterized by $\mathcal{K}\mathcal{L}$ bounds; see \cite[Ch. 11]{khalil} and \cite[Assumption 4.]{WangTeelNesic}. 

The incorporation of resets can resolve the robustness issues of system \eqref{ODEmomentum}. However, as the following example shows, a naive distributed and \emph{uncoordinated} implementation of resets in network games can hinder the potential advantages of using algorithms with dynamic momentum.

\begin{figure}[t!]
    \centering
    \includegraphics[width=0.6\linewidth]{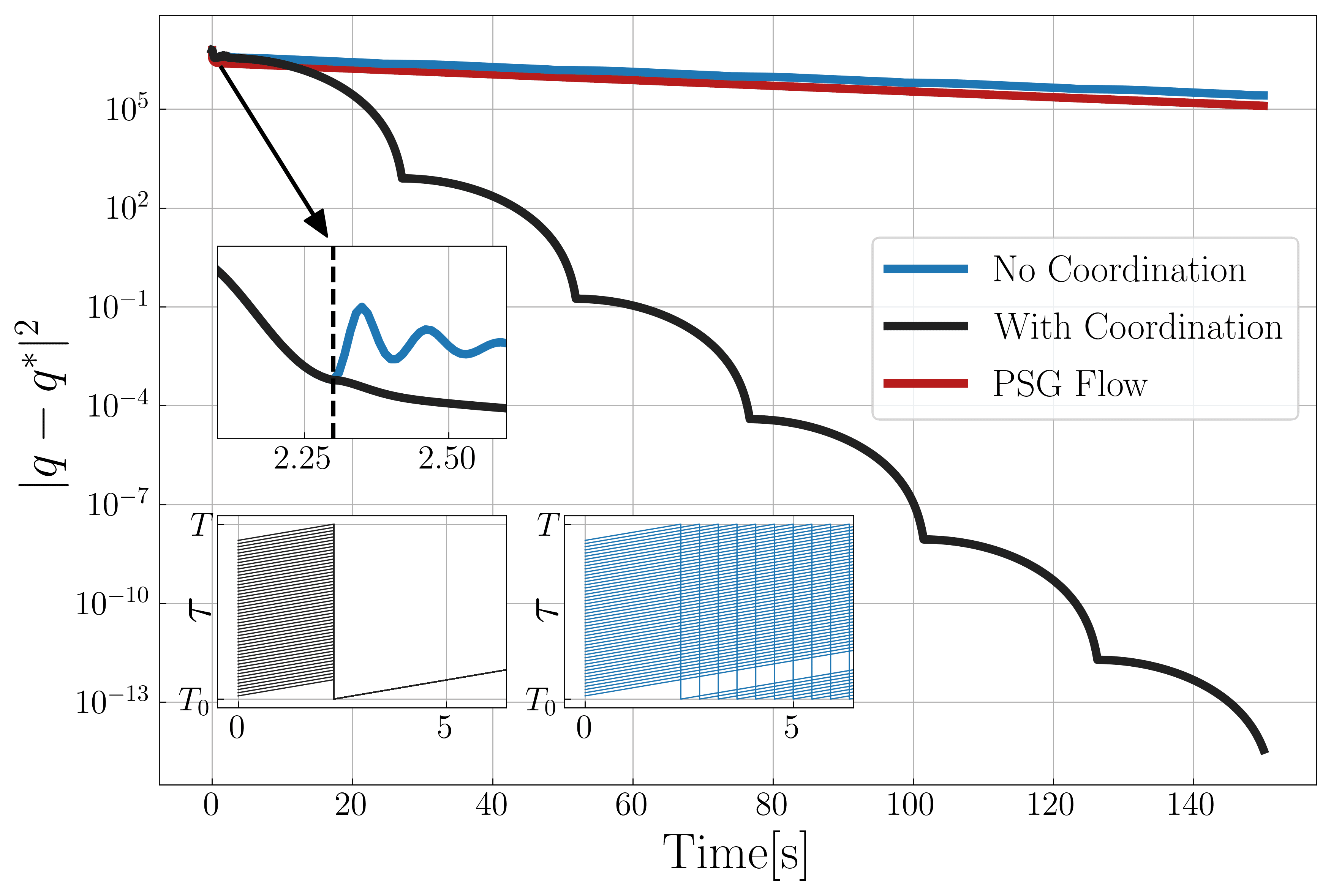}
    \caption{Coordinated vs non-coordinated resets in a quadratic \tcb{$\kappa$-strongly monotone potential-game} with $\kappa = 0.01, \ell=100$ and $n=30$. The insets show the evolution of the states $\tau_i$ with and without coordination mechanisms.}
    \label{fig:example2}
    \vspace{-0.5cm}
\end{figure}

\vspace{0.06cm}
\begin{example}\textsl{(Slow Convergence and Uncoordinated Resets)}\label{negative_example2}
Consider a \emph{distributed} implementation of system \eqref{ODEmomentum} in a \tcb{$\kappa$-strongly monotone potential-game} with 30 players and  $\kappa=0.01$. Each player $i$ implements its own states $(q_i,p_i,\tau_i)$, with dynamics $\dot{q}_i=\frac{2}{\tau_i}(p_i-q_i),~\dot{p}_i=-2\tau_i\frac{\partial\phi_i}{\partial q_i}$, and $\dot{\tau}_i=\eta$, with $\eta=\frac{1}{2}$. Also, players implement periodic resets of $(\tau_i,p_i)$ every 25 seconds (in their own local time reference frame) via the individual jump maps $\tau_i^+=0.1$ and $p_i^+=q_i$. While this periodic reset strategy has been shown to guarantee fast convergence in \emph{centralized} optimization problems, e.g., \cite[Thm. 1]{Poveda_Li:2019_CDC}, 
%
%
Figure \ref{fig:example2} shows the emerging behavior in noncooperative games with distributed and uncoordinated resets. As shown in blue, the solutions of \eqref{ODEmomentum} actually converge to the NE, but at a slower rate compared to the standard PSG flow. We also show a black trajectory corresponding to players implementing \emph{coordinated} resets. Here, the acceleration properties of the momentum-based dynamics can be fully exploited; c.f. Theorem \ref{theorem1}. \QEDB
\end{example}

\vspace{0.05cm}
%
%

%
%
%
Examples \ref{negative_example1} and \ref{negative_example2} were developed for potential games. Next, we show that in non-potential games the solutions of \eqref{ODEmomentum} might not even converge to the NE of the game, even when resets are implemented.
%

\begin{figure}[t!]
    \centering
    \includegraphics[width=0.6\linewidth]{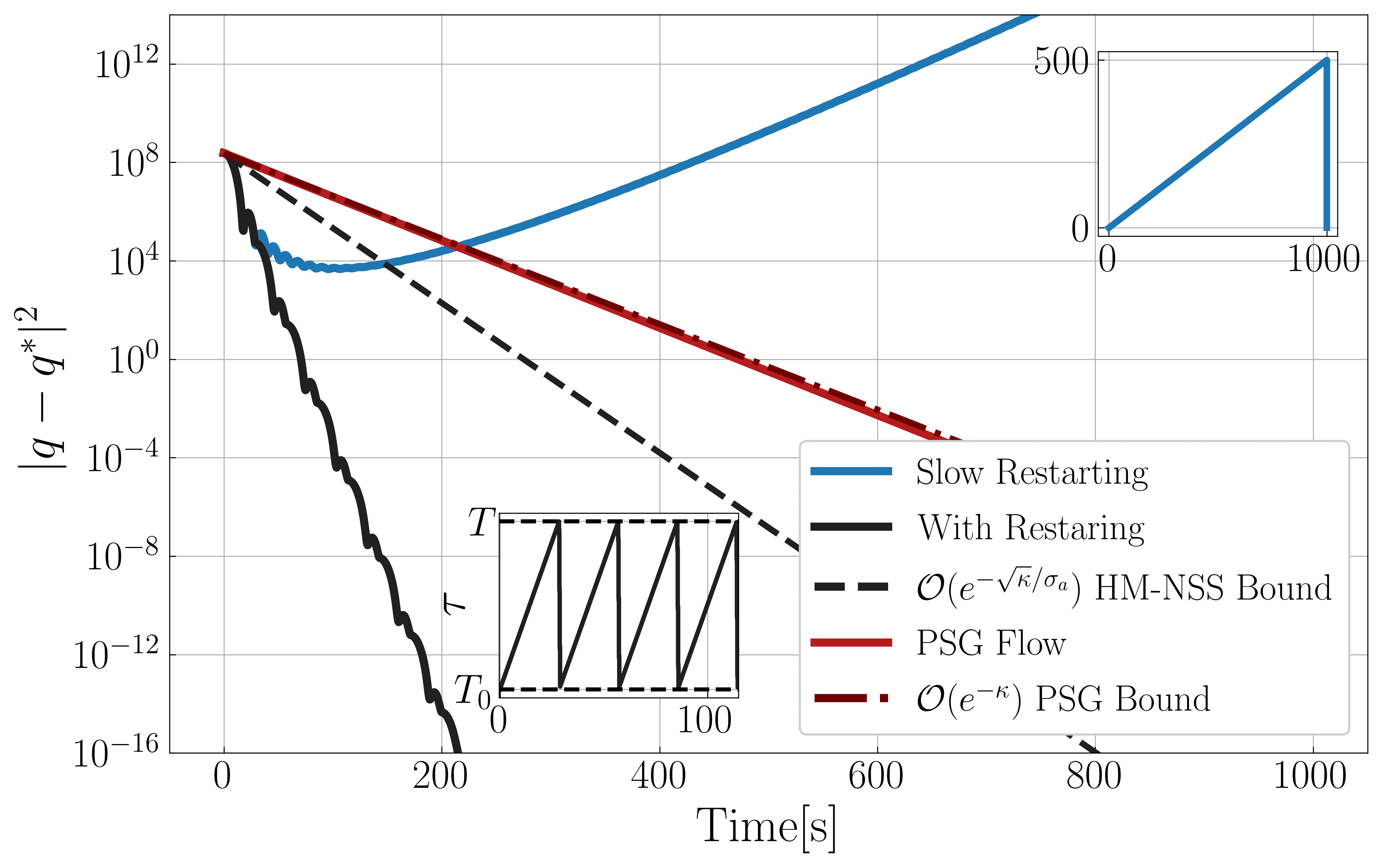}
    \caption{Lack of convergence of trajectories of \eqref{ODEmomentum} in a non-potential \tcb{$\kappa$-strongly monotone quadratic game} with $\kappa = 0.02,~\ell=0.0214,~n=30,~T_0=0.1,~T = 3.74$. The black line shows athe trajectory of the proposed hybrid controller.with resets.}
    \label{fig:example3}
    \vspace{-0.5cm}
\end{figure}

\vspace{0.05cm}
\begin{example}(\emph{Lack of Convergence in Non-Potential Games})\label{example3} We consider a non-potential \tcb{$\kappa$-strongly monotone quadratic game} with $30$ players and $\kappa=0.02$. For this game, the standard PSG flow guarantees exponential convergence via \cite[Thm. 1]{rosen}. However, as shown in color blue in Figure \ref{fig:example3}, system \eqref{ODEmomentum} generates trajectories that diverge, even when resets are slowly implemented. The same plot shows in black color a trajectory that rapidly converges to the unique NE of the game. We will show that this stable and fast behavior can be guaranteed using a \emph{hybrid} algorithm with distributed coordinated resets that dissipate energy ``sufficiently often'' via suitable contraction properties; c.f., Theorem \ref{theorem:strong:np}. \QEDB
\end{example}
%

%
%
%
\vspace{0cm}

\section{DISTRIBUTED HYBRID NSS DYNAMICS WITH COORDINATED RESTARTING}
\label{Sec_R1}
To achieve robust NSS with dynamic momentum, we start by endowing each player $i\in\mathcal{V}$ with a state $x_i=(q_{i},p_{i},\tau_i)\in\R\times\R\times\R_{>0}$, and a gradient Oracle that provides real-time measurements of the partial derivative $\frac{\partial \phi_i(q)}{\partial q_i}$ at the overall action state $q\in\mathbb{R}^n$. The reset mechanisms of the players make use of three positive tunable parameters $(\eta,T_0,T)$, which satisfy $T>T_0>0$ and $1/2\ge\eta>0$, and which are selected \emph{a priori} by the system designer. The state $x_i$ evolves according to hybrid dynamics that are coordinated by a local timer $\tau_i$. In particular, the continuous-time dynamics of each player are
\begin{align}\label{decoupled_flows111}
\tau_i\in[T_0,T){\implies} \left(\begin{array}{c} 
\dot{q}_{i}\\
\dot{p}_{i}\\
\dot{\tau_i}
\end{array}
\right){=}F_i(x):=\left(\begin{array}{c} 
\frac{2}{\tau_i}(p_{i}-q_{i})\\
-2\tau_i\frac{\partial \phi_i(q)}{\partial q_i}\\
\eta
\end{array}
\right),
\end{align}
and the discrete-time dynamics are given by
\begin{equation}\label{decoupled_flows2}
\tau_i=T{\implies} \left(\begin{array}{c} 
q^+_{i}\\
p^+_{i}\\
\tau^+_i
\end{array}
\right){=}R_i(x_i)\coloneqq \left(\begin{array}{c} 
q_{i}\\
\alpha_i p_{i}+(1-\alpha_i)q_{i}\\
T_0
\end{array}
\right).
\end{equation}
%
In \eqref{decoupled_flows2}, the parameters 
$\alpha_i\in\{0,1\}$ model the different individual \emph{reset policies} of the players. The choice $\alpha_i=0$ leads to resets of the form $p_{i}^+=q_{i}$, which corresponds to $\dot{q}^+_i=0$, i.e., the momentum of player $i$ is reset. On the other hand, $\alpha_i=1$ corresponds to keeping $p_{i}$ constant.  
%
%

Since players have access to Oracles that provide real-time evaluations of their gradient, they can implement the hybrid dynamics \eqref{decoupled_flows111}-\eqref{decoupled_flows2} in a fully decentralized way by running their own timers $\tau_i$ to coordinate the flows \eqref{decoupled_flows111} and the jumps \eqref{decoupled_flows2}. However, as shown in Example \ref{negative_example2}, lack of coordination between the resets of the players can hinder the acceleration properties expected from using dynamic momentum, even in potential games when all players implement the same reset policy $\alpha_i$, c.f., Figure \ref{fig:example2}. To address this issue, we proceed to endow each player with a distributed hybrid coordination mechanism for the resets. 
%

\vspace{-0.3cm}
\subsection{Coordinated Distributed Resets}
\label{sec_coordination_mechanism}
The coordination mechanism of each player $j\in\mathcal{V}$ uses a set-valued coordination mapping $\mathcal{C}_j:\R_{\geq0}\rightrightarrows\R_{\geq0}$, defined as
\begin{equation}\label{coordination_mechanism}
\mathcal{C}_j(\tau_j)\coloneqq \left\{\begin{array}{cl}
T&\text{if}~\tau_j\in(T_0+r_j,T]\\
\{T_0,T\}&\text{if}~\tau_j=T_0+r_j\\
T_0&\text{if}~\tau_j\in[T_0,T_0+r_j)
\end{array}\right.,
\end{equation}
where the individual parameter $r_j>0$ satisfies $r_j\in\left(0,\frac{T-T_0}{n}\right)$. Using $\mathcal{C}_j$, the coordination mechanism works as follows: whenever the timer of player $i$ satisfies $\tau_i=T$, the following two events happen:
\begin{enumerate}
    \item Player $i$ resets its own state $x_i$ using the dynamics \eqref{decoupled_flows2}, and
    \item Player $i$ sends a pulse to its neighbors $j\in \mathcal{N}_i$, who proceed to update their state $x_j$ as follows:
    \begin{equation}\label{update1}
    q_{j}^+=q_{j},~~~~~p_{j}^+=p_{j},~~~~~\tau_j^+\in \mathcal{C}_j(\tau_j).
    \end{equation}
\end{enumerate}
Since player $i$ can only signal its neighbors, the rest of the players $j\notin\mathcal{N}_i$ will keep their states constant after the above two events, i.e., $x_{j}^+=x_{j},$ for all $j\notin\mathcal{N}_i$. 

The combination of continuous-time dynamics with momentum \eqref{decoupled_flows111}, and the set-valued discrete-time dynamics that model the coordinated resets leads to a HDS of the form \eqref{HDS}, where multiple resets can happen simultaneously (in the continuous-time domain) when more than two players satisfy the condition $\tau_i=T$. To ensure that this system has suitable robustness properties we need to guarantee that small disturbances in the states, including in $\tau_i$, do not lead to drastic changes in the behavior of the players. This property can be asserted by working with \emph{well-posed} HDS in the sense of \cite[Ch. 7]{HDS}. Roughly speaking, for a HDS to be well-posed, a suitable (graphically) convergent sequence of solutions of the overall system must also converge (in a graphical sense) to another solution of the hybrid system. In the context of \eqref{decoupled_flows111}-\eqref{update1}, we need to guarantee, among others, that for each $\tau_0\in[T_0,T]$, and each graphically convergent sequence of solutions $\{\tau_k\}_{k\in\mathbb{N}}$ with individual components $\tau_{i,k}$ satisfying
\begin{equation}\label{limitcondition1}
0\leq \tau_{1,k}(0,0)\leq\hdots\leq \tau_{n,k}(0,0)<\tau_0,~~\forall~k\in\mathbb{N},
\end{equation}  
and
\begin{equation}\label{limitcondition2}
\lim_{k\to\infty}\tau_{1,k}(0,0)=\hdots=\lim_{k\to\infty}\tau_{n,k}(0,0)=\tau_0,
\end{equation}
the sequence $\{\tau_k\}_{k\in\mathbb{N}}$ must converge (graphically) to a function $\tilde{\tau}$ that is also a solution starting from the initial condition $\tilde{\tau}_1(0,0)=\tilde{\tau}_2(0,0)=\hdots=\tilde{\tau}_n(0,0)=\tau_0$. Thus, when $\tau_0=T$, the above conditions imply that players will reset their timers $\tau_{i,k}$ sequentially with smaller and smaller times between resets as $k\to\infty$. It follows that in the limit, resets must also be sequential with no time between resets. Since the sequence is determined by the initial conditions, a well-posed model of the coordination mechanism must take into account \emph{every} possible order of sequential resets of the timers $\tau_i$. \tcb{In other words, if multiple players simultaneously satisfy the condition $\tau_i=T$, then we need to consider all possible sequential resets induced by such players.} 
As discussed in \cite{Sync_Poveda}, this behavior is unavoidable in well-posed multi-agent HDS with decentralized discrete-time dynamics.

\vspace{-0.2cm}
\subsection{Well-Posed Hybrid NSS Dynamics}
To formalize the above discussion, we proceed to construct a suitable jump map and a jump set that describe the behavior of the overall NSS dynamics. Specifically, we introduce a new set-valued mapping $G^0:\R^{3n}\rightrightarrows\R^{3n}$, which is defined to be non-empty only when $\tau_i=T$ and $\tau_j\in[T_0,T)$ with $j\neq i$, for each $i\in\mathcal{V}$, and has elements given by
\begin{align}\label{ginitialmap}
G^0(x)&\coloneqq \Big\{(v_1,v_2,v_3)\in\R^{3n}:(v_{1,i},v_{2,i},v_{3,i})=R_i(x_i),\notag\\
&~~~~~~v_{1,j}=q_{j},v_{2,j}=p_{j},v_{3,j}\in \mathcal{C}_j(\tau_j),~\forall~j\in\mathcal{N}_i,\notag\\
&~~~~~~~v_j=x_j,\forall~j\notin\mathcal{N}_i\Big\},
\end{align}
where $x\coloneqq (x_1,x_2,\cdots,x_n)$, and where the \emph{reset map} $R_i$ and the \emph{coordination mapping} $\mathcal{C}_j$ are defined in \eqref{decoupled_flows2} and \eqref{coordination_mechanism}, respectively. Using the construction \eqref{ginitialmap}, the jump map of the overall hybrid system is defined as
\begin{equation}\label{jump_map}
x^+\in G_1(x)\coloneqq \overline{G^0}(x),
\end{equation}
where $\overline{G^0}$ is the outer-semicontinuous hull of $G^0$, \cite[pp. 154]{Rockafellar}, i.e., the unique set-valued mapping that satisfies $\text{graph}(G_1)=\text{cl}(\text{graph}(G^0))$. By construction, the mapping $G_1$ is locally bounded and outer-semicontinuous in $\R^n\times\R^n\times[T_0,T]^n$. Moreover, it preserves the sparsity properties of the graph $\mathbb{G}$. It also guarantees that any pair of resets of the form \eqref{decoupled_flows2} happen sequentially, thus satisfying conditions \eqref{limitcondition1}-\eqref{limitcondition2}. 

\begin{figure}
    \centering
    \includegraphics[width=0.55\linewidth]{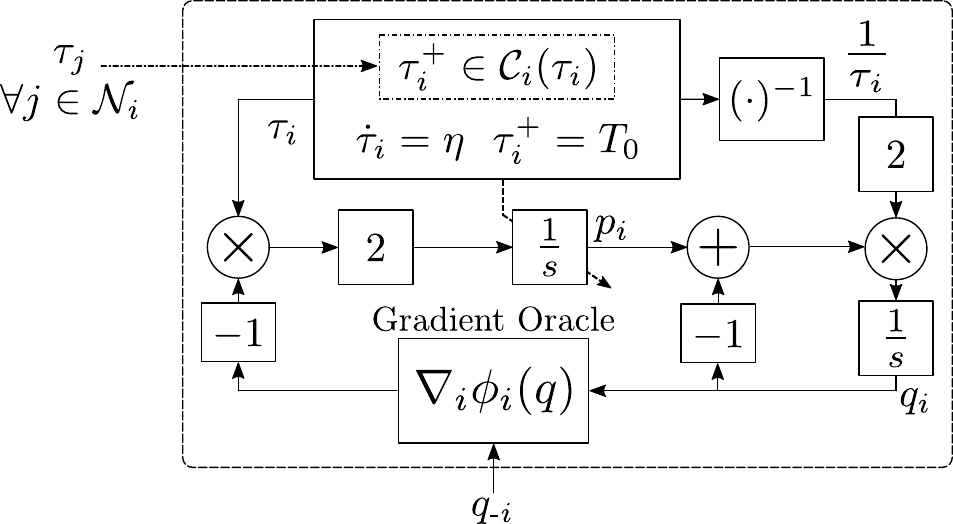}
    \caption{Scheme of Individual HM-NSS dynamics. Periodic coordinated resets restart the state $p_i$ and the timer $\tau_i$.}
    \label{blockdiagram1}
    \vspace{-0.4cm}
\end{figure}

\vspace{0.1cm}
Using the jump map \eqref{jump_map}, we can now define the \emph{hybrid momentum-based-}NSS (HM-NSS) dynamics  $\mathcal{H}_1\coloneqq (C_1,F_1,D_1,G_1)$, with overall state $x=(p,q,\tau)\in\R^{3n}$, and vectorial continuous-time dynamics:
\begin{align}\label{flowmap00}
\left(\begin{array}{c} 
\dot{q}\\
\dot{p}\\
\dot{\tau}
\end{array}
\right)=F_1(x)=\left(\begin{array}{c} 
2\mathcal{D}(\tau)^{-1}(p-q)\\
-2\mathcal{D}(\tau)\mathcal{G}(q)\\
\eta\mathbf{1}_n
\end{array}
\right),
\end{align}
where $p\coloneqq (p_1,p_2,\ldots,p_n)$ and $\tau\coloneqq (\tau_1,\tau_2,\ldots,\tau_n)$. The flow set $C_1$ is defined as:
\begin{align}\label{flow_set0}
C_1\coloneqq &\Big\{x\in\R^{3n}:q\in\R^n,~p\in\R^n,~\tau\in[T_0,T]^n \Big\},
\end{align}
the jump map $G_1$ is given by \eqref{jump_map}, and the jump set is
\begin{equation}\label{jump_set11}
D_1\coloneqq \Big\{x\in\R^{3n}:~x\in C_1,~~\max_{i\in\mathcal{V}}\tau_i=T\Big\}.
\end{equation}
Figure \ref{blockdiagram1} shows a block-diagram representation of the hybrid dynamics of each player.  

The next lemma is fundamental for our results. It asserts absence of finite escape times and Zeno behavior. It also guarantees completeness of solutions and fixed-time synchronization of timers. All proofs are presented in Section \ref{sec_proofs}.

\vspace{0.1cm}
\begin{lemma}\label{lemmawellposed}
The HDS \eqref{jump_map}-\eqref{jump_set11} is well-posed in the sense of \cite[Def. 6.29]{HDS}. Moreover, under Assumption \ref{lipschitzassumption}, every maximal solution of $\mathcal{H}_1$ is complete, and there are at most $n$ jumps in any continuous time interval of length $\frac{1}{\eta}(T-T_0)$. Furthermore, for each solution $x$ and for all $(t,j)\in\mathcal{T}(x)$, we have that $x(t,j)\in \mathcal{A}_{\text{sync}}\coloneqq \left(\{T_0,T\}^n\right)\cup \left(\mathbf{1}_n\cdot[T_0,T]\right)$, where
\begin{equation}\label{usefulsetoftimes}
\mathcal{T}(x)\coloneqq \left\{(t,j)\in\text{dom}(x):t+j\geq T^*\right\},
\end{equation}
and $T^*\coloneqq (T-T_0)/\eta+n$.
\QEDB
\end{lemma}
%
\vspace{0.1cm}

 The qualitative behavior of system $\mathcal{H}_1$ will depend on the choice of parameters $(\eta,T,T_0)$, which characterize the frequency and the minimum and maximum values of the momentum coefficient $\tau$. Different choices of $(\eta,T,T_0)$ will lead to different \emph{reset conditions} (RCs). In turn, as hinted in Example 3, and in contrast to standard optimization \cite{zero_order_poveda_Lina}, different types of games will require different RCs to guarantee convergence to the Nash set. These RCs will be defined in terms of the following \emph{condition numbers} of the game, the reset mechanism, and the graph, respectively:
\begin{equation}\label{condition_numbers}
\sigma_{\phi}\coloneqq \frac{\ell}{\kappa},~~~~\sigma_r\coloneqq \frac{T}{T_0},~~~\sigma_{\mathcal{L}}=\frac{\lambda_{\max}(\mathcal{L})}{\lambda_{2}(\mathcal{L})},
\end{equation}
where $\ell$ is given by Assumptions \ref{lipschitzassumption} or \ref{cocoerciveassumption}, $\kappa$ is given in Definition \ref{main_assumptions}, and $\lambda_{2}(\mathcal{L}),\lambda_{\max}(\mathcal{L})$ are the smallest positive and the largest eigenvalues, respectively, of the Laplacian $\mathcal{L}$.   
 

%
%

\vspace{-0.2cm}
\subsection{Main Stability Results}
\label{main_results_first}
%
We study the stability and convergence properties of the dynamics $\mathcal{H}_1$ with respect to the compact set
\begin{equation}\label{definitionofset}
\mathcal{A}\coloneqq \mathcal{A}_{pq}\times\mathcal{A}_{\text{sync}},
\end{equation}
where $\mathcal{A}_{qp}\coloneqq \{(q,p)\in \mathbb{R}^{2n}:p=q,~q\in\mathcal{A}_{\text{NE}}\}$.
%
%
The first RC that we consider is given by
\begin{equation}\tag{$\text{RC}_1$}
T^2-T_0^2>\frac{\rho_J}{2}\cdot \left(1-\underline{\alpha}\right),
\end{equation}
where $\rho_J\in\mathbb{R}_{\geq 0}$ is a parameter to be determined and $\underline{\alpha}=\min_{i\in\mathcal{V}}\alpha_i$. This condition will regulate how frequently players reset their states. Finally, we also introduce the constant 
\begin{equation}\label{gamma_parameterized}
\gamma(\rho_J)\coloneqq \left(1-\frac{1}{\sigma_r^2}-\frac{\rho_J}{2 T^2}\right).
\end{equation}
where $\sigma_r$ is defined in \eqref{condition_numbers}. This quantity will be instrumental to characterize the rates of convergence towards $\mathcal{A}$. 

%
%
\vspace{0.1cm}
\subsubsection{Results for Potential Games}
Our first result focuses on \tcb{monotone potential games} and \tcb{$\kappa$-strongly  potential game}s. 
%
%
%
\vspace{0.1cm}
\begin{thm}\label{theorem1}
Let $\tcb{\mathcal{G}}$ describe an \tcb{monotone potential game}. Suppose that Assumption \ref{lipschitzassumption} holds, and consider the HDS $\mathcal{H}_1$ under (RC$_1$). Then, the following holds: 
\begin{enumerate}[label={(i$_\arabic*$)}]
\item If $\alpha=\mathbf{1}_n$ and $\rho_J\geq0$ then the set $\mathcal{A}$ is R-UGAS. Moreover, and during flows, for any $i\in\mathcal{V}$ the potential function satisfies the bound
\begin{equation}\label{decrease_potential}
P(q(t,j))-P(\mathcal{A}_{\text{NE}})\leq \frac{c_j}{\tau_i^2(t,j)},~~~\forall~(t,j)\in \mathcal{T}(x),
\end{equation}
where $\{c_j\}_{j=0}^{\infty}\searrow 0^+$ depends on $x(0,0)$.

\vspace{0.1cm}

\item If $\alpha\in\{0,1\}^n$, $\mathcal{G}_\omega$ describes a \tcb{$\kappa$-strongly monotone potential game} and $\rho_J=\kappa^{-1}$, then the set $\mathcal{A}$ is R-UGES, and there exists $\lambda>0$ such that for each compact set $K_0\subset C_1\cup D_1$ there exists $M_0>0$ such that for all solutions $x$ with $x(0,0)\in K_0$, and for all $(t,j)\in\text{dom}(x)$ the following bound holds: 
\begin{align}\label{exp_property}
|q(t,j)-q^*|\leq M_0 e^{-\lambda(t+j)}.
\end{align}

\item  If $\alpha=\mathbf{0}_n$, $\mathcal{G}$ describes a \tcb{$\kappa$-strongly monotone potential game} and $\rho_J=\kappa^{-1}$, then the set $\mathcal{A}$ is R-UGES, and for each compact set $K_0\subset C_1\cup D_1$ there exists $M_0>0$ such that all solutions $x$ with $x(0,0)\in K_0$, and for all $(t,j)\in\dom{x}$ the following bound holds:
\begin{equation*}
|q(t,j)-q^*|\leq\sigma_r\sqrt{\sigma_\phi}\left(1-\gamma(\rho_J)\right)^{\frac{\alpha(j)}{2}}M_0,
\end{equation*}
where $\alpha(j)\coloneqq \max\{0,\lfloor\frac{j-n}{n}\rfloor\}$ and $\gamma(\rho_J)\in(0,1)$. \QEDB
%
%
\end{enumerate}
\end{thm}

\vspace{0.2cm}
%

The results of Theorem \ref{theorem1} establish robust NSS for $\mathcal{H}_1$ in monotone and strongly monotone potential games. Thus, unlike system \eqref{ODEmomentum}, for the hybrid dynamics $\mathcal{H}_1$ there exists a class $\mathcal{K}\mathcal{L}$ function $\beta$ such that a bound of the form \eqref{SGPASbound} holds under small bounded additive disturbances on the dynamics. This effectively rules out the instability observed in Figure \ref{fig:example1}. The bounds of Theorem \ref{theorem1} also establish suitable \tcb{semi-acceleration} properties in both monotone and strongly monotone games. 
Such bounds will eventually hold since the UGAS result also implies that for all times $(t,j)\notin \mathcal{T}(x)$, the trajectories remain (uniformly) bounded, and Lemma \ref{lemmawellposed} guarantees completeness of solutions. Specifically, solutions of $\mathcal{H}_1$ exhibit a ``transient phase'', where the momentum coefficients synchronize to each other, followed by a ``\tcb{semi}-acceleration phase'' where the system behaves as having one global momentum coefficient coordinating the overall network. Figures \ref{fig:example1} and \ref{fig:example2} illustrate the advantages of the hybrid NES dynamics $\mathcal{H}_1$ in potential games compared to the ODE \eqref{ODEmomentum}.
\begin{figure*}[t!]
    \centering
    \includegraphics[width=0.95\linewidth]{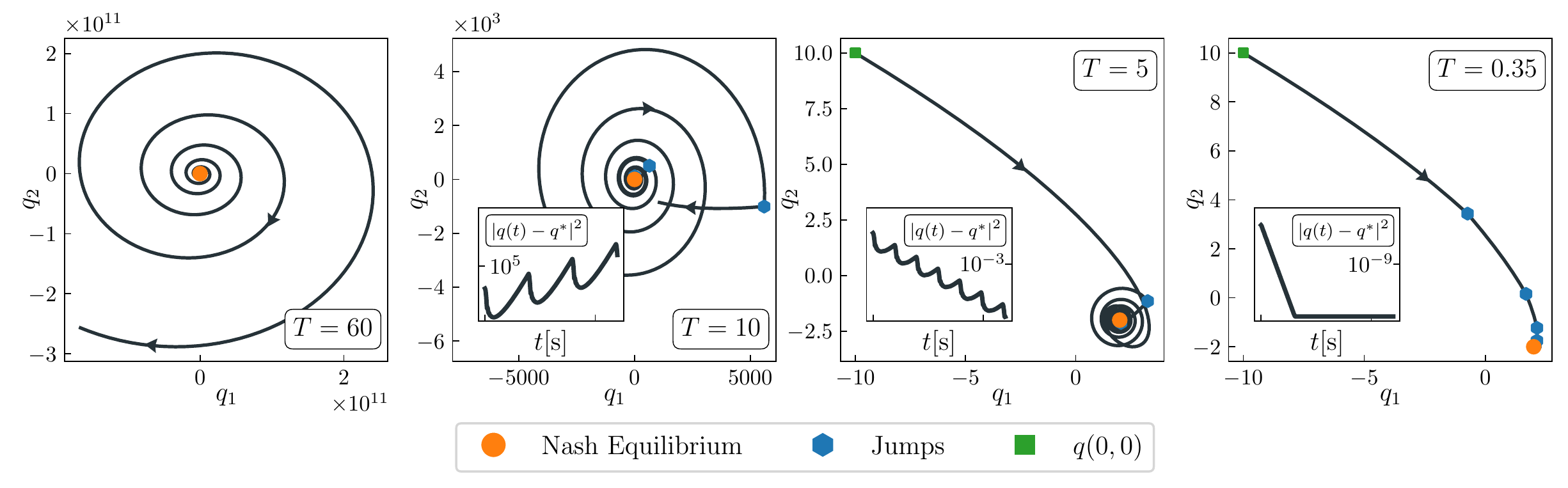}
    \caption{Phase plane plots showing the trajectories of the actions resulting from the HM-NSS dynamics in a non-potential 2-player \emph{\tcb{$\kappa$-strongly monotone quadratic game}} with $\kappa = 6,~\ell=6.2$ and $\tau(0,0)=0.1\cdot\vec{1}_2$.}
    \label{fig:phaseplot}
    \vspace{-0.3cm}
\end{figure*}
\vspace{0.1cm}
\begin{remark}\label{remark1}
 When all players implement the reset protocol $\alpha_i=1$, item (i$_1$) establishes a \emph{semi-acceleration} property of order $\mathcal{O}(1/\tau^2)$ that holds during intervals of flow in $\mathcal{T}(x)$. Since intervals of flow in $\mathcal{T}(x)$ have a length proportional to $T-T_0$, they can be made arbitrarily large by increasing $T$. Moreover, if all players initialize their coefficients as $\tau_i(0,0)=T_0$, then during the first interval of flow we have that $P(q(t,0))-P(q^*)\leq \frac{ d_0}{t^2}$, for all $(t,0)\in\text{dom}(x)$, where $d_0>0$ is fully determined by the initial conditions of the system and the properties $\tcb{\mathcal{G}}$. To the best knowledge of the authors, the result of Theorem \ref{theorem1}-(i$_1$) is the first in the literature that establishes R-UGAS \emph{and} this type of acceleration property in \emph{distributed} NES dynamics. Centralized convergence results without resets were recently studied independently in \cite{GaoCDC20}. \QEDB
\end{remark}

\vspace{0.05cm}
\begin{remark}\label{remark2}
For \emph{\tcb{$\kappa$-strongly monotone potential games}} games, the reset policy $\alpha_i=0,~\forall i\in\mathcal{V}$, guarantees exponential NSS with rate of convergence dictated by $1-\gamma\left(\kappa^{-1}\right)$. In this case, by borrowing results from the literature on centralized accelerated optimization \cite{Candes_Restarting,zero_order_poveda_Lina}, we can consider a ``quasi-optimal'' restarting parameter $T= e\sqrt{\frac{1}{2\kappa}+T_0^2}$, which guarantees exponential convergence of order $\mathcal{O}(e^{-\sqrt{\kappa}t})$ whenever $T_0\ll1$. Finally, the result of item (i$_2$) shows that the stability and convergence properties of $\mathcal{H}_1$ are robust to heterogeneous reset policies in the game.  \QEDB
\end{remark}

\vspace{0.1cm}

\subsubsection{Results for Non-Potential Games} When a potential function does not exist the analysis of the HDS $\mathcal{H}_1$ is more challenging. In this case, we now introduce the following state-dependent matrix parameterized by  $(\rho_F,\delta)\in\mathbb{R}_{>0}\times\mathbb{R}_{\geq0}$:
\begin{equation}\label{contractivematrix}
\mathcal{M}_\delta(q,\rho_F)\coloneqq I_n-\mathcal{S}_{\delta}\left(q,\rho_F\right)\mathcal{S}_{\delta}\left(q,\rho_F\right)^\top,    
\end{equation}
%
%
with $\mathcal{S}_{\delta}:\mathbb{R}^n\times\mathbb{R}_{>0}\to\mathbb{R}^{n\times n}$ given by the scaled matrix
\begin{equation*}
\mathcal{S}_{\delta}(q,\rho_F)\coloneqq \chi(\rho_F,\delta)^\frac{1}{2} \Big(\rho_F I_n-\partial \mathcal{G}(q)\Big),
\end{equation*}
where $\partial \mathcal{G}$ is the Jacobian of $\mathcal{G}$, and where the mapping $\chi:\mathbb{R}_{>0}\times\mathbb{R}_{\geq0}\to\mathbb{R}_{>0}$ is given by
\begin{equation*}
    \chi(\rho_F,\delta)=\frac{T^2}{1-\delta T^2}\cdot\frac{1}{\rho_F(1-\eta) - \delta\rho_F^2},
\end{equation*}
which is defined for all arguments such that $\delta^2T^2<1$ and $1-\eta>\delta\rho_F$. We use the following definition to extend \cite[Def. 4.1.2]{Bernstein_Book} to matrices of the form \eqref{contractivematrix}.

\vspace{0.1cm}
\begin{definition}\label{definitionSC}
The matrix-valued function $q\mapsto\mathcal{S}_\delta(q,\rho_F)$ is said to be 
$\rho_F$-\emph{Globally Contractive} \big($\rho_F$-GC\big) if $\mathcal{M}_\delta(q,\rho_F)\succ0$ for all $q\notin\mathcal{A}_{\text{NE}}$. ~~\hfill\QEDB
\end{definition}

\vspace{0.1cm}
Note that when $\mathcal{M}_\delta\succ0$, the coefficient $\chi$ characterizes the \emph{level of contraction} of $\mathcal{S}_{\delta}$. Indeed, $\mathcal{M}_\delta\succ0$ if and only if
\begin{equation}\label{contractioncondition}
\frac{1}{\chi(\rho_F,\delta)}\geq \sigma_{\max} \Big(\rho_FI_n-\partial \mathcal{G}(q)\Big)^2,    
\end{equation}
where $\sigma_{\max}(\cdot)$ is the maximum singular value of its argument \cite[Thm. 7.7.2]{HornBook}. Using the definition of $\chi$, and inequality \eqref{contractioncondition}, it can be observed that in order to ensure that $\mathcal{S}_\delta$ is $\rho_F$-GC for some pair $(\delta,\rho_F)$, the resetting parameter $T$ cannot be chosen arbitrarily large. Example 4 illustrates this point. 

\vspace{0.1cm}
\begin{example}
Consider a \tcb{$\kappa$-strongly monotone quadratic game} with $\kappa = 6$, and
\begin{equation}\label{pseudogradient:necessityofresets}
    \mathcal{G}(q) = \begin{pmatrix}
            6&1.5\\-1.5&6
    \end{pmatrix}\left(q-q^*\right),
\end{equation}
where $q^* = (2, -2)$. First, let $\delta=0$, and note that for this game $\mathcal{M}_0(q,\rho_F) = \mathcal{D}\left(m_0(\rho_F)\mathbf{1}_{2}\right)\in\R^{2\times 2}$, where
\begin{equation*}
    m_0(\rho_F)= 1-T^2\frac{4 (\rho_F -12) \rho_F +153}{4 (1-\eta) \rho_F }.
\end{equation*}
Notice that $4 (\rho_F -12) \rho_F +153>0$ for all $\rho_F\in \R_{> 0}$, and recall that $\eta \leq \frac{1}{2}$ by assumption. Thus, for all $\rho_F>0$ there exists $\overline{T}\in \R_{>0}$ such that $\mathcal{M}_0(q,\rho_F) \succ 0$ for all $T\in(0, \overline{T})$, and $\mathcal{M}_0(q,\rho_F) \preceq 0$ for $T\geq \overline{T}$. 
Similarly, when $\delta>0$ we have that if $\mathcal{S}_\delta$ is $\rho_F$-GC, then $\mathcal{S}_0$ is also $\rho_F$-GC. Thus, we can conclude that for every $\rho_F$ and $\delta\ge 0$ there exists $\overline{T}$ such that $\mathcal{S}_\delta$ is not $\rho_F$-GC for any $T\ge \overline{T}$. 
\QEDB 
\end{example}
%
%

\vspace{0.1cm}
Using the global contractivity property of Definition \ref{definitionSC}, we have the following result for non-potential games.

\vspace{0.1cm}
\begin{thm}\label{theoremstrictmono}
Let $\mathcal{G}$ describe a \tcb{strictly monotone game}, and suppose that Assumptions \ref{cocoerciveassumption} and \ref{regularassumption} hold. Consider the HDS $\mathcal{H}_1$ under (RC$_1$) with $\rho_J\geq0$ and with reset policy $\alpha=\mathbf{1}_n$. If $\mathcal{S}_0$ is $\ell$-GC then the set $\mathcal{A}$ is R-UGAS, \tcb{for every $i\in\mathcal{V}$}, and for all solutions $x$ the following bound holds during flows
\begin{equation}\label{decrease:nonpotential:strict}
\abs{\mathcal{G}\tcb{(q(t,j))}}^2\leq \frac{\tilde{c}_j}{\tau_i^2(t,j)},~~~\forall~(t,j)\in \mathcal{T}(x),
\end{equation}
where $\set{\tilde{c}_j}\searrow 0^+$ is a sequence parameterized by $x(0,0)$.
\QEDB
\end{thm}
\vspace{0.1cm}

Unlike Theorem 1, in non-potential games the $\rho_F$-global-contractivity of $\mathcal{S}_\delta$ plays a fundamental role in the stability analysis of $\mathcal{H}_1$. In particular, the $\ell$-GC property of $\mathcal{S}_{\delta}$ will guarantee a suitable dissipativity property during flows via Lyapunov-based tools. While in Theorem 2 this is only a sufficient condition, the plots of Figure \ref{fig:phaseplot} indicate that keeping $T$ ``sufficiently small'' is also a necessary condition to preserve stability in non-potential games. In this figure, we show the phase plane of solutions to $\mathcal{H}_1$ with different values of $T$, in a game $\mathcal{G}$ given by \eqref{pseudogradient:necessityofresets}. It can be observed in the left plots that if $T$ is not small, divergent trajectories emerge. As shown in the right plots, the instability is removed by implementing sufficiently frequent resets.

\vspace{0.1cm}
Next, we provide a sufficient reset condition that guarantees that $\mathcal{S}_0$ is $\ell$-GC in any cocoercive \tcb{strictly monotone game}.

\vspace{0.1cm}
\begin{lemma}\label{sufficientlemma}
Suppose that Assumption \ref{cocoerciveassumption} holds, and $\mathcal{G}$ describes a \tcb{strictly monotone game}. Let $(\eta,T,\ell)$ satisfy the scalar inequality:
\begin{align}
0<T^2<\frac{1-\eta}{2\ell}\tag{$\text{RC}_2$}\label{RC2}.
\end{align}
Then, $\mathcal{S}_0$ is $\ell$-GC. \QEDB
\end{lemma}

\vspace{0.1cm}
%
%

%
%


\begin{figure*}[t!]
    \centering
    \includegraphics[width=\linewidth]{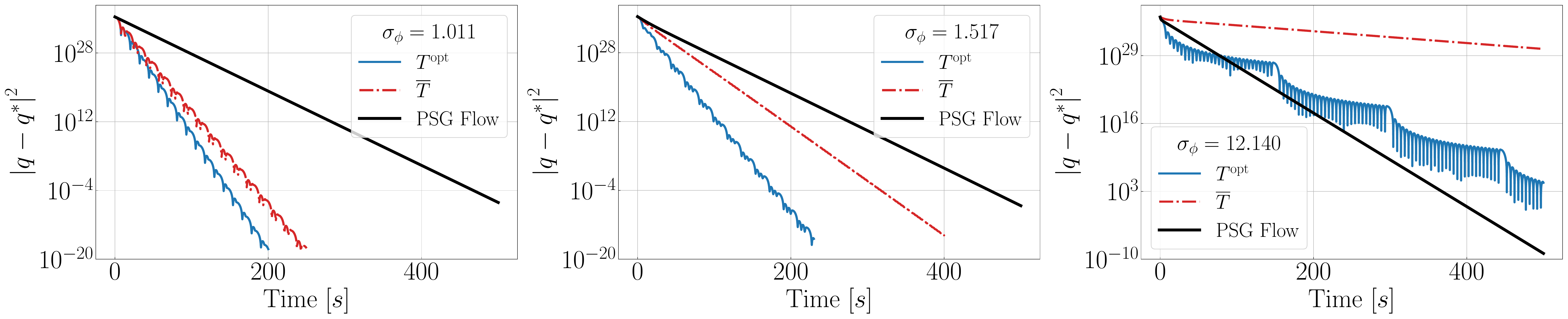}
    \caption{Distance to the Nash Equilibrium for trajectories resulting from the HM-NSS dynamics in  $5$-player $\kappa$-strongly monotone quadratic games with $\kappa=0.099$, and different condition numbers $\sigma_\phi$. }
    \label{fig:behaviorConditionNumbers}
\end{figure*}

We now turn our attention to games that are \tcb{$\kappa$-strongly monotone} and $\ell$-Lipschitz. For these games, we ask that the contractivity properties of $\mathcal{S}_{\delta}$ hold with $\delta>0$, and that (RC$_1$) holds with a particular value of $\rho_J$. We recall that the condition numbers $(\sigma_\phi,\sigma_r)$ are defined in \eqref{condition_numbers}.

\vspace{0.1cm}
\begin{thm}\label{theorem:strong:np}
Suppose that Assumption \ref{lipschitzassumption} holds and let $\mathcal{G}$ describe a \tcb{$\kappa$-strongly monotone game}. Consider the HDS $\mathcal{H}_1$ under (RC$_1$), and suppose that $\mathcal{S}_\delta$ is $(\sigma_\phi\ell)$-GC with $0< \delta < (1-\eta)/(\sigma_\phi\ell)$. Then, the following holds:
\begin{enumerate}[start=4,label={(i$_{\arabic*}$)}]
\item If $\alpha\in\{0,1\}^n$ and $\rho_J=0$, then $\mathcal{A}$ is R-UGES, and there exists $\lambda>0$ such that for each compact set $K\subset C_1\cup D_1$ there exists $M_0>0$ such that for all solutions $x$, with $x(0,0)\in K_0$, the bound \eqref{exp_property} holds.
\item If $\alpha=\mathbf{0}_n$ and $\rho_J=\sigma_\phi^2\kappa^{-1}$, then $\mathcal{A}$ is R-UGES and for each compact set $K_0\subset C_1\cup D_1$ there exists $M_0>0$ such that for all solutions $x$, with $x(0,0)\in K_0$, and for all $(t,j)\in \dom{x}$, the following bound holds:
\begin{equation*}
|q(t,j)-q^*|\leq\sigma_r\sigma_\phi\left(1-\gamma\left(\rho_J\right)\right)^{\frac{\alpha(j)}{2}}M_0,
\end{equation*}
where $\alpha(j)\coloneqq \max\{0,\lfloor\frac{j-n}{n}\rfloor\}$, and $\gamma\left(\rho_J\right)\in(0,1)$. \QEDB
\end{enumerate}
\end{thm}
%
%

\vspace{0.1cm}
Before commenting on the implications of Theorem \ref{theorem:strong:np}, we present a reset condition for \tcb{$\kappa$-strongly monotone game}s that is analogous to the one of Lemma \ref{sufficientlemma}.
\vspace{0.1cm}
\begin{lemma}\label{lemma_strong_suff1}
Suppose that Assumption \ref{lipschitzassumption} holds and that $\mathcal{G}$ describes a \tcb{$\kappa$-strongly monotone game}. Let $(\eta,T,\sigma_\phi\ell)$ satisfy the scalar inequality:
\begin{equation}\tag{$\text{RC}_3$}
0<T^2<\frac{1-\eta-\delta\sigma_\phi\ell}{\sigma_\phi\ell-\kappa+\delta(1-\eta-\delta\sigma_\phi\ell)},
\end{equation}
with $0\leq \delta < (1-\eta)/(\sigma_\phi\ell)$. Then
$\mathcal{S}_\delta$ is $(\sigma_\phi\ell)$-GC. \QEDB
\end{lemma}

\vspace{0.1cm}
Similar to (RC$_2$), the reset condition (RC$_3$) imposes an upper bound on the lengths of the intervals of flow of the HDS $\mathcal{H}_1$, now modulated by the \tcb{condition number $\sigma_\phi$ and the Lipschtiz constant $\ell$. }
%

\vspace{0.1cm}
\begin{remark}
When $\rho_J=\sigma_\phi^2\kappa^{-1}$, the conjunction of (RC$_1$) and (RC$_3$) imposes upper \emph{and} lower bounds for the reset times of the HDS $\mathcal{H}_1$ for all times $(t,j)\in\mathcal{T}(x)$. This result is in contrast to  potential games (and standard convex optimization problems) with periodic restarting where only a lower bound between resets is usually needed to achieve exponential convergence \cite{Candes_Restarting,zero_order_poveda_Lina}. Instead, Theorem \ref{theorem:strong:np} suggests that resets must occur in a particular frequency band: they should not occur too frequently (i.e., $T$ should not be too small) such that (RC$_1$) holds and the distance $|q-q^*|$ shrinks by a constant quantity after each interval of flow; however, resets should also happen frequently enough (i.e., $T$ should not be too large) such that $\mathcal{S}_{\delta}$ remains $(\sigma_\phi\ell)$-GC. \tcb{Whenever a resetting time $T$ is in such frequency band, i.e., whenever it satisfies (RC$_1$) and (RC$_3$), or (RC$_1$) and $\rho_F$-global contractivity of $\mathcal{S}_\delta$, we say that the $T$ is \emph{feasible}}.  
\QEDB
\end{remark}

\vspace{0.1cm}
The next lemma provides a sufficient condition to guarantee feasibility of the reset conditions of Theorem \ref{theorem:strong:np}. 

\vspace{0.1cm}
\begin{lemma}\label{sufficient_condition}
For any $\kappa>0,~\eta\leq1/2$ and $\sigma_\phi$ such that $\sigma_\phi^4-\sigma_\phi^2<2(1-\eta)$, there exists $(T,T_0)$ such that (RC$_1$) and (RC$_3$) simultaneously hold with $\rho_J=\sigma_\phi^2\kappa^{-1}$, provided that $\delta$ is sufficiently small. \QEDB
\end{lemma}

\begin{figure}[t!]
    \centering
    \includegraphics[width=0.6\linewidth]{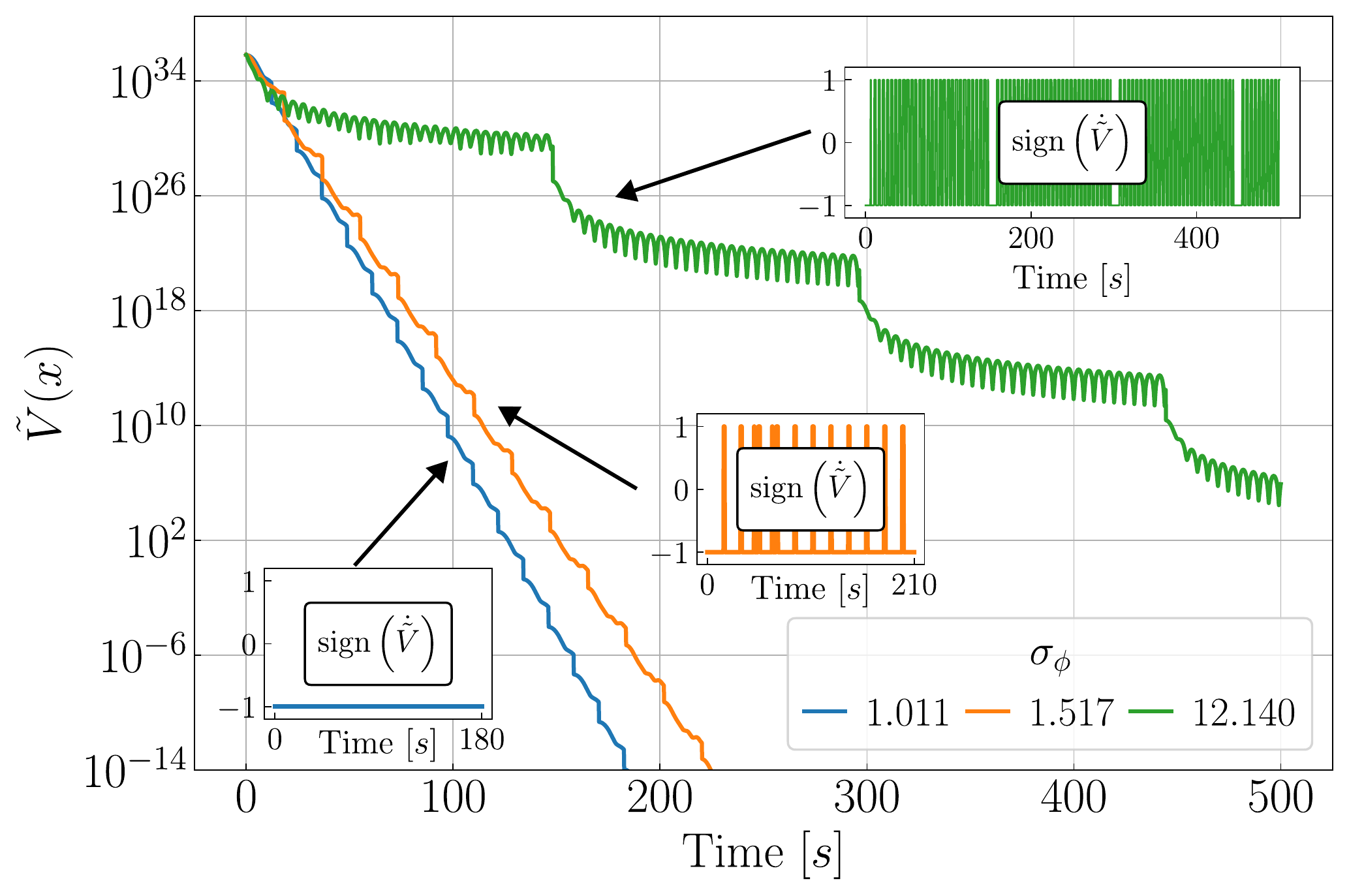}
    \caption{Lyapunov function $\tilde{V}(x)$ for trajectories resulting from the HM-NSS in $5$-player $\kappa$-strongly monotone quadratic games with $\kappa=0.009$, and different condition numbers $\sigma_\phi$. In all the cases the dynamics are implemented with the resetting period $T=T^{\text{op}}$ presented in Lemma \ref{quasi-optimal-restaring:np}.}
    \label{fig:dotVCondition}
\end{figure}

\vspace{0.1cm}
In Theorem \ref{theorem:strong:np}, the restarting policy $\alpha=\vec{0}_n$ leads to exponential NSS with rate of convergence characterized by $(1-\gamma(\sigma_\phi^2/\kappa))$. For this coefficient, one can choose a ``quasi-optimal'' restarting parameter $T$ to induce an acceleration-like property in \tcb{$\kappa$-strongly monotone game}s:

\vspace{0.1cm}
\begin{lemma}\label{quasi-optimal-restaring:np}
Under the Assumptions of Theorem \ref{theorem:strong:np}-(i$_5$), and for any $\nu>0$, the choice $T=T^{\text{opt}}:=e\sigma_\phi\sqrt{\frac{1}{2\kappa}+\frac{T_0^2}{\sigma_\phi^2}}$ guarantees that $\abs{q(t,j)-q^*}\leq \nu$ for all $t\ge t_{\nu^\text{opt}}$, where

\vspace{-0.3cm}
\begin{small}
\begin{align*}
        t^{\text{opt}}_\nu = \frac{1}{\eta}\left(e\sigma_\phi\sqrt{\frac{1}{2\kappa} + \frac{T_0^2}{\sigma_\phi^2}} - T_0\right)\ln\left(\frac{\sigma_\phi\sigma_rM_0}{\nu}\right),
\end{align*}
\end{small}

\vspace{-0.3cm}\noindent 
and $M_0$ is a constant that depends on $\abs{q(0,0)-q^*}$. Moreover, the convergence is of order $\mathcal{O}(e^{-\sqrt{\kappa}/\sigma_\phi})$ as $T_0\to0^+$. \QEDB
\end{lemma}
\begin{remark}\label{remark:optimal}
\tcb{The result of Lemma \ref{quasi-optimal-restaring:np} showcases the exponential bound induced by the HM-NSS dynamics: as $\sigma_\phi\to 1$, the convergence is of order $\mathcal{O}(e^{-\sqrt{\kappa}t})$, which is advantageous in games with low curvature and moderate condition number, see Figure \ref{fig:example3}. However, as $\sigma_\phi$ increases, the theoretical convergence bound deteriorates. In Figure \ref{fig:behaviorConditionNumbers}, we confirm this fact by implementing the HM-NSS in different $\kappa$-strongly non-potential monotone quadratic games with low curvature ($\kappa\ll 1$) and different condition numbers $\sigma_\phi\in\set{1.011,~ 1.517,~ 12.140}$. For these games, we compare two resetting times: 1) $T=T^{\text{opt}}$, which is only feasible for $\sigma_\ell=1.011$, and  2) $T=\overline{T}$, where $\overline{T}$ is the biggest resetting time found to guarantee that $\mathcal{S}_\delta$ is $(\sigma_\phi\ell)$-GC, with $\delta>0$ approaching to $0$, and is feasible for all $\sigma_\ell$. Experimentally, we observe that performance deteriorates slower when using $\overline{T}$ instead of $T^{\text{opt}}$, and acceleration can be achieved for a wider range of values of $\sigma_\ell$. However, when $T=T^{\text{opt}}$ is not feasible we cannot use the results of Theorem \ref{theorem:strong:np} to guarantee stability certificates under the implementation of the HM-NSS dynamics.
Whether or not a small condition number is a \emph{necessary} condition to achieve acceleration in games remains an open question.}
\end{remark}
\tcb{It is possible to find additional conditions on the game and the parameters of $\mathcal{H}_1$ to guarantee that $T^{\text{opt}}$ satisfies (RC$_1$) and (RC$_3$), i.e., that is \emph{feasible}. However, such conditions are rather involved and unintuitive, and therefore are omitted for brevity. Yet, we note that in Example 3 the quasi-optimal restarting $T^{\text{opt}}$ can be verified to be feasible. We also note that, based on numerical experiments, our theoretical bound is conservative, see Figure \ref{fig:example3}.}
\begin{remark}
\tcb{The stability results of Theorem 3-(i$_3$) are obtained by guaranteeing strong-decrease of a suitably chosen Lyapunov function $\tilde{V}$ during flows and jumps of the HM-NSS dynamics. In Figure \ref{fig:dotVCondition} we show the value of $\tilde{V}$ for different trajectories resulting from the HM-NSS dynamics in the same $\kappa$-strongly monotone games studied in Remark \ref{remark:optimal}, and using the resetting time $T^{\text{opt}}$ of Lemma \ref{quasi-optimal-restaring:np}. As seen in the figure, the Lyapunov function does not exhibit strong decrease during flows \bigg(with $\text{sign}\left(\dot{\tilde{V}}\right)$ changing multiple times\bigg) for the non-feasible values of $T^{\text{opt}}$, which correspond to the condition numbers $\sigma_{\phi} = 1.157$ and $\sigma_{\phi}=12.140$. However, the function does experience an overall decrease during the simulation time due to the strong decrease enforced through the jumps of the hybrid dynamics. Experimentally, we observe that to obtain better transient performance (with $T$s that are outside the feasible set), one could allow the Lyapunov function to increase during flows, provided this increase is compensated via jumps. This indicates interesting opportunities to attain even better behavior via adaptive resetting policies that reset whenever the Lyapunov function increases. However, the development of those techniques in a distributed way is challenging and falls out of the scope of the current paper.}
\end{remark}
\vspace{0.1cm}

\vspace{0.1cm}
We finish this section with a result for quadratic games.
\begin{prop}\label{propositionquadratic}
Suppose that Assumption \ref{lipschitzassumption} holds and $\mathcal{G}$ describes a \tcb{$\kappa$-strongly monotone quadratic game}. Consider the hybrid dynamics $\mathcal{H}_1$ under the reset condition (RC$_1$) with $\rho_J=\kappa^{-1}$ and reset policy $\alpha=\mathbf{0}_n$. If $\mathcal{S}_\delta$ is $\kappa$-GC, then $\mathcal{A}$ is R-UGES, and for each compact set $K_0\subset C_1\cup D_1$ there exists $M_0>0$ such that all solutions $x$ with $x(0,0)\in K_0$ satisfy
\begin{equation*}
|q(t,j)-q^*|\leq\sigma_r\sqrt{\sigma_\phi}\left(1-\gamma\left(\kappa^{-1}\right)\right)^{\frac{\alpha(j)}{2}}M_0,
\end{equation*}
where $\alpha(j)\coloneqq \max\{0,\lfloor\frac{j-n}{n}\rfloor\}$, and $\gamma\left(\kappa^{-1}\right)\in(0,1)$. \QEDB
%
%
\end{prop}
%


%

\vspace{0.1cm}
It is not difficult to show that the assumptions of Proposition \ref{propositionquadratic} can be satisfied when $T_0$ is sufficiently small and $\sigma_\phi<\sqrt{2}$. However, it is also not difficult to find games for which this condition is violated, yet (RC$_1$) holds and $\mathcal{S}_\delta$ is still $\kappa$-GC. In fact, for quadratic games, condition \eqref{contractivematrix} simplifies to $0\prec I_n-k_1(\kappa I_n-A)(\kappa I_n-A)^\top$, with $k_1=\frac{T^2}{\kappa(1-\eta)}$.\vspace{0.1cm}
\tcb{
\begin{remark}
For some games, a potential function $P$ might be available but its gradient $\nabla P$ might only be monotone if a suitable vector of weights $\omega\in \R^n$ is chosen so that $ D(\omega)\nabla P$ is strictly or $\kappa$-strongly monotone. In such cases, the results of Theorem \ref{theorem1} are not applicable, but those of Theorems \ref{theoremstrictmono} and \ref{theorem:strong:np} hold provided that $\mathcal{G}\coloneqq D(\omega)\nabla P$ satisfies the suitable assumptions on its continuity and contractivity of $S_\delta$. Such weighted potential games have been
recently studied in \cite{arcak2020dissipativity} n the context of congestion games. \QEDB
\end{remark}}

\begin{figure*}[b]
	\vspace{-5pt}
	\hrulefill
	\begin{small}
	\begin{equation}\label{upperBoundEpsilon}
	\varepsilon^*_\delta \coloneqq \frac{1}{2\sigma_{\mathcal{L}}\sqrt{n}}\left(1+\sigma_r^2\frac{\max\set{\frac{1}{T^2}+4\frac{\ell }{T\lambda_{\max}(\mathcal{L})} ~,~2 + 2\frac{\ell }{T\lambda_{\max}(\mathcal{L})}}}{\delta\min\set{1,\zeta^2}}\right)^{-1}
	\end{equation}
	\end{small}
	\vspace{-5pt}
\end{figure*}
\section{HYBRID MOMENTUM-BASED NSS WITH PARTIAL INFORMATION}
\label{Sec_R2A}
In Section \ref{Sec_R1}, we assumed that each player $i$ had access to an individual Oracle able to generate measurements of the partial gradient $\frac{\partial \phi_i(\cdot)}{\partial q_i}$ at the \emph{overall} state $q$ of the game. In this section, we now relax this assumption by considering Oracles that provide \emph{evaluations} of the partial derivatives. In this case, to carry out the gradient evaluations each player needs to estimate the actions of the other players by extending the hybrid dynamics $\mathcal{H}_1$. 

\vspace{-0.1cm}
\subsection{Individual Multi-Time Scale Hybrid Dynamics}
To achieve distributed NSS over graphs with partial information, we proceed to endow each player $i\in\mathcal{V}$ with an auxiliary state $\mathbf{e}^i$ that serves as an individual estimation of the actions of the other players:
$$\mathbf{e}^i\coloneqq (\mathbf{e}_1^i,\mathbf{e}_2^i,\ldots,q_i,\ldots,\mathbf{e}_{n-1}^i,\mathbf{e}_n^i)\in\R^{n}.$$ 
Since players do not need to estimate their own action, it is also convenient to introduce the auxiliary state $\vec{e}_{-i}^i\in\mathbb{R}^{n-1}$ which contains the same entries of $\mathbf{e}^i$ with the exception of $q_i$, which is removed. Using this notation, we now assume that each player $i$ has access to individual \emph{gradient Oracles} characterized by mappings of the form $(q_i,\mathbf{e}^i_{-i})\mapsto\hat{\mathcal{G}}_i(q_i,\mathbf{e}^i_{-i})$, which satisfy $\hat{\mathcal{G}}_i(q_i,q_{-i})=\frac{\partial \phi_i(q)}{\partial q_i}$. Following similar notation used in the literature \cite{PavelGames}, we define the matrices
\begin{align*}
\mathcal{Q}_i&\coloneqq \left[\begin{array}{ccc}
I_{(i-1)} & 0_{(i-1)\times 1} & 0_{(i-1)\times (n-i)}\\
0_{(n-i)\times (i-1)} & 0_{(n-i)\times 1} & I_{(n-i)}\\
\end{array}\right],\\
\mathcal{P}_i&\coloneqq \left[0_{1\times (i-1)}~~1~~0_{1\times (n-i)} \right].
\end{align*}
%
%
By using these definitions, each player $i$ now implements the following accelerated augmented continuous-time NSS dynamics: 
\begin{align}\label{continuostimewithestimation}
\left(\begin{array}{c} 
\vspace{0.1cm}\dot{q}_{i}\\
\vspace{0.1cm}\dot{p}_{i}\\
\vspace{0.1cm}\dot{\tau_i}\\
\dot{\vec{e}}_{-i}^i
\end{array}
\right)=\left(\begin{array}{c} 
\frac{2}{\tau_i}(p_{i}-q_{i})+\mathcal{P}_i\sum_{j\in\mathcal{N}_i}(\mathbf{e}^i-\mathbf{e}^j)\\
-2\tau_i\hat{\mathcal{G}}_i(q_i,\mathbf{e}^i_{-i})\\
\eta\\
-\frac{1}{\varepsilon}\mathcal{Q}_i\sum_{j\in\mathcal{N}_i}(\mathbf{e}^i-\mathbf{e}^j)
\end{array}
\right),
\end{align}
where $\varepsilon>0$ is a new tunable parameter. These dynamics  are implemented whenever the state $\tau_i$ satisfies $\tau_i\in[T_0,T)$. 
The momentum-based dynamics \eqref{continuostimewithestimation} implement a dynamic consensus mechanism with state $\vec{e}_{-i}^i$. This mechanism uses a high gain $\frac{1}{\varepsilon}$ to induce a time-scale separation in the flows of the hybrid algorithm. In particular, if the states $\mathbf{e}^i$ were to instantaneously achieve their steady state value, the flows \eqref{continuostimewithestimation} would reduced to the flows \eqref{flowmap00}. 

When players are uncoordinated, the individual resets are triggered by the condition $\tau_i=T$, and are given by
\begin{equation}\label{uncoordinatedresets2}
x_i^+=R_i(x_i),~~~\vec{e}^{i+}_{-i}=\vec{e}_{-i}^i,
\end{equation}
where $R_i$ is defined in \eqref{decoupled_flows2}. 
\begin{figure}
    \centering
    \includegraphics[width=0.65\linewidth]{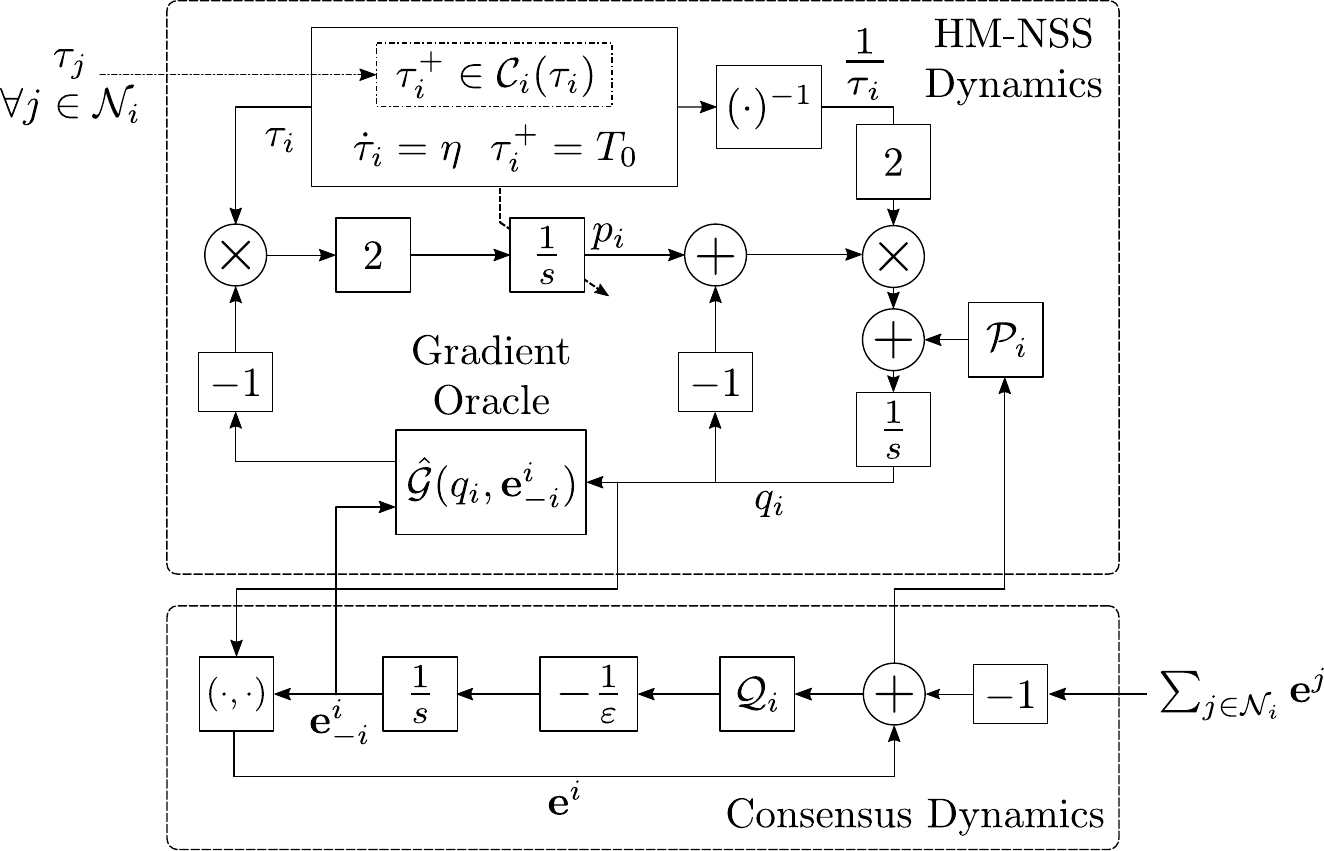}
    \caption{Scheme of Individual HM-NSS dynamics for games with partial information. Consensus dynamics are implemented to estimate the actions of other players. }
    \label{blockdiagram2}
    \vspace{-0.3cm}
\end{figure}
However, lack of coordination between resets can induced the same issues discussed in Example \ref{negative_example2}. To avoid this, we will incorporate the hybrid coordinated restarting mechanism described in Section \ref{sec_coordination_mechanism}. Figure \ref{blockdiagram2} shows a block-diagram representation of the multi-time scale hybrid dynamics of each player.
\vspace{-0.1cm}
\subsection{Well-Posed Coordinated Hybrid System with Partial Information}
To write the coordinated HDS in vectorial form, we introduce the matrices $\mathcal{Q}\coloneqq D(\mathcal{Q}_i)\in \R^{(n^2-n)\times n^2}$ and $\mathcal{P}\coloneqq D(\mathcal{P}_i)\in \R^{n\times n^2}$, and note that $q=\mathcal{P}\mathbf{e}\in\R^{n^2-n}$. Additionally, we define the state $\hat{q}\coloneqq\mathcal{Q}\vec{e}$, such that using $\mathcal{P}\mathcal{P}^\top =  I_n,~\mathcal{Q}\mathcal{Q}^\top= I_{n^2 -n},\text{ and }\mathcal{P}\mathcal{Q}^\top = 0$,
we can write $\mathbf{e}=\psi(q,\hat{q}):=\mathcal{P}^\top q+\mathcal{Q}^\top \hat{q}$, where $\vec{e}=\left(\vec{e}^1,\cdots,\vec{e}^n\right)$, and express the overall hybrid NSS dynamics as a HDS \eqref{HDS} with state $(x,\hat{q})$, where $x\coloneqq (q,p,\tau)\in\R^{3n}$, and data $\mathcal{H}_2=(C_2,F_2,D_2,G_2)$. The flows are given by
\begin{align}\label{flowEstimation}
    \begin{pmatrix}
        \dot{q}\\\dot{p}\\\dot{\tau}\\\dot{\hat{q}}
    \end{pmatrix}{=}F_2(x,\hat{q})&{\coloneqq}\begin{pmatrix}
        2D(\tau)^{-1}(p-q) - \mathcal{P}\vec{L}\psi(q,\hat{q}) \\
        -2D(\tau)\hat{\mathcal{G}}(\psi(q,\hat{q}))\\
    \eta\vec{1}_n\\
        -\frac{1}{\varepsilon}\mathcal{Q}\vec{L}\psi(q,\hat{q})
    \end{pmatrix},
\end{align}
where $\vec{L} \coloneqq \mathcal{L}\otimes I_n$ denotes the communication matrix of the graph $\mathbb{G}$. The continuous-time dynamics in \eqref{flowEstimation} are allowed to evolve whenever $(x,\hat{q})$ belongs to the flow set:
\begin{align}\label{flowSetEstimation}
C_2\coloneqq &\Big\{(x,\hat{q})\in\R^{n^2 + 2n}:~q\in\R^n,p\in\R^n,\notag\\
&~~~~~~~~~~~~~~~~~~\tau\in[T_0,T]^n,~\hat{q}\in\R^{n^2-n} \Big\}.
\end{align}
On the other hand, the jump set is defined as:
\begin{equation}\label{jumpSetEstimation}
D_2\coloneqq \Big\{(x,\hat{q})\in\R^{n^2 + 2n}:x\in C,~\max_{i\in\mathcal{V}}\tau_i=T\Big\},
\end{equation}
and the discrete-time dynamics of the algorithm are given by:
\begin{equation}\label{jumpEstimation}
(x^+,\hat{q}^+)\in G_2(x,\hat{q})\coloneqq G_1(q,p,\tau)\times\{\hat{q}\},
\end{equation}
where $G_1$ is defined as in \eqref{jump_map}. Thus, the jump map affects only the momentum coefficients $\tau$ and the state $p$ via the reset policy $\alpha$. Similar to Lemma \ref{lemmawellposed}, the next lemma follows directly by construction of the HDS.

\vspace{0.1cm}
\begin{lemma}\label{lemma4}
For the HDS $\mathcal{H}_2\coloneqq (C_2,F_2,D_2,G_2)$, all the properties of Lemma \ref{lemmawellposed} still hold. \QEDB
\end{lemma}
\vspace{0.1cm}

We will study the stability properties of the HDS $\mathcal{H}_2$ with respect to the following compact set:
\begin{equation}\label{compactgraphset}
    \mathcal{A}_{\mathbb{G}}\coloneqq \mathcal{A}\times\{\mathcal{Q}(\vec{1}_n\otimes q^*)\},
\end{equation}
where $\mathcal{A}$ was defined in \eqref{definitionofset}. In this case, we will use the following restricted reverse-Lipschitz assumption, also used in \cite{povedaBianchi2020} for optimization, and in \cite{gadjov2022exact} for NES with static inertia.

\vspace{0.1cm}
\begin{assumption}\label{assumption:reverse:lipschitz}
There exists $\zeta>0$ such that $\abs{\mathcal{G}(q) -\mathcal{G}(q^*)} \ge \zeta\abs{q-q^*}$, for all $q\in\R^n$.\QEDB
\end{assumption}
\vspace{0.1cm}
Next, we have the following result, which leverages items (i$_1$)-(i$_5$) of Theorems 1-3.

\vspace{0.1cm}
\begin{thm}\label{theorem:partialInformation}
Let $\mathcal{G}$ describe a \tcb{strictly monotone game}. Suppose that Assumptions \ref{cocoerciveassumption}, \ref{regularassumption} and \ref{assumption:reverse:lipschitz} hold, and consider the HDS $\mathcal{H}_2$ under (RC$_1$). If $\mathcal{S}_\delta$ is $\ell-$GC with $0<\delta<(1-\eta)/\ell$, then under any of the conditions (i$_1$)-(i$_5$) the following holds:
 \begin{enumerate}[(a)]
 \item For all $\varepsilon\in(0,\varepsilon_\delta^*)$, where $\varepsilon_\delta^*$ is given by \eqref{upperBoundEpsilon}, 
 the set $\mathcal{A}_{\mathbb{G}}$ is R-UGAS.
 \item For each $(\hat{t},\hat{j},\nu)\in\mathbb{R}_{>0}^3$ and each compact set $K_x\times K_{\hat{q}}\subset C_2\cup D_2$, there exists $\varepsilon^{**}$ such that for each $\varepsilon\in(0,\varepsilon^{**})$ and each solution of $\mathcal{H}_2$ with $x(0,0)\in K_x$ and $\hat{q}(0,0)\in K_{\hat{q}}$, there exists a solution $\tilde{x}$ of system $\mathcal{H}_1$ with $\tilde{x}\in K_x$ such that $x$ and $\tilde{x}$ are $(\hat{t},\hat{j},\nu)$-close.  \QEDB 
 \end{enumerate}
%
%
%
%
%
%
%

\end{thm}

\begin{figure}[t!]
    \centering
    \includegraphics[width=0.65\linewidth]{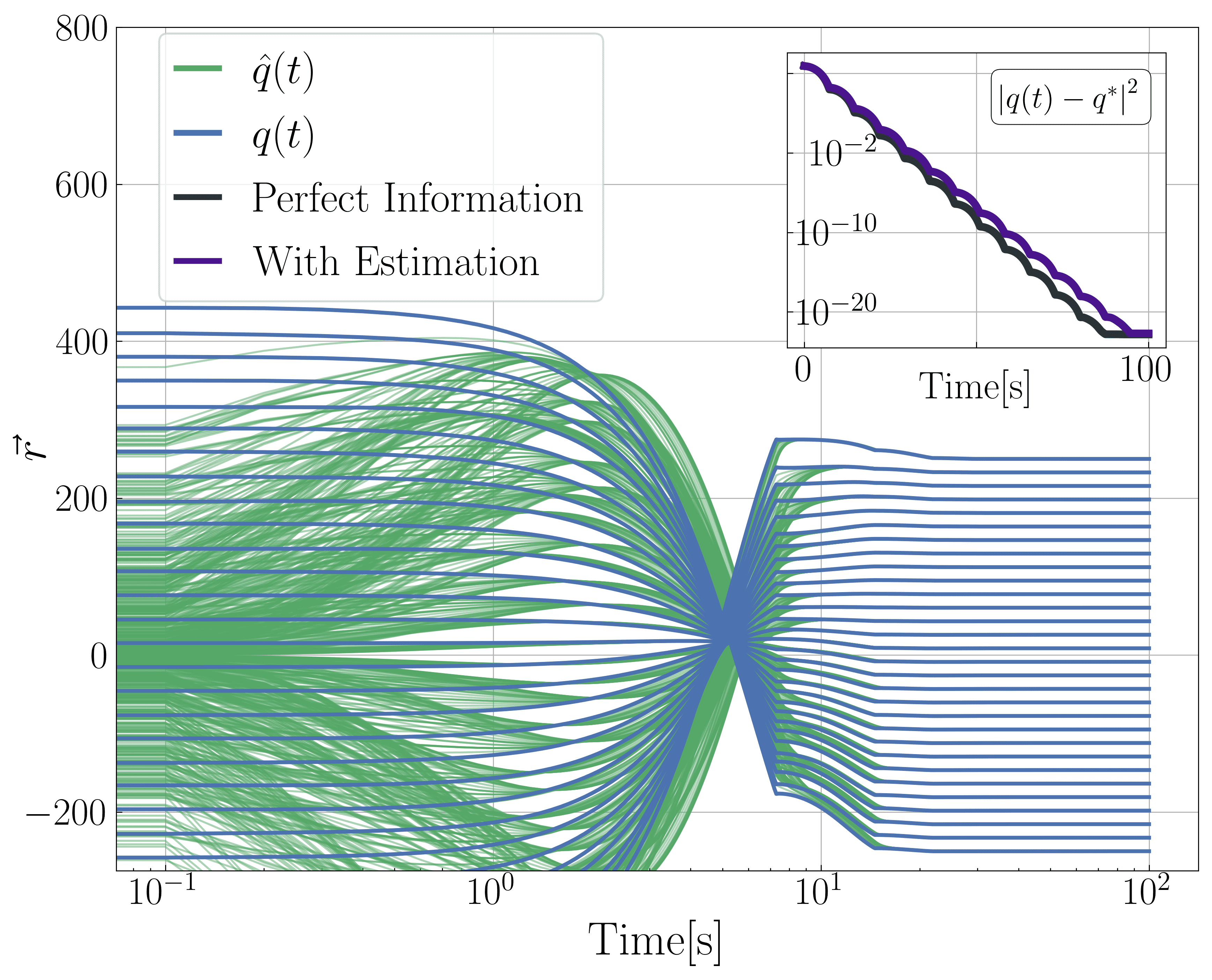}
    \caption{Trajectories of actions $q$ and estimates $\hat{q}$  in a non-potential \tcb{$\kappa$-strongly monotone quadratic game} with $n=30$, $\kappa=0.01,~\ell=0.1$, $\tau_s(0,0)=0.1\cdot\vec{1}_n$, and $\varepsilon=5\times 10^{-3}$. The inset shows the distance to the NE.}
    \label{fig:partialInformationPlot}
    \vspace{-0.3cm}
\end{figure}

\vspace{0.1cm}
Item (a) of Theorem \ref{theorem:partialInformation} establishes  robust stability and convergence properties for the hybrid NSS dynamics $\mathcal{H}_2$. On the other hand, item (b) establishes that on compact sets of initial conditions and on compact time domains, the trajectories $x$ will behave as the trajectories of the ``full-information'' system $\mathcal{H}_1$ as $\varepsilon\to0^+$ in \eqref{flowEstimation}. In particular, by combining items (a) and (b), we recover the convergence bounds of Theorems \ref{theorem1}, \ref{theoremstrictmono}, and \ref{theorem:strong:np}, now in a semi-global practical sense as $\varepsilon\to0^+$. 
%
%
%
%
%
%
%
%
This behavior is illustrated in Figure \ref{fig:partialInformationPlot}, which shows the trajectories $q$ and $\hat{q}$ in a \tcb{$\kappa$-strongly monotone game}. As it can be observed, the solutions of $\mathcal{H}_2$ approximate the solutions of $\mathcal{H}_1$ as $\varepsilon\to0^+$.

\vspace{0.1cm}
\begin{remark}
Assumption \ref{assumption:reverse:lipschitz} always holds for \emph{\tcb{$\kappa$-strongly monotone game}s} with $\zeta = \kappa$. Hence, for these games one can compute an alternative expression of $\varepsilon_{\delta}^*$ by substituting  Assumptions \ref{cocoerciveassumption}-\ref{assumption:reverse:lipschitz} in Theorem \ref{theorem:partialInformation} by Assumption \ref{lipschitzassumption} when $\mathcal{S}_{\delta}$ is $(\sigma_\phi\ell)$-GC. Moreover, to guarantee that $\mathcal{S}_\delta$ is $\ell$-GC, required in Theorem \ref{theorem:partialInformation}, a suitable upper bound for $T$ can be obtained by mirroring the derivations of Lemma \ref{lemma_strong_suff1}, which we omit here due to space limitations.
\QEDB
\end{remark}
%

To our best knowledge, Theorem \ref{theorem:partialInformation} is the first result in the literature that establishes robust convergence and stability properties for decentralized momentum-based NSS algorithms over graphs. Note that the stable incorporation of the multiple-time scale consensus mechanism is enabled by the use of resets, since otherwise no $\mathcal{K}\mathcal{L}$ bound (or strong Lyapunov function) would exist for the slow dynamics (also called ``reduced dynamics'' in the singular perturbation literature \cite[Ch. 11]{khalil}) of the flows, which are precisely given by \eqref{ODEmomentum}. 

\section{MODEL-FREE NSS WITH MOMENTUM}
\label{Sec_ESC}
In the previous sections, we studied algorithms that made use of gradient Oracles with full or partial information. In this section, we dispense with this assumption by designing accelerated \emph{model-free} hybrid NSS dynamics, suitable for applications where players have access only to real-time \emph{measurements} of the signals that correspond to their cost functions $\phi_i$ (e.g., the difference between the individual cost and revenue in a market), which are generated by the game. Such payoff-based algorithms can be designed using tools recently developed in the context of hybrid equilibrium seeking control \cite{PoTe16,zero_order_poveda_Lina}. 
%

%
\vspace{-0.1cm}
\subsection{Model-Free NSS Dynamics}
To achieve model-free NSS, each player $i$  generates an individual probing signal $t\mapsto\mu_{i}(t)$, obtained as the solution of a dynamic oscillator with state $\mu_i\coloneqq (\tilde{\mu}_{i},\hat{\mu}_{i})\in\R^2$, evolving on the unit circle $\mathbb{S}^1$ according to
\begin{equation}\label{oscillator}
\dot{\mu}_i=\frac{1}{\varepsilon_p}\mathcal{R}_i \mu_i,~~\mu_i\in\mathbb{S}^1,~~~\mathcal{R}_i\coloneqq 2\pi\varsigma_i\left[\begin{array}{cc}
0 & 1\\
-1 & 0
\end{array}\right],
\end{equation}
where $\varepsilon_p$ and $\varsigma_i$ are positive tunable parameters. Note that $\mathbb{S}^1$ is forward invariant under the dynamics of $\mu_i$. Using this probing signal, each player implements the flows:
\begin{align}\label{decoupled_flows1}
\left(\begin{array}{c} 
\dot{q}_{i}\\
\dot{p}_{i}\\
\dot{\tau_i}\\
\end{array}
\right)=\left(\begin{array}{c} 
\frac{2}{\tau_i}(p_{i}-q_{i})\\
-\frac{4}{\varepsilon_a}\tau_i \phi_i(q+\varepsilon_a\tilde{\mu})e^\top_1\mu_i\\
\eta
\end{array}
\right),
\end{align}
where $\mu=(\mu_1,\mu_2,\ldots,\mu_n)\in\mathbb{R}^{2n}$, and where $\tilde{\mu}$ is the vector that contains the odd components of $\mu$. The dynamics \eqref{decoupled_flows1} use real-time \emph{measurements} of the cost $\phi_i$, and are implemented whenever $\tau_i\in[T_0,T)$. Conversely, when $\tau_i=T$ and players are uncoordinated, they reset their states according to the dynamics $x^+_{i}=R_i(x_i)$,~ $\mu^+_i=\mu_i$, where $R_i$ is defined as in \eqref{decoupled_flows2}. The constant  $\varepsilon_a>0$ is also a tunable parameter. 

We impose the following assumption on the parameters $\varsigma_i$ of \eqref{oscillator}, which is standard in the averaging-based NES literature \cite{Frihauf12a,Ye:16,Poveda:15,Dither_ReUse}.

\vspace{0.1cm}
\begin{assumption}
For all $i$, $\varsigma_i$ is a positive rational number, $\varsigma_i\neq \varsigma_j$, $\varsigma_i\neq 2\varsigma_j$, $\varsigma_i\neq 3\varsigma_j$, for all $i\neq j\in\mathcal{V}$. \QEDB
\end{assumption}
\vspace{0.1cm}
%

%
%
%
As in the model-based case, an uncoordinated implementation of the model-free hybrid dynamics can be detrimental to the stability and/or transient performance of the algorithm. Thus, we incorporate the hybrid coordination mechanism described in Section \ref{sec_coordination_mechanism} to coordinate the resets of the players. Using the set-valued coordination mechanism, the overall  discrete-time dynamics of the algorithm are modeled by the difference inclusion
\begin{equation}\label{jumpES}
(x^+,\mu^+)\in G_{3}(x,\mu)\coloneqq G_1(x)\times\{\mu\},
\end{equation}
where $G_1$ is given by \eqref{jump_map}. This jump map will preserve the sequential nature of the resets needed to guarantee a well-posed HDS that satisfies \eqref{limitcondition1} and \eqref{limitcondition2}. Using $\bar{\mathbf{\phi}}\coloneqq (\phi_1,\phi_2,\ldots,\phi_N)$, the continuous-time dynamics of the model-free hybrid NSS algorithm can be written in vector form as:
\begin{align}\label{decoupled_flows11}
\left(\begin{array}{c} 
\dot{q}\\
\dot{p}\\
\dot{\tau}\\
\dot{\mu}
\end{array}
\right)=F_{3}(x,\mu)\coloneqq \left(\begin{array}{c} 
2D(\tau)^{-1}(p-q)\\
-\frac{4}{\varepsilon_a}D(\tau)\bar{\phi}(q+\varepsilon_a\tilde{\mu})\tilde{\mu}\\
\eta\\
\frac{1}{\varepsilon_p}D(\mathcal{R}_i) \mu
\end{array}
\right),
\end{align}
and the flow and jump sets are defined as:
\begin{equation}\label{sets_ES}
C_{3}\coloneqq C_1\times\mathbb{T}^n~~\text{and}~~D_{3}\coloneqq D_1\times\mathbb{T}^n.
\end{equation} \noindent 
\begin{figure}[t!]
    \centering
    \includegraphics[width=0.55\linewidth]{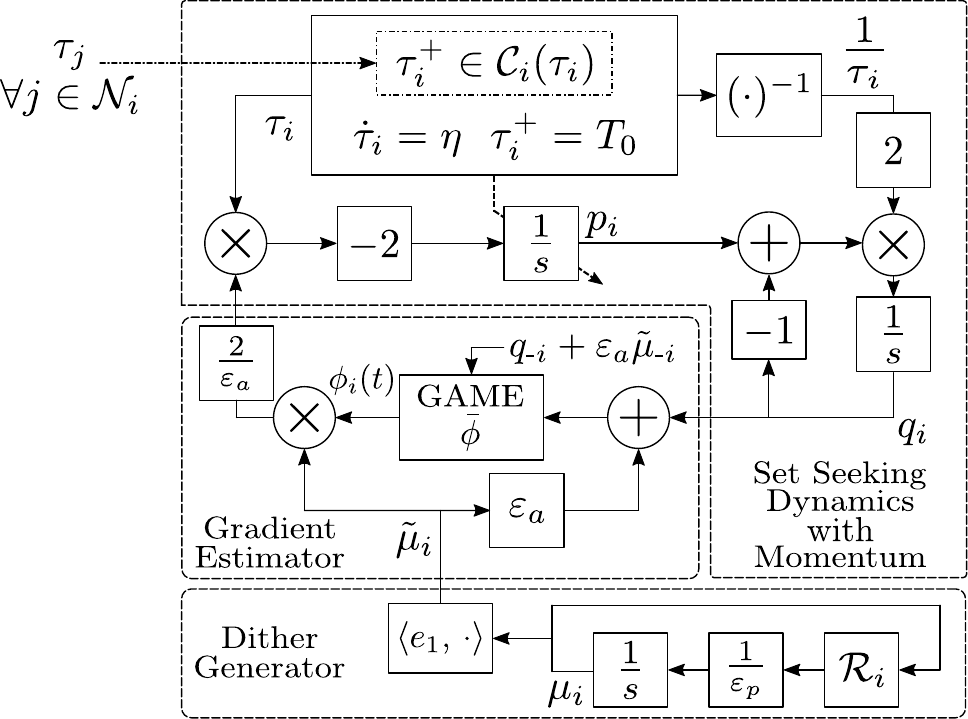}
    \caption{Scheme of Individual Model-Free HM-NSS dynamics with real-time measurements of the cost.}
    \label{blockdiagram3}
    \vspace{-0.4cm}
\end{figure}

\vspace{-0.3cm}
\noindent 
Figure \ref{blockdiagram3} shows a scheme of the proposed model-free NSS dynamics.

\vspace{-0.3cm}
\subsection{Semi-Global Practical Stability Results}
The data $\mathcal{H}_3=(C_3,F_3,D_3,G_3)$ defines the third hybrid NSS dynamics considered in this paper. The stability and convergence properties of $\mathcal{H}_3$ are given in the following theorem, which also leverages items (i$_1$)-(i$_5$) of Theorems 1-3.
%

\vspace{0.1cm}
\begin{thm}\label{theoremES}
Let $\mathcal{G}$ describe a \tcb{strictly monotone game}, and consider the HDS $\mathcal{H}_{3}$ under (RC$_1$). Then, under any of the conditions (i$_1$)-(i$_5$) the following holds:
\begin{enumerate}[(a)]
\item The set $\mathcal{A}\times\mathbb{T}^n$ is SGPAS as $(\varepsilon_p,\varepsilon_a)\to0^+$.
\item For each $(\hat{t},\hat{j},\nu)\in\mathbb{R}^3$ and each compact set $K_x\subset C_1\cup D_1$,  $~\exists~\varepsilon_a^*>0$ s.t. $\forall~\varepsilon_a\in(0,\varepsilon_a^*)$ $\exists~\varepsilon_p^*>0$ s.t. $\forall~\varepsilon_p\in(0,\varepsilon_p^*)$, and for each trajectory $x$ of system $\mathcal{H}_{3}$ with $x(0,0)\in K_x$ there exists a solution $\tilde{x}$ of system $\mathcal{H}_1$ such that $x$ and $\tilde{x}$ are $(\hat{t},\hat{j},\nu)$-close. \QEDB
\end{enumerate}
\end{thm}

\vspace{0.1cm}
Similar to Theorem \ref{theorem:partialInformation}, the result of Theorem \ref{theoremES} establishes two key properties: First, for any desired precision $\nu>0$, and any compact set of initial conditions $K_x$, every solution of the HDS $\mathcal{H}_{3}$ initialized in $K_x$ will satisfy a bound of the form\footnote{We note that $|\mu(t,j)|_{\mathbb{T}^n}=0$ for all $(t,j)$ in the domain of the solutions.}
\begin{equation}
|x(t,j)|_{\mathcal{A}}\leq \beta(|x(0,0)|_{\mathcal{A}},t+j)+\frac{\nu}{2},    
\end{equation}
with $\beta\in\mathcal{K}\mathcal{L}$, provided the parameters $\varepsilon_a$ and $\varepsilon_p$ are selected sufficiently small.
\tcb{Moreover, item (b) implies that by selecting $\varepsilon_a$ and $\varepsilon_p$ sufficiently small, the trajectories $x$ of $\mathcal{H}_{3}$ will recover all the fast convergence bounds established in Section \ref{main_results_first}, modulo a small residual error and on compact sets.}

\begin{figure}[t!]
    \centering
    \includegraphics[width=0.65\linewidth]{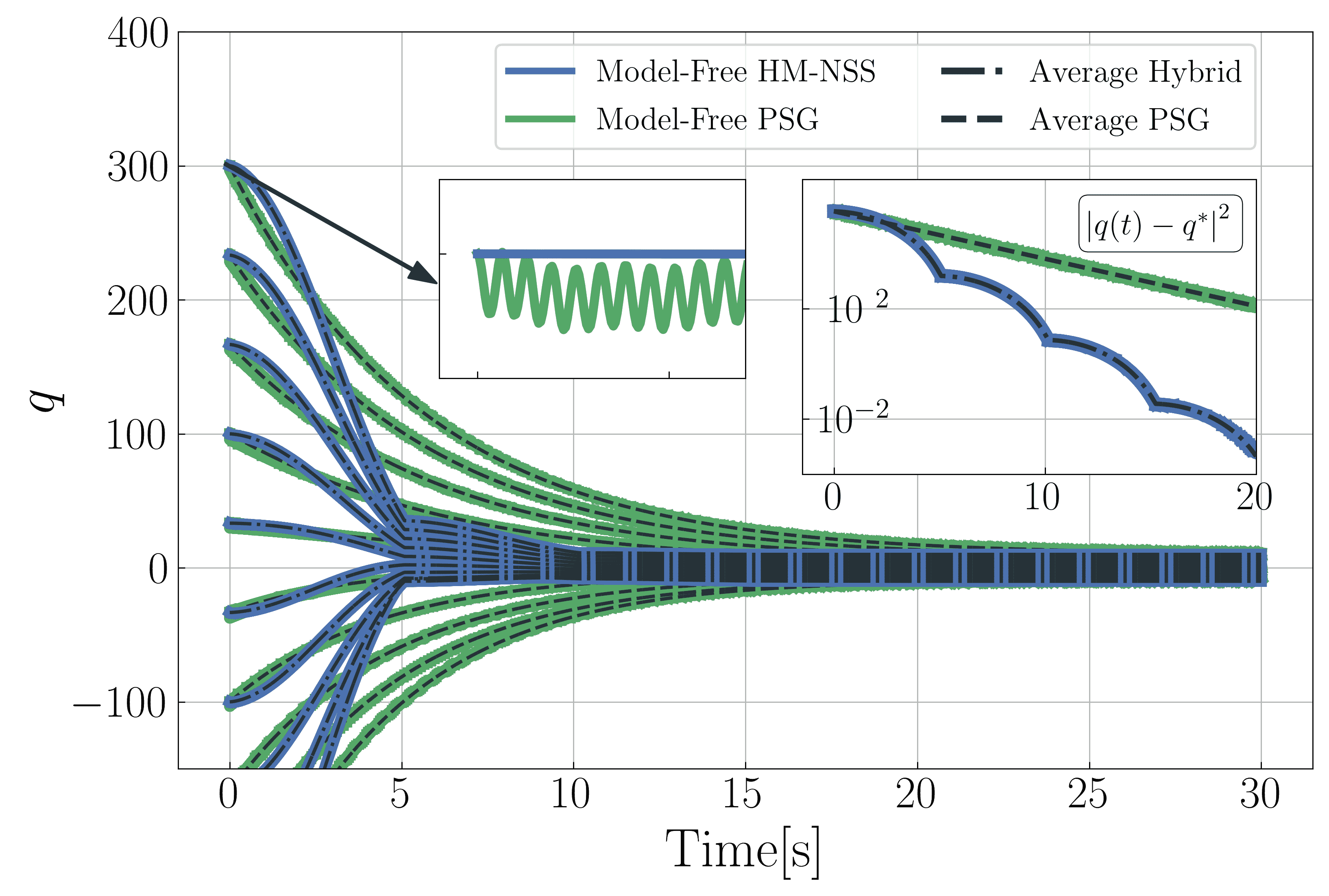}
    \caption{Trajectories of $\mathcal{H}_1$ and $\mathcal{H}_3$ in a non-potential \tcb{$\kappa$-strongly monotone quadratic game} with $\kappa = 0.197,~\ell=0.2$ and $n=10$.}
    \label{fig:es}
    \vspace{-0.4cm}
\end{figure}

\vspace{0.1cm}
\begin{remark}
The model-free dynamics $\mathcal{H}_3$ are based on averaging theory for (perturbed) hybrid systems \cite{WangTeelNesic,zero_order_poveda_Lina}. In particular, as $\varepsilon_a,\varepsilon_p\to0^+$ the trajectories of $\mathcal{H}_3$ behave as their average hybrid dynamics (modulo a small perturbation), which are precisely given by $\mathcal{H}_1$. Both dynamics are set-valued, which differs from existing results in the literature of model-free Nash seeking \cite{Frihauf12a}. Figure \ref{fig:es} shows a numerical experiment in a \tcb{$\kappa$-strongly monotone quadratic game} where a solution of $\mathcal{H}_3$ is compared to a solution of the model-free dynamics of \cite{Frihauf12a}. \QEDB
\end{remark}


We finish this section by commenting on the extensions of system $\mathcal{H}_3$ to applications where players could have access to an individual ``Black-Box Oracle'' that allows them to \emph{evaluate} (as opposed to measure) their local cost $\phi_i$ at their current state $q_i$, using estimations of the actions of the other players (in a similar spirit as in Section \ref{Sec_R2A}) and without knowledge of the mathematical form of $\phi_i$ (e.g., using dynamic simulators). In this case, we can follow the same approach of Section \ref{Sec_R2A}, incorporating an auxiliary estimation state $\hat{q}$. In this, the vectorial hybrid system $\mathcal{H}_4=(C_4,F_4,D_4,G_4)$ will have flow map given by
\begin{align}\label{decoupled_flows_es111}
\left(\begin{array}{c} 
\dot{q}\\
\dot{p}\\
\dot{\tau}\\
\dot{\mu}\\
\dot{\hat{q}}
\end{array}
\right)=F_{4}(\zeta)\coloneqq \left(\begin{array}{c} 
2D(\tau)^{-1}(p-q)- \mathcal{P}\vec{L}\psi(q,\hat{q})\\
-\frac{4}{\varepsilon_a}D(\tau)\bar{\phi}(\psi(q+\varepsilon_a\tilde{\mu},\hat{q}))\tilde{\mu}\\
\eta\\
\frac{1}{\varepsilon_p}D(\mathcal{R}_i) \mu\\
 -\frac{1}{\varepsilon_c}\mathcal{Q}\vec{L}\psi(q,\hat{q})
\end{array}
\right),
\end{align}
jump map $G_{4}(x,\mu,\hat{q})\coloneqq G_1(x)\times\{\mu\}\times\{\hat{q}\}$, flow set $C_{4}\coloneqq C_1\times\mathbb{T}^n\times\mathbb{R}^{n^2-n}$ and jump set $D_{4}\coloneqq D_1\times\mathbb{T}^n\times\mathbb{R}^{n^2-n}$. For this hybrid system, a result like Theorem \ref{theoremES}-(a) also holds, now with respect to the set $\mathcal{A}\times \mathbb{T}^n\times\{\mathcal{Q}(\vec{1}_n\otimes q^*)\}$ and as $(\varepsilon_p,\varepsilon_a,\varepsilon_c)\to0^+$. Similarly, a result like Theorem \ref{theoremES}-(b) holds by noting that the average hybrid dynamics of $\mathcal{H}_4$ are precisely given by the HDS $\mathcal{H}_2$ (modulo a small perturbation on the gradient), whose robust stability properties were already established in Section \ref{Sec_R2A}. Thus, we can follow exactly the same steps of the proof of Theorem \ref{theorem:partialInformation} to obtain an equivalent result.

\section{ANALYSIS AND PROOFS}
\label{sec_proofs}

\subsection{Proofs of Section \ref{Sec_R1}}
\label{proofsSection4}

\tcb{In this section we present condensed proofs of the Theorems and Lemmas introduced in the previous sections.
For additional steps and equations we refer the reader to the extended manuscript \cite{extended_manuscript}.}

%
%
\noindent
\textbf{Proof of Lemma \ref{lemmawellposed}:} Well-posedness follows directly by \cite[Thm. 6.30]{HDS}, since $F_1$ is continuous, $C_1$ and $D_1$ are closed sets, and $G_1$ is outer-semicontinuous  (OSC) and locally bounded (LB) in $D_1$. To rule out finite escape times it suffices to study the behavior of the states $(q,p)$. Using Assumption \ref{lipschitzassumption}, the form of \eqref{flowmap00}, and the fact that $\mathcal{G}(q^*)=0$, we have that $|\dot{q}|\leq \frac{2}{T_0}|p-q|$, $|\dot{p}|\leq 2T \ell|q-q^*|~\forall q^*\in\mathcal{A}_{NE}$, which implies that $|(\dot{q},\dot{p})|\leq \tilde{\ell} |(q,p)-(q^*,q^*)|$, with $\tilde{\ell}\coloneqq 2\sqrt{3}\max\{\frac{1}{T_0},T \ell\}$. \tcb{Then, by using the fact that $\diff{\abs{u}}{t}\leq \abs{\diff{u}{t}}$, it follows that $\diff{\abs{(q-q^*,p-q^*)}}{t}\leq \tilde{\ell}\abs{(q-q^*, p-q^*)}\quad\forall q^*\in\mathcal{A}_{NE},$
which, by the Gronwall-Bellman inequality implies that
$
    \abs{(q-q^*, p-q^*)}\leq \abs{(q(0)-q^*, p(0)-q^*)}e^{\tilde{\ell}t}\quad\forall q^*\in\mathcal{A}_{NE},
$
and all $(t,\cdot)\in\text{dom}((q,p)).$} Thus, the flow map \eqref{flowmap00} does not generate finite escape times. Moreover, since $\tau^+\in\{T_0,T\}^n$, we have that $G_1(D)\subset C_1\cup D_1$, which implies that solutions do not stop due to jumps leaving the set $C_1\cup D_1$. Next, note that the synchronization mechanism that governs the state $\tau$ is decoupled from the dynamics of the states $(q,p)$, and can be written as a HDS:
\begin{subequations}
\begin{align}\label{hybridsynchronization}
&\tau\in C_\tau\coloneqq[T_0,T]^n,~~~~~~~~~~~~~~~~~~~~~~~\dot{\tau}=\eta\mathbf{1}_n,\\
&\tau\in D_{\tau}:=\set{\tau\in C_\tau:\max_i\tau_i=T},~~\tau^+ \in G_{\tau}(\tau),
\end{align}
\end{subequations}
where $G_{\tau}(\tau)$ is the projection of $G_1$ into the $\tau$-component, which is independent of $(p,q)$. This hybrid system is well-posed by construction, and by \cite[Thm. 1]{Sync_Poveda} it renders uniformly fixed-time stable the set $\mathcal{A}_{\text{sync}}$, with a convergence bound $T^*$ given by $T^*\coloneqq \frac{1}{\eta}(T-T_0)+n$, $\forall~\tau(0,0)\in [T_0,T]^n$. Moreover, by \cite[Thm. 1]{Sync_Poveda}, each solution has at most $n$ jumps in any interval of length $L\coloneqq \frac{1}{\eta}(T-T_0)$, and, for any pair of hybrid times $(t,j), (s,i)\in\dom{\tau}$ with $t+j\ge s+i \ge T^*$ the following dwell-time condition holds $L + t-s \ge \left\lfloor\frac{j-i}{n}\right\rfloor L$, where $\lfloor\cdot\rfloor$ denotes the floor function. Thus, any solution $\tau$ of system \eqref{hybridsynchronization} is complete and also satisfies $|\tau(t,j)|_{\mathcal{A}_{\text{sync}}}=0$ for all $t+j\geq T^*$ such that $(t,j)\in\text{dom}(\tau)$.  Since the states $(q,p)$ of $\mathcal{H}_1$ evolve in $\R^n\times\R^n$, for each $\tau(0,0)\in [T_0,T]^n$ the hybrid time domains of system \eqref{jump_map}-\eqref{jump_set11} are the same hybrid time domains of system \eqref{hybridsynchronization}. This equivalence, plus the above properties, establish the result. \hfill $\blacksquare$

The previous Lemma directly implies the following:

\vspace{0.1cm}
\begin{lemma}\label{lemma:fixedsynchron}
Let $\nu>0$ and consider the HDS $\mathcal{H}_1$ with restricted flow and jump sets given by:
\begin{align*}
C_{\nu}\coloneqq &\Big\{x\in\R^{3n}:(p,q)\in \{(q^*,q^*)\}+\nu\mathbb{B},~\tau\in[T_0,T]^n \Big\},\\
D_{\nu}\coloneqq &\Big\{x\in\R^{3n}:~x\in C_{\nu},~~\max_{i\in\mathcal{V}}\tau_i=T\Big\},
\end{align*}
and jump map $G_1$ with values intersected with the set $C_{\nu}\cup D_{\nu}$. Then, the restricted system $\mathcal{H}_{\nu}=\{F_1,C_{\nu},G_{\nu},D_{\nu}\}$ renders UGFxS the set $\mathcal{A}_{\nu}\coloneqq (\{(q^*,q^*)\}+\nu\mathbb{B})\times \mathcal{A}_{\text{sync}}$. \QEDB
\end{lemma}
\vspace{0.1cm}

With Lemmas \ref{lemmawellposed} and \ref{lemma:fixedsynchron} at hand, we can proceed to analyze the HDS $\mathcal{H}_1$ by studying the HDS  $\mathcal{H}_v$ with data intersected with the set  $\mathcal{A}_{\nu}$. We denote this new restricted HDS as $\mathcal{H}_s\coloneqq \{F_s,C_s,G_s,D_s\}$, and we note that any compact set $\mathcal{A}'\subset\mathbb{R}^n\times\mathbb{R}^n$ such that $\mathcal{A}'\times\mathcal{A}_{\text{sync}}$ is UGAS for this system will also be UGAS for $\mathcal{H}_v$ thanks to the hybrid reduction principle \cite[Cor. 7.24]{HDS}. Moreover, since $\nu$ is arbitrary and independent of any parameter, and $\mathcal{H}_1$ has no finite escape times, the set $\mathcal{A}'\times\mathcal{A}_{\text{sync}}$ will also be UGAS for $\mathcal{H}_1$. Thus, in the following we focus on studying the stability properties $\mathcal{H}_s$.

\vspace{0.1cm}
\subsubsection{Proofs for \tcb{Potential-games}} For simplicity, we first present the proofs for potential games.
%
%

\vspace{0.05cm}
\begin{lemma}\label{lemma:convex:potential}
Under the conditions of Theorem \ref{theorem1}-(i$_1$), system $\mathcal{H}_s$ renders UGAS the set $\mathcal{A}$ given by \eqref{definitionofset}.  \QEDB
\end{lemma}

\vspace{0.1cm}\noindent
\textsl{Proof:} By using the potential function $P$ we define the error $\tilde{P}(q)\coloneqq P(q)-P(\mathcal{A}_{\text{NE}})$, and we consider the Lyapunov function
\begin{equation}\label{Lyapunov_function1}
V(x)=V_1(x)+V_2(x)+V_3(x),
\end{equation}
where the smooth functions $V_i$ are defined as follows:

\begin{subequations}\label{LyapunovfunctionsV}
\begin{align}
    \tcb{V_1(x)}&\tcb{\coloneqq \frac{1}{4}|p-q|^2,~V_2(x)\coloneqq \frac{1}{4}|p|_{\mathcal{A}_{NE}}^2},\label{V1V2}\\
V_3(x)&\coloneqq \frac{|\tau|^2}{n} \tilde{P}(q),
\end{align}
\end{subequations}
where $\abs{z}_{, \mathcal{A}_{NE}}^2=\min_{s\in\mathcal{A}_{NE}}\abs{z-s}^2$ and $\abs{z}^2 = z^\top z$.
By our definition of \tcb{potential-games}, and the construction of $V_1$ and $V_2$, the function $V$ is radially unbounded and positive definite with respect to the compact set $\mathcal{A}\cap \left(C_s\cup D_s\right)$. During flows in  $C_s$ the time derivative of $V$ satisfies:
%
%
\begin{align}\label{dotV1}
\dot{V}(x)&=\diffp{V(x)}{x}\dot{x}=\frac{\partial V(x)}{\partial q}\dot{q}+\frac{\partial V(x)}{\partial p}\dot{p}+\frac{\partial V(x)}{\partial \tau}\dot{\tau}.
\end{align}
The first term of \eqref{dotV1} is given by:
\begin{align*}
\frac{\partial V(x)}{\partial q}\dot{q}&=\left(-\frac{1}{2}(p-q)+\frac{|\tau|^2}{n}\tcb{\nabla P(q)}\right)^\top \dot{q}\notag\\
&=-(p-q)^\top D(\tau)^{-1}(p-q)\ldots\notag\\
&~~~+2\frac{|\tau|^2}{n}\tcb{\nabla P(q)}^\top D(\tau)^{-1}(p-q),\notag\\
&=-\frac{1}{\tau_s}(p-q)^\top (p-q)+2\tau_s (p-q)^\top \tcb{\nabla P(q)},
\end{align*}
where in the last equality we used the fact that in the set $C_s\backslash D_s$ we have $\tau=\tau_s\mathbf{1}_n$ with $\tau_s\in[T_0,T]$, and that points in $C_s \cap D_s$ lead to solutions that cannot flow. Similarly, 
\begin{align*}
\frac{\partial V(x)}{\partial p}\dot{p}&=\left(\frac{1}{2}(p-q)+\frac{1}{2}\big(p-\Pi_{\mathcal{A}_{NE}}(p)\big)\right)^\top \dot{p}\notag\\
&=-(p-q)^\top D(\tau)\tcb{\mathcal{G}(q)}\ldots\notag\\
&~~~-(p-\Pi_{\mathcal{A}_{NE}}(p))^\top  D(\tau)\tcb{\mathcal{G}(q)},\notag\\
&=-\tau_s\Big((p-q)+(p-\Pi_{\mathcal{A}_{NE}}(p))\Big)^\top \tcb{\mathcal{G}(q)},
\end{align*}
where $\Pi_{\mathcal{A}_{NE}}(p)$ is the projection of $p$ on $\mathcal{A}_{NE}$. Moreover, it follows that
\begin{align*}
\frac{\partial V(x)}{\partial \tau}\dot{\tau}&=\left(\tilde{P}(q) \frac{\tau}{n} \right)^\top \dot{\tau}=\tilde{P}(q)\frac{\tau^\top \mathbf{1}_n}{n}=2\tau_s \tilde{P}(q)\eta.
\end{align*}
Combining the above inequalities, and using the definition of $\tcb{\mathcal{G}(q) = \nabla P(q)}$ and $0<\eta\leq 1/2$, we obtain:
\begin{equation}\label{inequalityflowsLyapunov}
\dot{V}(x)\leq -\frac{1}{\tau_s}|p-q|^2-\tau_s\left(\big(q-\Pi_{\mathcal{A}_{NE}}(p)\big)^\top \tcb{\mathcal{G}}(q)-\tilde{P}(q)\right).
\end{equation}
\tcb{Using the Lipschitz property of $\mathcal{G}$ and the convexity of $\tilde{P}$ (implied by the monotonicity of $\tcb{\mathcal{G}}$ \cite[Thm. 12.17]{Rockafellar}), it follows that \cite[Thm. 5.8]{beck2017first}:
\begin{align*}
    \frac{1}{2\ell}\abs{\mathcal{G}(q)}^2 &\leq  \Big(q-\Pi_{\mathcal{A}_{NE}}(p) \Big)^\top\mathcal{G}(q) - \tilde{P}(q),
\end{align*}
and thus, from \eqref{inequalityflowsLyapunov}, we obtain during flows that
\begin{align}\label{decrease_flows}
\dot{V}(x)
&\leq -\frac{1}{\tau_s}|p-q|^2-\frac{\tau_s}{2\tcb{\ell}}|\tcb{\mathcal{G}}(q)|^2.
\end{align}
}
\noindent Since $|\tcb{\mathcal{G}}(q)|=0$ if and only if $q\in\mathcal{A}_{\text{NE}}$, during flows we have $\dot{V}(x)<0$ for all $x\in C_s\backslash\mathcal{A}$.

Similarly, during each jump the change of $V$, denoted $\Delta^{j+1}_{j} V(x)\coloneqq V(x(t,j+1))-V(x(t,j))$, satisfies $\Delta^{j+1}_j V(x)=\Delta_j^{j+1}V_3(x)$ and
\begin{align}
\Delta_j^{j+1}V_3(x)&=\frac{\tilde{P}(q)}{n}\Big(|\tau^+|^2-|\tau|^2 \Big)=\frac{\tilde{P}(q)}{n}\sum_{i=1}^n(\tau_i^{2+}-\tau_i^2),\notag\\
&=-\frac{\tilde{\varepsilon}}{n}\tilde{P}(q)
\leq0,\label{lastjumpsync}
\end{align}
 for some $\tilde{\varepsilon}>0$, where the last equality follows by the definition of $G_1$ in \eqref{jump_map} and the two following facts: first,  if $x\in D_{s}$, we have two possible cases for all players $i\in\mathcal{V}$: a) if $\tau_i=T_0$, then $\tau^+_i=T_0$; b) if $\tau_i=T$ then $\tau^+_i\in\{T_0,T\}$; second, if $x\in D_{s}$, we have that in each jump one and only one player $i$ satisfies $\tau_i=T$ and $\tau^+_i=T_0$. Therefore, since $T>T_0$ \tcb{there exists $\tilde{\varepsilon}>0$ such that $T_0^2-T^2=-\tilde{\varepsilon}$}, leading to \eqref{lastjumpsync}. This implies that $V$ does not increase during each reset triggered by a player. Given that the hybrid time domains of $\mathcal{H}_s$ are intervals of flow of duration $\frac{1}{\eta}(T-T_0)$, followed by $n$ consecutive jumps, we can use \eqref{lastjumpsync} $n$ times to obtain:
\begin{align*}
\Delta_j^{j+n} V(z)=\sum_{k=1}^n \Delta^{j+k}_{j+{k-1}} V(z)= -\varepsilon \tilde{P}(q)\leq 0,~~\forall~x\in D_{s}.
\end{align*}
By \cite[Prop. 3.27]{HDS}, the periodic strong decrease of $V$ during flows, and its non-increase during jumps, imply that  $\mathcal{H}_s$ renders UGAS the set $\mathcal{A}$. \hfill $\blacksquare$
\vspace{0.05cm}
\begin{lemma}\label{lemma:mixed:strong:potential}
Under the conditions of Theorem \ref{theorem1}-(i$_2$), system $\mathcal{H}_s$ renders UGES the set $\mathcal{A}$.
\end{lemma}

\vspace{0.05cm}
\noindent 
\textsl{Proof:} Let $V$ given by \eqref{Lyapunov_function1}, where $\mathcal{A}_{\text{NE}}=\{q^*\}$ is now a singleton by the strong monotonicity of the pseudogradient. During flows, we have \eqref{inequalityflowsLyapunov} with $\Pi_{\mathcal{A}_{NE}}(p)=q^*$, which using the strong monotonicity of $\tcb{\mathcal{G}}$ leads to
\begin{equation}\label{exponential_decrease_flows}
\dot{V}(x)\leq -\frac{1}{\tau_s}|p-q|^2-\tau_s\frac{\tcb{\kappa}}{2}|q-q^*|^2\leq -\tcb{\lambda} V(x)
\end{equation}
where we used the global Lipschitz property of $\tcb{\mathcal{G}}$, and the quadratic upper bound of \eqref{Lyapunov_function1}, with 
\begin{equation}
\tcb{\lambda}\coloneqq\frac{2}{3\Delta}\frac{\min\{1,0.25T_0T\tcb{\kappa}\}}{\max\{1,2T^2\ell\}}\approx\frac{1}{12T}\frac{1}{\sigma_r}\frac{1}{\sigma_\phi},
\end{equation}
where the approximation holds when $T_0$ is sufficiently small, and $T$ is sufficiently large (but finite). Thus, during each interval of flow, $V$ satisfies the $t$-time bound
\begin{equation}
    V(t,j)\leq V(t_j,j)e^{-\tcb{\lambda}(t-t_j)},\label{lyapunov:exponential:decrease:mixed:p:partial}
\end{equation}
%
for all $(t,j)\in \dom{x}$ such that $j=k\cdot n$ for some $k\in \N$.
To characterize the behavior of $V$ during jumps, let $\Theta$ and $I$ be the set of indices of players who implement $\alpha_i=0$, and  $\alpha_i=1$, respectively. It follows that after the $n$ consecutive jumps that proceed the intervals of flow of every solution of $\mathcal{H}_s$, we have:
\begin{align}
\hspace{-0.1cm}\Delta_j^{j+n} V(x)&\leq -\frac{1}{4}\left(|p-q|^2+|p-q^*|^2\right)\notag\\
&~~~+\frac{1}{4}\sum_{i\in I}\omega_i\left[(p_i-q_i)^2+(p_i-q_i^*)^2\right]\notag\\
&~~~+\frac{1}{4}\sum_{i\in\Theta}\omega_i[q_i-q_{i}^{*}]^2-\frac{\tcb{\kappa}}{2}(\tau_s^2-T_0^2),\notag\\
&\leq-\frac{1}{4}\sum_{i\in\Theta}\omega_i\left((p_i-q_i)^2+(p_i-q_i^*)^2\right)\ldots\notag\notag\\
&~~~-\frac{1}{2}\left(\tcb{\kappa}(\tau_s^2-T_0^2)-\frac{1}{2}\right)\sum_{i\in\Theta}(q_i-q_i^*)^2\ldots\notag\\
&~~~-\frac{\tcb{\kappa}}{2}(\tau^2_s-T_0^2)\sum_{i\in I}(q_i-q_i^*)^2\leq0,\label{deltaV:mixed:potential}
 \end{align}
where we used the fact that strong monotonicity of $\tcb{\mathcal{G}}$ implies that $\tilde{P}$ is strongly convex \cite[Thm. 12.17]{Rockafellar}, and the condition $T^2-T_0^2>\frac{1}{2\tcb{\kappa}}$ implied by (RC$_1$) with $\rho_J=\kappa^{-1}$. Therefore, it follows that $ V(t_j,j)\leq V(t_{j},j-n)e^{-\tcb{\lambda} (t-t_j)}$ for all $j\ge n$ and $t\ge t_j$. Since each interval of flow has length $L=(T-T_0)/\eta$, it follows that $ \tilde{V}(t_j,j)\leq \tilde{V}(t_{j-n} + L, j-n)e^{-\tcb{\lambda} L}e^{-\tcb{\lambda} (t-t_j)}$. Iterating, and using \eqref{lyapunov:exponential:decrease:mixed:p:partial}:
\begin{equation}\label{lyapunov:bound:p:mixed:thm}
    V(t,j)\leq V(0,0)e^{-\tcb{\lambda}\left(\lfloor \frac{j}{n}\rfloor-1\right)L}e^{-\tcb{\lambda}(t-t_j)}.
\end{equation}
By $\kappa$-strong convexity of $\tilde{P}$ and $\ell$-Lipschitz continuity of $\mathcal{G}_\omega$, $V$ satisfies:
\begin{equation}\label{v:quadraticBounds}
    \min\set{\frac{\tcb{1}}{4}, \frac{\kappa T_0^2}{2}}\abs{x}_\mathcal{A}^2\leq V(x) \leq \frac{1}{2}\max\set{\tcb{1} + \frac{\ell^2}{\kappa},~\frac{\tcb{3}}{2}}\abs{x}_\mathcal{A}^2,
\end{equation}
and thus, from \eqref{lyapunov:bound:p:mixed:thm}, we obtain:
\begin{equation}\label{state:bound:p:mixed}
    \abs{x(t,j)}_\mathcal{A} \leq \tcb{c} \abs{x(0,0)}_\mathcal{A}e^{-\frac{\tcb{\lambda}}{2}\left(t-(T-\max\tau(0,0))/\eta\right)},
\end{equation}
with $\tcb{c}>0$. 
Moreover, using the structure of the hybrid time domains, all hybrid times $(t,j)\in \dom{x}$ satisfy 
\begin{equation}\label{dwell:time:with:lambda}
    -\frac{\lambda}{2}t\leq -\frac{1}{3n}\lambda\left(t+j\right)+\frac{\lambda L}{3},
\end{equation}
for all $\lambda>0$. Hence, we obtain:
\begin{equation}
    \abs{x(t,j)}_\mathcal{A}\leq\tcb{\hat{c}_{s}}\abs{x(0,0)}_{\mathcal{A}}e^{-\frac{\tcb{\lambda}}{3n}(t+j)},\label{uges:p:mixed}
\end{equation}
where $\tcb{\hat{c}_{s}\coloneqq c} e^{\tcb{\lambda} L\left(\frac{1}{3} + \frac{1}{2\eta} \right)}$. Inequality \eqref{uges:p:mixed} implies that $\mathcal{A}$ is UGES under $\mathcal{H}_s$. \hfill $\blacksquare$

\vspace{0.1cm}
\begin{lemma}\label{lemma:homogeneous:strong:potential}
Under the conditions of Theorem \ref{theorem1}-(i$_3$), system $\mathcal{H}_s$ renders UGES the set $\mathcal{A}$.
\end{lemma}

\vspace{0.05cm}\noindent 
\textsl{Proof:} Using the Lyapunov function $V$ given by \eqref{Lyapunov_function1}, and the fact that $\mathcal{A}_{\text{NE}}=\{q^*\}$, we obtain again inequality \eqref{exponential_decrease_flows} during flows. Since now $\alpha= \vec{0}_n$, during jumps we have
\begin{align}
\Delta_j^{j+n} V(x)&\leq -V_1(x)-V_2(x)-\frac{1}{2}(\tau_s^2-T_0^2) \tilde{P}(q)\ldots\notag\\
&~~~+\frac{1}{4}|q-q^*|^2.\label{triggeringN}
\end{align}
By strong convexity of $\tilde{P}$, we can further bound \eqref{triggeringN} as:
%
%
\begin{align}
\Delta_j^{j+n} V(x)&\leq -V_1(x)-V_2(x)-\gamma\left(\kappa^{-1}\right) V_3(x)\notag\\
&\leq -\gamma \left(\kappa^{-1}\right) V(x),\label{jumps_exponetial}
\end{align}
where $\gamma(\cdot)$ is given by \eqref{gamma_parameterized}, which under (RC$_1$) satisfies $\gamma\left(\kappa^{-1}\right)\in(0,1)$. Thus, by \cite[Thm. 1]{LyapunovExponentialHDS}, inequalities \eqref{exponential_decrease_flows} and \eqref{jumps_exponetial}, and the quadratic upper and lower bounds of $V$, we obtain that $\mathcal{H}_s$ renders UGES the set $\mathcal{A}$. \hfill $\blacksquare$

With Lemmas \ref{lemma:convex:potential}-\ref{lemma:homogeneous:strong:potential} at hand for system $\mathcal{H}_s$, we can now proceed to proof the three main items of Theorem \ref{theorem1} for the original hybrid dynamics $\mathcal{H}_1$.

\vspace{0.1cm}\noindent 
\textbf{Proof of Theorem \ref{theorem1}}: \emph{(a) Stability:} By the hybrid reduction principle \cite[Cor. 7.24]{HDS}, UGAS of $\mathcal{A}$ for system $\mathcal{H}_s$ (established in Lemmas \ref{lemma:convex:potential}, \ref{lemma:mixed:strong:potential} and \ref{lemma:homogeneous:strong:potential}), and UGFxS of $\mathcal{A}_\nu$ for system $\mathcal{H}_{\nu}$, imply that $\mathcal{A}$ is UGAS for system $\mathcal{H}_{\nu}$. Moreover, since the choice of $\nu>0$ is arbitrary, and has no effect on the dynamics of the system, and since the trajectories of the original HDS $\mathcal{H}_1$ are complete and bounded, the compact set $\mathcal{A}$ is also UGAS for system $\mathcal{H}_1$. This establishes UGAS of $\mathcal{A}$ under the conditions of items (i$_1$), (i$_2$) and (i$_3$). For items (i$_2$) and (i$_3$), UGES follows by the exponential convergence bounds of Lemmas \ref{lemma:mixed:strong:potential}-\ref{lemma:homogeneous:strong:potential} and the fixed-time synchronization of $\tau$. R-UGAS and R-UGES follow directly by robustness results of well-posed HDS, specifically by \cite[Thm. 7.21]{HDS}.\\

\vspace{0.02cm}\noindent 
\emph{(b) Convergence Bounds:} For any solution $x$ and all $(t,j)\in\mathcal{T}(x)$ we have that  $|\tau(t,j)|_{\mathcal{A}_{\text{sync}}}=0$. Thus, for such times the trajectories of $\mathcal{H}_1$ satisfy the Lyapunov inequalities established in Lemmas \ref{lemma:convex:potential}-\ref{lemma:homogeneous:strong:potential}. To establish \eqref{decrease_potential}, we use inequality \eqref{decrease_flows}, which implies that for each $(t,j),(s,j)\in\mathcal{T}(x)$, such that $t>s$, we have $V(t,j)\leq V(s,j)$. Since $V_3\leq V$, and using $s_j\coloneqq \min~\{t\in\R_{\geq0},~(t,j)\in\mathcal{T}(x)\}$, we obtain
 \begin{equation}\label{bound_potentialgames1}
\tilde{P}(q(t,j))\leq \frac{2n}{\tau^\top \tau}V(s_j,j)=\frac{c_j}{\tau_s^2},~~\forall t~>s_j,
\end{equation}
where $c_j\coloneqq 2 V(s_j,j)$. Using the fact $V$ is non-increasing during flows and jumps, and also converges to zero, we obtain that $\{c_j\}_{j=0}^{\infty}\searrow0^+$. 
To obtain the convergence bound of item (i$_2$), we first note that from the proof of Lemma \ref{lemmawellposed} it follows that $\diff{\abs{x}_\mathcal{A}}{t}\leq \tilde{\ell}\abs{x}_\mathcal{A}$ for all $(t,j)\in \dom{x}$, where $\tilde{\ell}= 2\sqrt{2}\max\set{\frac{1}{T_0},T\ell}$. In particular, this implies that
\begin{equation}\label{potential:bound:tosync}
    \abs{x(t_s,j_s)}_\mathcal{A}\leq \abs{x(0,0)}_\mathcal{A}e^{\tilde{l}(T-\max\tau(0,0))/\eta},
\end{equation}
where $t_s,j_s$ are the smallest times for which $\abs{\tau(t,j)}_{\mathcal{A}_{sync}}=0$ for all $t+j\ge t_s + j_s$. Note that $x(t_s,j_s)\in C_s\cup D_s$, and hence \eqref{state:bound:p:mixed} holds with $\abs{x(0,0)}_{\mathcal{A}}$ replaced by $\abs{x(t_s,j_s)}_{\mathcal{A}}$, i.e., $\abs{x(t,j)}_{\mathcal{A}}$ satisfies:
\begin{equation}\label{potential:exponential:auxiliar}
        \abs{x(t,j)}_\mathcal{A} \leq \tcb{c} \abs{x(t_s,j_s)}_\mathcal{A}e^{-\frac{\tcb{\lambda}}{2}\left(t-(T-\max\tau(0,0))/\eta\right)},
\end{equation}
for all $t+j\ge t_s + j_s$. Therefore, inequality \eqref{potential:exponential:auxiliar}, together with \eqref{potential:bound:tosync}, implies:
\begin{equation}\label{state:bound:p:texp}
    \abs{x(t,j)}_\mathcal{A}\leq \tcb{\tilde{c}}\abs{x(0,0)}_{\mathcal{A}}e^{-\frac{\tcb{\lambda}}{2}t},
\end{equation}
for all $(t,j)\dom{x}$, and where $\tcb{\tilde{c}}=\tcb{c} e^{\left(\frac{\tcb{\lambda}}{2}+\tilde{l}\right)L}$. Using \eqref{dwell:time:with:lambda} and \eqref{state:bound:p:texp} gives: 
\begin{equation}\label{exponential:decrease:mixed:bound:potential}
    \abs{x(t,j)}_\mathcal{A}\leq\hat{c}_\omega\abs{x(0,0)}_{\mathcal{A}}e^{-\frac{\tcb{\lambda}}{3n}(t+j)},
\end{equation}
with $\hat{c}_\omega=\tcb{\tilde{c}} e^{\tcb{\lambda} L /3}$ which establishes the bound in \eqref{exp_property}. This also implies that $\mathcal{H}_1$ renders $\mathcal{A}$ UGES under the conditions of Theorem \ref{theorem1}-(i$_2$). Finally, to establish the convergence bound of item (i$_3$), we note that \eqref{jumps_exponetial} implies $V(x(t,j+n))\leq (1-\gamma\left(\kappa^{-1}\right))V_3(x(t,j))$. Since $V_3(x)\leq V(x)$ for all $(t,j)\in\mathcal{T}(x)$, $V$ does not increase during flows, and using the periodicity of the hybrid time domains, we obtain:
\begin{equation}\label{auxiliaryboundexponential}
V_3(t,j_s+kn)\leq (1-\gamma\left(\kappa^{-1}\right))^k V_3(t_s,j_s),~~ \forall~k\in\mathbb{Z}_{\geq0},
\end{equation}
for all $t\in(t_s+(k-1)L,t_s+kL)$, where $(t_s,j_s)$ denotes the first hybrid time after which the timers $\tau$ flow synchronized. By Lemma \ref{lemmawellposed}, such times are uniformly bounded as $0\leq t_s+j_s\leq 2T^*$. Using \eqref{auxiliaryboundexponential}, the definition of $V_3$, as well as strong convexity and smoothness of $\tilde{P}$, we obtain:
\begin{equation}\label{exponentialdecreasejumpsA}
|q(t,j_s+kn)-q^*|\leq \sigma_r\sqrt{\frac{\tcb{\ell}}{\tcb{\kappa}}}(1-\gamma\left(\kappa^{-1}\right))^{\frac{k}{2}} |q(t_s,j_s)-q^*|,
\end{equation}
for all $k\in \mathbb{Z}_{\ge 0}$.
Finally, since by Lemma \ref{lemmawellposed} all solutions are bounded, it follows that for each compact set of initial conditions $K_0$ there exists $M_0>0$ such that $|x(t,j)|_{\mathcal{A}}\leq  M_0$ for all $(t,j)\in\text{dom}(x)$ such that $0\leq t \leq t_s$ and $0\leq j\leq j_s$. This bound, combined with \eqref{exponentialdecreasejumpsA}, implies the bound of the theorem via the change of variable $j=j_s+kn$ and the upper bound $n\leq j_s\leq 2n$.  \hfill $\blacksquare$
\vspace{0.1cm}
\subsubsection{Proofs for Non-Potential Games} As before, we divide the proof in different lemmas.

\vspace{0.1cm}\noindent 
\begin{lemma}\label{lemma:strict:np}
Consider the HDS $\mathcal{H}_s$ under the Assumptions of Theorem \ref{theoremstrictmono}. Then, the set $\mathcal{A}$ is UGAS. \QEDB
\end{lemma}

\vspace{0.1cm}\noindent 
\textsl{Proof:} First, by Assumption \ref{regularassumption} and strict convexity of $\phi_i$ on $q_i$  for all $i\in\mathcal{V}$, which follows by strict monotonicity of the pseudo-gradient, existence of the NE is guaranteed via \cite[Cor 4.2]{BasarDNGT}. With this at hand, we consider the Lyapunov function \tcb{ $\tilde{V}=V_1+V_2+\tilde{V}_3$, where $V_1$ and $V_2$ are defined in \eqref{LyapunovfunctionsV}, and $\tilde{V}_3$ is given by}:
\begin{equation}\label{lyapunov2}
\tilde{V}_3(x)\coloneqq  c_o\frac{\abs{\tau}^2\abs{\mathcal{G}(q)}^2}{2n},
\end{equation}
where $c_o$ corresponds to the cocoercivity constant of $\mathcal{G}$.
By construction and Assumption \ref{cocoerciveassumption}, we have that $V$ is radially unbounded,  and also positive definite with respect to $\mathcal{A}\cap (C_s\cup D_s)$. Using $c_o$-cocoercivity of $\mathcal{G}$, inequality \eqref{inequalityflowsLyapunov} now becomes
\begin{align*}
\dot{\tilde{V}}(x)&\leq -\frac{1}{\tau_s}\abs{p-q}^2-2 \tau_s(p-q)^\top \mathcal{G}(q)\notag\\
&~~~+2c_o\tau_s \mathcal{G}^\top(q)\diffp{\mathcal{G}}{q}(q)(p-q)+\tcb{ c_o(1-\eta)\tau_s\abs{\mathcal{G}(q)}^2},
\end{align*}
which can be written in quadratic form as
\begin{align}\label{quadraticboundLyapunov1}
  \dot{\tilde{V}}(x)&\leq  -\tau_s \tilde{x}^\top M_{1/c_o}(q,\tau_s)\tilde{x},
\end{align}
with $\tilde{x}\coloneqq \big((p-q), \mathcal{G}(q)\big)$, and 
\begin{align*}
    M_{1/c_o}(q,\tau_s)\coloneqq  \begin{pmatrix}
    \frac{1}{\tau_s^2}I_n &  I_n -c_o\partial \mathcal{G}(q)^\top\\
    I_n -c_o\partial \mathcal{G}(q) &c_o(1-\eta)I_n \tageq{\label{M1}}
\end{pmatrix}.
\end{align*}
Since $\eta\leq \frac{1}{2}$ by design, $c_o=1/\ell$, and $\tau_s\in[T_0,T]$, under the conditions of Theorem \ref{theoremstrictmono}, we have that  $M_{\ell}(q,\tau_s)\succ 0$ for all $\tau_s\in[T_0,T]$ and $q\neq q^*$ whenever
\begin{equation}\label{inequaequivalent11}
    0\prec I_n -\frac{T^2}{\ell(1-\eta)}\big(\ell I_n - \partial  \mathcal{G}(q)^\top \big)\big(\ell I_n -\partial \mathcal{G}(q)\big).
\end{equation}
The expression in \eqref{inequaequivalent11} is precisely \eqref{contractivematrix} with $\rho_F=\ell$ and $\delta=0$. Thus, since by assumption $\mathcal{S}_0$ is $\ell$-GC, it follows that \eqref{inequaequivalent11} holds. Finally, note that when $q=q^*$ inequality \eqref{quadraticboundLyapunov1} reduces to $\dot{\tilde{V}}(x)\leq -\frac{1}{\tau_s}|p-q|^2\leq0$.


On the other hand, after the $n$ consecutive jumps that proceed each interval of flow, the change of $V$ is
\begin{equation}\label{jumpsinequalitystrict}
\Delta_j^{j+n} \tilde{V}(z)= \frac{c_o}{2}\abs{\mathcal{G}(q)}^2\left(T_0^2 - T^2\right)\leq0.
\end{equation}
Finally, we show that no complete solution $x$ of $\mathcal{H}_s$ keeps $\tilde{V}$ in a non-zero level set. In particular, since for all $(q,p,\tau)\in\mathbb{R}^n\backslash\{q^*\}\times\mathbb{R}^n\times[T_0,T]$ we have that $\dot{V}<0$, it suffices to consider the case $q=q^*$, which leads to $\dot{V}=0$ only when $p=q$, i.e., when $(p,q)\in\mathcal{A}$.
Therefore, no solution that flows can keep $\tilde{V}$ constant in a non-zero level set. Since the flows are periodic UGAS of $\mathcal{A}$ follows now by the hybrid invariance principle \cite[Thm. 8.8]{HDS}.   \hfill $\blacksquare$

\vspace{0.2cm}\noindent 
\textbf{Proof of Theorem \ref{theoremstrictmono}}: \emph{(a) Stability Properties:} Follows by the same ideas used in the proof of the stability properties of Theorem \ref{theorem1}-(i$_1$), but using Lemma \ref{lemma:strict:np} instead of Lemma \ref{lemma:convex:potential}.

\emph{(b) Convergence Bounds:}  Follows by the same steps used in the proof of Theorem 1-(i$_1$), substituting \eqref{bound_potentialgames1} by
%
 \begin{equation}\label{bound_np_games1}
\abs{\mathcal{G}(q)}^2\leq \frac{2\ell n}{\tau^\top \tau}\tilde{V}_3(s_j,j)=\frac{\tilde{c}_j}{\tau_s^2}.
\end{equation}
where $\tilde{c}_j\coloneqq 2\ell \tilde{V}_3(s_j,j)$. \hfill$\blacksquare$
%
\vspace{0.1cm}
\begin{lemma}\label{lemma:mixed:np}
Consider the HDS $\mathcal{H}_s$ under the Assumptions of Theorem \ref{theorem:strong:np}-(i$_4$). Then, the set $\mathcal{A}$ is UGES. \QEDB
\end{lemma}

\vspace{0.1cm}
\noindent
\textsl{Proof:} We consider the Lyapunov function $V$ used in the proof of Lemma \ref{lemma:strict:np}, with $\tilde{V}_3$ given by \eqref{lyapunov2} and  $c_o=\kappa/\ell^2$. 
The time derivative of $\tilde{V}$ now satisfies
\begin{align*}
  \dot{\tilde{V}}(x)&\leq  -\tau_s \tilde{x}^\top M_{\sigma_\phi\ell}(q,\tau_s)\tilde{x},
\end{align*}
with $\tilde{x}\coloneqq \big((p-q), \mathcal{G}(q)\big)$.
By assumption we know that $\mathcal{S}_\delta$ is $(\sigma_\phi\ell)$-GC, which is equivalent to:
\begin{align*}
    0\prec I_n - \left(\frac{T^2}{1-T^2\delta}\right)\frac{\left(\sigma_\phi\ell I_n-\partial\mathcal{G}(q)^\top\right)\left(\sigma_\phi\ell I_n-\partial\mathcal{G}(q)\right)}{\sigma_\phi\ell(1-\eta)-\sigma_\phi^2\ell^2\delta}.
\end{align*}
In turn, when $0<\delta <(1-\eta)/\sigma_\phi\ell$ and $0<\eta\leq 1/2$, the above inequality directly implies that $M_{\sigma_\phi\ell}(q,\tau_s)\succ \delta I_n$, for all $\tau_s\in[T_0,T]$ and all $q\neq q^*$. Thus, for such points, and during flows, we have  $\dot{\tilde{V}}\leq -\delta(\abs{p-q}^2 + \abs{\mathcal{G}}(q))$. Using $\kappa$-strong-monotonicity and $\kappa/\ell^2$-cocoercivity of $\mathcal{G}$ we conclude
\begin{align}\label{dotV:uges:strong:np}
    \dot{\tilde{V}}(x)\leq -\lambda\tilde{V}(x),\text{ with }\lambda = \frac{4\delta}{\max\set{3,2(\frac{1}{\kappa^2}+\frac{\kappa}{\ell^2}T^2)}}.
\end{align}


On the other hand, during jumps, using (RC$_1$), the definition of $\tilde{V}_3$, and the Reset Policy $\alpha\in\{0,1\}^n$, the change of $\tilde{V}$ is
%

\vspace{-0.1cm}
\begin{small}
\begin{align}\label{deltaV:mixed:np}
\Delta_j^{j+n} \tilde{V}\leq -\frac{1}{4}\sum_{i\in \Theta}\left[(p_i-q_i)^2+(p_i-q_i^*)^2\right]-\gamma(\sigma_\phi^2\kappa^{-1}) \tilde{V}_3(x),
\end{align}
\end{small}

\vspace{-0.1cm}\noindent 
where $\gamma(\sigma_\phi^2\kappa^{-1})\in(0,1)$ is given by \eqref{gamma_parameterized}, and $\Theta$ is defined in the proof of Lemma \ref{lemma:mixed:strong:potential}. Thus, it follows that $\Delta_j^{j+n} \tilde{V}\leq 0$. Moreover, by the $\kappa$-strong monotonicity and $\ell$-Lipschitz continuity of $\mathcal{G}$, $\tilde{V}$ satisfies the quadratic bounds $\underline{c}\abs{x}_\mathcal{A}^2\leq\tilde{V}(x)\leq \overline{c}\abs{x}_{\mathcal{A}}^2 $, where:
\begin{equation}\label{vtilde:quadraticBounds}
\underline{c} \coloneqq \min \set{\frac{1}{4},~\frac{\kappa T_0^2}{2\sigma_a^2}},~\overline{c}\coloneqq \max\set{\frac{3}{4},~\frac{1}{2}+ \frac{\kappa T^2 \ell^2}{2}}.
\end{equation}
The exponential decrease of $V$ during the flows (which are periodic), the non-increase of $V$ during the jumps, and the quadratic upper and lower bounds of $\tilde{V}$, imply that $\mathcal{H}_s$ renders UGES the set $\mathcal{A}$.
\hfill $\blacksquare$

\begin{lemma}\label{lemma:strong:np}
Consider the HDS $\mathcal{H}_s$ under the Assumptions of Theorem \ref{theorem:strong:np}-(i$_5$). Then, the set $\mathcal{A}$ is UGES. \QEDB
\end{lemma}

\vspace{0.1cm}\noindent 
\textsl{Proof:} Consider the Lyapunov function $\tilde{V}$ used in the proof of Lemma \ref{lemma:strict:np}, which still satisfies \eqref{dotV:uges:strong:np}. 
During jumps, the reset policy $\alpha=\mathbf{0}_n$ implies that $\Theta = \mathcal{V}$ in \eqref{deltaV:mixed:np}, leading to $\Delta_j^{j+n} \tilde{V}(x)\leq -V_1(x)-V_2(x)-\gamma(\sigma_\phi^2\kappa^{-1})V_3(x)\leq -\gamma(\sigma_\phi^2\kappa^{-1}) \tilde{V}(x)$.  Therefore, by \cite[Thm. 1]{LyapunovExponentialHDS}, and the quadratic upper and lower bounds of $\tilde{V}$, 
system $\mathcal{H}_s$ renders UGES the set $\mathcal{A}$. \hfill $\blacksquare$

\vspace{0.2cm}\noindent 
\textbf{Proof of Theorem \ref{theorem:strong:np}}: \emph{(a) Stability Properties:} Follows by using using Lemmas \ref{lemma:mixed:np} and \ref{lemma:strong:np} in conjunction with the same ideas used in the proof of Theorem 1.

\emph{(b) Convergence Bounds:} We follow the same steps of the proof of Theorem 1, using now $\tilde{V}_3$ instead of $V_3$. For item (i$_4$), this leads to the following bound instead of \eqref{exponential:decrease:mixed:bound:potential}:
\begin{equation*}
     \abs{x(t,j)}_\mathcal{A}\leq\hat{c}\abs{x(0,0)}_{\mathcal{A}}e^{-\frac{\lambda}{3n}(t+j)},
\end{equation*}
where $\lambda$ are defined in \eqref{dotV:uges:strong:np}, $\hat{c}\coloneqq \sqrt{\overline{c}/\underline{c}}\cdot e^{\left(\frac{5}{6}\lambda+\tilde{l}\right)L}$, and $\underline{c}$ and $\overline{c}$ are given in \eqref{vtilde:quadraticBounds}. Finally, for item (i$_5$), we obtain the following bound instead of \eqref{exponentialdecreasejumpsA}:
%
%
\begin{equation*}
|q(t,j_s+kn)-q^*|\leq \sigma_r\sigma_\phi\left(1-\gamma\left(\sigma_\phi^2\kappa^{-1}\right)\right)^{\frac{k}{2}} |q(t_s,j_s)-q^*|, 
\end{equation*}
from here, the proof follows the exact same steps.
\hfill$\blacksquare$

\vspace{0.1cm}\noindent 
\textbf{Proof of Lemmas \ref{sufficientlemma} and \ref{lemma_strong_suff1}:} We first show Lemma \ref{lemma_strong_suff1}. Using $c_o=\kappa/\ell^2$ we have that (RC$_3$) can be equivalently written as 
%
%
\begin{equation}\label{T:bound:c0:preMatrix}
\tilde{\alpha} > 1 - 2c_o \kappa + c_o^2\ell^2,
\end{equation}
with $\tilde{\alpha} \coloneqq\left(\frac{1}{T^2}-\delta\right)(c_o(1-\eta)-\delta).$
By using the fact that $\mathcal{G}$ is $\ell$-Lipschitz continuous, and hence that $\partial \mathcal{G}(q)^\top \partial \mathcal{G}(q) \prec \ell^2I_n$ \cite{gadjov2022exact}, and the monotonicity properties of the pseudogradient, which imply that 
$\partial\mathcal{G}(q)+\partial\mathcal{G}(q)^\top\succ 2\kappa I_n$ \cite[Prop 2.3.2 c)]{facchinei2007finite}, from \eqref{T:bound:c0:preMatrix} it follows that
\begin{align*}
0 &\prec (\alpha -1) I +  c_o\left(\partial \mathcal{G}(q) + \partial \mathcal{G}(q)^\top\right)  -c_o^2 \partial \mathcal{G}(q)^\top \partial \mathcal{G}(q).
\end{align*}
and hence that 
\begin{align*}
    0\prec I - \left(\frac{T^2}{1-T^2\delta}\right)\frac{\left(I_n-c_o\partial \mathcal{G}(q)^\top\right)\left(I_n-c_o\partial \mathcal{G}(q)\right)}{c_o(1-\eta) - \delta},
\end{align*}
which implies, whenever $0\leq \delta <c_0(1-\eta)$, that $\mathcal{S}_\delta$ is $(1/c_o)$-GC. Lemma 2 follows by the same arguments, using $c_o=1/\ell$ and letting $\kappa\to 0^+$.  
\hfill $\blacksquare$

\vspace{0.1cm}\noindent 
\textbf{Proof of Lemma \ref{sufficient_condition}:} To satisfy (RC)$_1$ and (RC)$_3$ we need 
\begin{align*}
 T_0^2 + \frac{\ell^2}{2\kappa^3} < \frac{\kappa(1-\eta)-\delta\ell^2}{\ell^2-\kappa^2+\delta(\kappa(1-\eta)-\delta\ell^2)}.
\end{align*}
Since $T_0^2>$ 0, it is necessary that
\begin{align*}
    &\frac{\ell^2}{2\kappa^3} < \frac{\kappa(1-\eta)-\delta\ell^2}{\ell^2-\kappa^2+\delta(\kappa(1-\eta)-\delta\ell^2)},
\end{align*}
which in turn is equivalent to 
\begin{align*}
        1  < 2\frac{(1-\eta)}{\sigma_\phi^4} + \frac{1}{\sigma_\phi^2}\left(1 -\delta\left(2\kappa + \frac{1-\eta}{\kappa} -\delta\frac{\ell^2}{\kappa^2}\right)\right).\tageq{\label{necessary:joint}}
\end{align*}
Since by assumption $\sigma_\phi^4-\sigma_\phi^2<2(1-\eta)$, there exists $\delta>0$ sufficiently small such that \eqref{necessary:joint} holds.  \hfill$\blacksquare$

\vspace{0.1cm}\noindent 
\textbf{Proof of Lemma \ref{quasi-optimal-restaring:np}: }
The convergence bound of Theorem \ref{theorem:strong:np}-(i$_5$) implies the slightly looser bound $\forall~(t,j)\in\text{dom}(x)$:
\begin{equation}\label{exponentialbound:np:simplified}
    |q(t,j)-q^*|\leq \sigma_r\sigma_\phi\left(1-\gamma\left(\rho_J\right)\right)^{\frac{1}{2}\left\lfloor\frac{j-n}{n}\right\rfloor}M_0.
\end{equation}
Since $(0,0)\in \dom{x}$, by using the structure of the hybrid time domain, with $(i,s)=(0,0)$, we obtain
\begin{equation*}
\left\lfloor\frac{j-n}{n}\right\rfloor \leq \frac{t}{L},
\end{equation*}
where $L=(T-T_0)/\eta$. Hence, using \eqref{exponentialbound:np:simplified}, $\abs{q-q^*}$ satisfies
\begin{equation}\label{exponentialbound:np:simplified2}
    \abs{q(t,j)-q^*}\leq  \sigma_\phi\sigma_r \left(1-\gamma\left(\sigma_\phi^2\kappa^{-1}\right)\right)^{\frac{\eta}{2}\frac{t}{T-T_0}}M_0.
\end{equation}
By minimizing the right hand side of \eqref{exponentialbound:np:simplified2} with respect to the restarting parameter $T$, we find that the minimum value is achieved when $(1-\gamma(\sigma_\phi\kappa^{-1})) = \frac{1}{e^2}$. Solving for $T$, we obtain precisely  $T^{\text{opt}}$. Using this restarting parameter, for any $\nu>0$ the error $\abs{q(t,j)-q^*}\leq \nu$ is obtained when $\sigma_\phi\sigma_r e^{-\eta\frac{t}{T^{\text{opt}}-T_0}}M_0 \leq \nu$.
%
%
This inequality holds precisely for all $t\geq t_\nu^{\text{opt}}$.
Finally, using $(1-\gamma(\sigma_\phi\kappa^{-1})) = \frac{1}{e^2}$ in \eqref{exponentialbound:np:simplified2}: 
%
\begin{equation}\label{optimal:convergence}
    \abs{q(t,j)-q^*}\leq  \sigma_\phi\sigma_r e^{-\eta\frac{t}{T^{\text{opt}}-T_0}}M_0.
\end{equation}
As $T_0\to 0^+$, and using the value of $T^{\text{opt}}$:
\begin{equation}
    \abs{q(t,j)-q^*}\leq  \sigma_\phi\sigma_r e^{-t\frac{\eta\sqrt{2\kappa}}{e\sigma_\phi}}M_0,
\end{equation}
which gives the convergence bound of order $\mathcal{O}(e^{-\sqrt{\kappa}/\sigma_\phi})$.

\vspace{0.1cm}\noindent 
\textbf{Proof of Proposition \ref{propositionquadratic}:} We consider the Lyapunov function $\hat{V}=V_1+V_2+\hat{V}_3$, where $V_1$ and $V_2$ are given by \eqref{V1V2}, and $\hat{V}_3$ is defined as $\hat{V}_3=\kappa\frac{|\tau|^2}{2n}|q-q^*|^2$. During flows we now have $\dot{\tilde{V}}(x)\leq  -\tau_s \hat{x}^\top M_Q(\tau_s)\hat{x}$,
%
%
%
where $\hat{x}=(p-q,q-q^*)$ with $M_Q$ given by
\begin{equation}
M_Q(\tau_s)=\begin{pmatrix}
        \frac{1}{\tau_s^2}I_n & \left( A-\kappa I_n\right)\\
        \left( A^\top-\kappa I_n\right) & \kappa(1-\eta) I_n
    \end{pmatrix}.
\end{equation}
This matrix is positive definite if and only if the following matrix inequality holds for all $\tau_s\in[T_0,T]$:
\begin{equation}
0 \prec \frac{\kappa(1-\eta)}{\tau_s^2}I_n -
    \left(\kappa  I_n-A\right)\left(\kappa I_n-A^\top\right).
\end{equation}
Hence, it suffices to verify the condition $0 \prec I_n -\frac{T^2}{\kappa(1-\eta)}
    \left(\kappa  I_n-A\right)\left(\kappa I_n-A^\top\right),$
%
which is equivalent to $\kappa$-GC of $\mathcal{S}_0$. It follows that $M_Q(\tau_s)\succ 0$ for all $\tau\in[T_0,T]$, and $\dot{\hat{V}}\leq -\tau_s c |\hat{z}|^2$ during flows, for some $c>0$. On the other hand, during jumps, the policy $\alpha=\mathbf{0}_n$ leads to:
\begin{align*}
\Delta_j^{j+n} V(z)&\leq-V_1(z)-V_2(z)\left(\frac{\kappa}{2}(\tau_s^2-T_0^2)-\frac{1}{4\kappa^2}\right)|q-q^*|^2.
\end{align*}
Using (RC$_2$), we obtain  $\Delta_j^{j+n} V(z)\leq-V_1(z)-V_2(z)-\gamma(\kappa^{-1}) V_3(z)$, with $\gamma$ as in \eqref{gamma_parameterized}.  UGES of $\mathcal{A}$ follows by the same arguments of the proof of Theorem 3-(i$_5$).
\hfill $\blacksquare$

\vspace{-0.2cm}
\subsection{Proofs of Section 3}
%
To prove Theorem \ref{theorem:partialInformation}, we present two auxiliary lemmas:
\begin{lemma}\label{lemma:consensus:UGASsynchro}
Consider the assumptions of Theorem \ref{theorem:partialInformation}, and let $\mathcal{H}_{2,s}=\{C_{2,s},F_{2,s},D_{2,s},G_{2,s}\}$ be obtained by intersecting the data of $\mathcal{H}_2$ with $\mathcal{A}_{2,\nu} \coloneqq \mathcal{A}_\nu\times (\mathcal{Q}(\vec{1}_n\otimes q^*)+\nu\mathbb{B})$, where $\mathcal{A}_\nu = (\{(q^*,q^*)\}+\nu\mathbb{B})\times \mathcal{A}_{\text{sync}}$. Then $\mathcal{H}_{2,s}$ renders UGAS the set $\mathcal{A}\times\{\mathcal{Q}(\vec{1}_n\otimes q^*)\}$.
\end{lemma}
\vspace{0.1cm}
\noindent
\textsl{Proof:} 
Consider the change of variable $\theta=\hat{q}-h(q)$, with $h(q)\coloneqq \mathcal{Q}(\mathbf{1}_n\otimes q)$ and let
\begin{align}
W(q,\theta,\varepsilon)&\coloneqq -\mathcal{Q}\vec{L}\mathcal{Q}^\top \theta  - \varepsilon\mathcal{Q}\left(\vec{1}_n\otimes 2D(\tau)^{-1}(p-q)\right)\notag\\
&~~~~+\varepsilon\mathcal{Q}\left(\vec{1}_n\otimes\mathcal{P}\vec{L}\mathcal{Q}^\top \theta\right).\label{flowEstimationTheta}
\end{align}
This change of coordinates leads to a HDS $\mathcal{H}_\vartheta$ with state $\vartheta\coloneqq(x,\theta)$, where $x=(q,p,\tau)$, and data $\mathcal{H}_\vartheta= (C_{2,\vartheta}, F_{2,\vartheta}, D_{2,\vartheta}, G_{2,\vartheta})$, where $C_{2,\vartheta},D_{2,\vartheta}$ and $G_{2,\vartheta}$ are obtained directly from \eqref{flowSetEstimation}, \eqref{jumpSetEstimation}, and \eqref{jumpEstimation} respectively via the change of coordinates, and where the 
flow map is defined by $F_{2,\vartheta}(\vartheta)\coloneqq (U(x,\theta + h(q)),~W(q,\theta,\varepsilon)/\varepsilon)$ where:
\begin{align}
U(x,\theta + h(q)) &= \begin{pmatrix}
    2D(\tau)^{-1}(p-q) -\mathcal{P}\vec{L}\mathcal{Q}^\top \theta\\
    -2D(\tau)\hat{\mathcal{G}}(\vec{1}_n\otimes q + \mathcal{Q}^\top \theta)\\
    \eta \vec{1}_n
\end{pmatrix}.\label{flowEstimationWoTheta}
\end{align}
Let $\mathcal{H}_{\vartheta,s}$ be the HDS that results from intersecting the data of $\mathcal{H}_{\vartheta}$ with $\mathcal{A}_\nu \times (\nu\mathbb{B})$, with $\nu>0$. Note that studying the stability of $\mathcal{A}\times\set{\mathcal{Q}(\vec{1}_n\otimes q^*}$ under $\mathcal{H}_{2,s}$, is is equivalent to analyzing the stability properties of the compact set $\mathcal{A}_{\mathbb{G},\theta}=\mathcal{A}\times \set{0}^{n^2-n}$ under $\mathcal{H}_{\vartheta,s}$. For this last system, we consider the Lyapunov function
\begin{equation}\label{new_lyapunov222}
V_{\mathbb{G}}(\vartheta)=(1-d) \tilde{V}(x) + d\cdot V_{\theta}(\theta),\text{ with }d\in (0,1),~~
\end{equation}
where $\tilde{V}$ is defined as in Lemma \ref{lemma:strict:np}, and $V_\theta(\theta)\coloneqq\frac{1}{2}|\theta|^2$. By using the proof of Lemma \ref{lemma:strict:np}, and noting that $\hat{\mathcal{G}}(\vec{1}\otimes q)= \mathcal{G}(q)$, we obtain:
\begin{equation}\label{dotV:consensus}
     \diffp{\tilde{V}(x)}{x}U(x,h(q))\leq  -\tau_s \tilde{x}^\top M_{\ell}(q,\tau_s)\tilde{x},
\end{equation}
with $\tilde{x}\coloneqq \big((p-q), \mathcal{G}(q)\big)$ and $M_\ell$ given by \eqref{M1} with $c_o = \frac{1}{\ell}$. Under the assumptions of Theorem \ref{theorem:partialInformation} we know that 
\begin{align*}
    0\prec I_n - \left(\frac{T^2}{1-T^2\delta}\right)\frac{\left(\ell I_n-\partial\mathcal{G}(q)^\top\right)\left(\ell I_n-\partial\mathcal{G}(q)\right)}{\ell(1-\eta)-\ell^2\delta},
\end{align*}
and thus that $M_{\ell}(q,\tau_s)\succ \delta I_n~\forall \tau_s\in [T_0,T]$. Hence, from \eqref{dotV:consensus} and letting $\xi(x)\coloneqq \left(\abs{p-q}^2 + \abs{q-q^*}^2\right)^{1/2} $ we obtain that
\begin{align}\label{dotV:consensus:reduced}
     \diffp{\tilde{V}(x)}{x}\dot{x}\leq -T_0\delta\min\set{1,\zeta}\xi^2(x),
\end{align}
where we have also used the bound of Assumption \ref{assumption:reverse:lipschitz}. Additionally, it also follows that
\begin{align}
    &\diffp{\tilde{V}}{x}\left(U(x,\theta + h(q)) - U(x)\right)\leq c_1 \left(\abs{p-q} + \abs{q-q^*}\right)\abs{\theta},\label{dotV:consensus:mixed}\\
    &c_1 \coloneqq 
    \frac{T^2\lambda_{\max}(\mathcal{L})}{\sqrt{2} }\max\set{\frac{1}{T^2}+\frac{4\ell }{T\lambda_{\max}(\mathcal{L})},2+\frac{2\ell}{T\lambda_{\max}(\mathcal{L})}}.\notag
\end{align}
On the other hand, by the fact that the underlying communication graph is undirected and connected, it follows that  $\mathcal{Q}\vec{L}\mathcal{Q}^\top$ is positive definite \cite[Lemma 6]{PavelGames}, and, moreover
\begin{align}\label{dotVtheta:consensus:boundaryLayer}
    \diffp{V_\theta}{\theta}W(q,\theta, 0) \leq -\frac{\lambda_2(\mathcal{L})}{n}\abs{\theta}^2.
\end{align}
We also have that
\begin{equation}\label{dotVtheta:consensus:mixed}
\left(\frac{\partial V_\theta}{\partial x}{-} \diffp{V_\theta}{\theta}\diffp{h}{x}\right)U(x,\theta{+}h(q))\leq c_2\psi(x)\abs{\theta} {+} c_3\abs{\theta}^2,
\end{equation}
where $c_2\coloneqq 2\sqrt{2n}/T_0$ and $c_3\coloneqq 2\sqrt{n}\lambda_{\max}(\mathcal{L})$. Hence, using \eqref{dotV:consensus:reduced}-\eqref{dotVtheta:consensus:mixed} it follows that the time derivative of $V_{\mathbb{G}}$ satisfies $\dot{V}_\mathbb{G}\leq -(\xi(x),\theta)^\top\Lambda_\varepsilon(\xi(x),\theta)$ with
\begin{align*}
\Lambda_\varepsilon \coloneqq
\begin{pmatrix}
        (1-d)T_0\epsilon\min\set{1,\zeta^2} & -\frac{1}{2}(1-d)c_1  - \frac{1}{2}c_2\\ -\frac{1}{2}(1-d)c_1 - \frac{1}{2}c_2 & d\left(\frac{\lambda_2(\mathcal{L})}{\varepsilon n}-c_3 \right)
\end{pmatrix},
\end{align*}
which is positive definite provided that $\varepsilon\in (0,\varepsilon_\delta^*)$ where $\varepsilon_\delta^*$ is as defined in \eqref{upperBoundEpsilon}. Note moreover, that if $\varepsilon$ satisfies this condition there exists $k_\varepsilon>0$ such that
\begin{equation}\label{dotV:consensus:strongDecreaseFlows}
    \dot{V}_{\mathbb{G}} \leq -k_\varepsilon\left(\abs{p-q}^2 + \abs{q-q^*}^2 + \abs{\theta}^2\right).
\end{equation}
Leveraging the results regarding the change of the Lyapunov function $\tilde{V}$ during jumps presented in the proofs of Lemmas \ref{lemma:strict:np}, \ref{lemma:mixed:np} and \ref{lemma:strong:np}, given that (RC$_1$) is satisfied with $\rho_J=0$ by assumption, and since $V_\theta^+(\theta) = V_\theta(\theta)$ for all $\theta$ whenever $\vartheta\in D_{2,\vartheta}$, it follows that $\Delta^{j+n}_j V_{\mathbb{G}}(\vartheta) \leq 0$ for any resetting policy $\alpha\in \set{0,1}^n$.
This inequality and \eqref{dotV:consensus:strongDecreaseFlows} imply that $\mathcal{H}_{\vartheta,s}$ renders the set $\mathcal{A}_{\mathbb{G},\theta}$ UGAS  via \cite[Prop. 3.27]{HDS}. The stability results for $\mathcal{H}_{2,s}$ follow directly by the change of cooordinates $\hat{q} = \theta + h(q)$ and the described result for $\mathcal{H}_{\vartheta,s}$.\hfill$\blacksquare$

\begin{lemma}\label{lemma:H2:absenceOfFiniteEscapeTimes}
Every solution of $\mathcal{H}_2$ is complete.
\end{lemma}
\textsl{Proof:}  Since $\tau$ is restricted to a compact set, it suffices to study the behavior of the states $(q,p,\hat{q})$ or equivalently of $(q,p,\theta)$. Hence, considering the dynamics in \eqref{flowEstimationTheta} and \eqref{flowEstimationWoTheta} it follows that $\abs{\dot{q}}\leq \tilde{\ell}_q\left(\abs{p-q}+\abs{\theta}\right)$, $\abs{\dot{p}}\leq \tilde{\ell}_p\left(\abs{q-q^*}+\abs{\theta}\right)$, and  $\abs{\dot{\theta}}\leq \tilde{\ell}_\theta \left(\abs{p-q}+\abs{\theta}\right)$,
%
where 
$\tilde{\ell}_q \coloneqq \max\set{\frac{2}{T_0},\lambda_{\max}(\mathcal{L})},~~\tilde{\ell}_p\coloneqq 2\ell T\sqrt{N}$ and $\tilde{\ell}_\theta \coloneqq\max\set{\lambda_{\max}(\mathcal{L})\left(\frac{1}{\varepsilon} + \sqrt{N}\right), \frac{\sqrt{N}}{T_0}}$. Using these inequalities we obtain:
\begin{align*}
   \diff{\abs{(q-q^*,~p-p^*,~\theta)}}{t}\leq \abs{(\dot{q},\dot{p},\dot{\theta})}\leq \hat{\ell}\abs{(q-q^*,p-q^*,\theta)},
\end{align*}
with $\hat{\ell}=2\sqrt{3}\max\set{\tilde{\ell}_q,\tilde{\ell}_p,\tilde{\ell}_\theta}$, which by the Gronwall-Bellman inequality implies that the continuous time dynamics of $\mathcal{H}_2$ do not generate finite escape times. Since $G_2(D_2)\subseteq C_2 \cup D_2$, solutions do not stop due to jumps. Therefore, every maximal solution of $\mathcal{H}_2$ is complete. \hfill$\blacksquare$

\vspace{0.1cm}
\textbf{Proof of Theorem \ref{theorem:partialInformation}:} (a) Let $\mathcal{H}_{2,\nu}$ be defined from $\mathcal{H}_2$ by following the same procedure described in the statement of Lemma \ref{lemma:fixedsynchron}. Since the addition of the state $\hat{q}$ and its associated dynamics do not affect the synchronization dynamics, $\mathcal{H}_{2,\nu}$ renders UGFxS the set $\mathcal{A}_{2,\nu}$, where $\mathcal{A}_{2,\nu}$ is as defined in the previous Lemma. Therefore, by the hybrid reduction principle \cite[Cor. 7.24]{HDS}, UGAS of $\mathcal{A}\times\set{\mathcal{Q}(\vec{1}_n\otimes q^*)}$ for system $\mathcal{H}_{2,s}$, established in Lemma \ref{lemma:consensus:UGASsynchro}, implies that $\mathcal{A}\times\set{\mathcal{Q}(\vec{1}_n\otimes q^*)}$ is UGAS for system $\mathcal{H}_{2,\nu}$. Since the choice of $\nu>0$ is arbitrary and since solutions of $\mathcal{H}_{2}$ are complete and bounded, using Lemma \ref{lemma:H2:absenceOfFiniteEscapeTimes} we obtain that the compact set \tcb{$\mathcal{A}\times\set{\mathcal{Q}(\vec{1}_n\otimes q^*)}$} is also UGAS for system $\mathcal{H}_2$.\\
(b) Let $\nu>0$, and $K_0\coloneqq K_x\times K_{\hat{q}}\subset \R^{3n}\times \R^{n^2-n}$ be an arbitrary compact set. Moreover, define $\overline{v}\coloneqq \max_{\vartheta\in K_0} V_\mathbb{G}(\vartheta)$, where $V_\mathbb{G}$ is as defined in \eqref{new_lyapunov222}. Notice that $\overline{v}$ exists since $V_\mathbb{G}$ is continuous and $K_0$ is compact by assumption. It follows that $K_0\subseteq L_{V_\mathbb{G}}\left(\overline{v}\right)$, where $L_{f}\left(c\right)$ represents the $c$-sublevel set of the function $f:\R^m\to \R$. Since $V_\mathbb{G}$ is radially unbounded by construction and Assumption \ref{cocoerciveassumption}, it follows that $L_{V_\mathbb{G}}\left(\overline{v}\right)$ is compact. Let $K^V\coloneqq L_{V_\mathbb{G}}\left(\overline{v}\right)$ and define the HDS $\mathcal{H}_{2,K}=(F_2,~C_2\cap K^V,~G_2,~D_2\cap K^V)$. Notice that under $\mathcal{H}_{2,K}$,  $\hat{q}$ evolves in a compact set. Moreover, by the arguments presented in the proof of item (a), $\mathcal{H}_{2,K}$ renders $K^V$ strongly forward invariant for any $\varepsilon\in (0,\varepsilon^*_\delta)$. Hence, using Lemma \ref{lemma:H2:absenceOfFiniteEscapeTimes}, it follows that, given any arbitrary compact set $\tilde{K}_x\times \tilde{K}_{\hat{q}}\subset K^V$, every solution to $\mathcal{H}_{2,K}$ with $(x(0,0),\hat{q}(0,0))\in\tilde{K}_x\times \tilde{K}_{\hat{q}}$ is complete. Therefore, 
by \cite[Thm. 1]{WangTeelNesic}, for any pair $\hat{t},\hat{j}>0$ there exists $\tilde{\varepsilon}\in(0,\varepsilon_{\delta}^*)$ such that for each $\varepsilon\in (0,\tilde{\varepsilon}]$ and each solution $z$ to $\mathcal{H}_{2,K}$, with $z(0,0)\in K_x\times K_{\hat{q}}$, there exists a solution $x$ to $\mathcal{H}_{1}$
such that $x$ and $z$ are $(\hat{t}, \hat{j}, \nu)-$close. The result follows by using 
$\varepsilon^{**}= \min\set{\tilde{\varepsilon}, \varepsilon_\delta^*}$. \hfill$\blacksquare$

\vspace{0.1cm}\noindent
\textbf{Proof of Theorem 5:} The proof uses tools recently developed for hybrid extremum seeking control \cite{PoTe16,zero_order_poveda_Lina}. Specifically, we show that all the assumptions needed to apply \cite{zero_order_poveda_Lina} are satisfied. In particular, using a Taylor expansion of the form $\phi_i(q+\varepsilon_a\tilde{\mu})\tilde{\mu}_i=\tilde{\mu}_i\phi_i(q)+\varepsilon_a\tilde{\mu}_i\tilde{\mu}^\top\nabla \phi_i(q)+\tilde{\mu}_i\mathcal{O}(\varepsilon_a^2)$, and the fact that $|\tilde{\mu}_i|\leq 1$ for all $i\in\mathcal{V}$, and that $\frac{1}{\tilde{L}}\int_{0}^{\tilde{L}}\tilde{\mu}_i(t)\tilde{\mu}(t)^\top dt=e_i$, where $\tilde{L}=2\pi \text{LCM}\{1/\varsigma_1,1/\varsigma_2,\ldots,1/\varsigma_n\}$ and LCM denotes the least common multiple, the average dynamics of $\mathcal{H}_3$ are precisely given by $\mathcal{H}^A_3=(C_1,F^A_1,D_1,G_1)$, where $G_1,C_1$ and $D_1$ are given by \eqref{jump_map}, \eqref{flow_set0}, and \eqref{jump_set11}, respectively, and $F_1^A$ is given by: 
\begin{equation}
F^A_1(x)=\left(\begin{array}{c} 
2\mathcal{D}(\tau)^{-1}(p-q)\\
-2\mathcal{D}(\tau)\left(\mathcal{G}(q)+\mathcal{O}(\varepsilon_a)\right)\\
\eta\mathbf{1}_n
\end{array}
\right).
\end{equation}
It follows that, on compact sets, we have
\begin{equation}
F_1^A(x)\in \overline{\text{con}}F_1(x+k\varepsilon_a\mathbb{B})+k\varepsilon_a\mathbb{B},
\end{equation}
for some $k>0$, where $F_1$ was defined in \eqref{flowmap00}. Thus, any solution of the average dynamics $\mathcal{H}^A_3$ is also a solution of an inflated HDS generated from $\mathcal{H}_1$. By \cite[Thm. 7.21]{HDS}, we conclude that, under the Assumptions of Theorems \ref{theorem1}-\ref{theorem:strong:np}, system $\mathcal{H}^A_3$ renders SGPAS as $\varepsilon_a\to0^+$ the compact set $\mathcal{A}$. Since $\mathcal{H}^A_3$ and $\mathcal{H}_1$ are nominally well-posed, all the assumptions needed to apply \cite[Thm.7]{zero_order_poveda_Lina} are satisfied, and we can conclude that $\mathcal{H}_3$ renders SGPAS as $(\varepsilon_p,\varepsilon_a)\to0^+$ the compact set $\mathcal{A}\times \mathbb{T}^n$. Item (b) follows directly by \cite[Prop. 6]{zero_order_poveda_Lina}. \hfill$\blacksquare$

\vspace{-0.1cm}
\section{CONCLUSIONS}
\label{sec_conclusions}
In this paper, we introduced a class of Nash set-seeking algorithms with dynamic momentum for the efficient solution of non-cooperative games with finitely many players. The algorithms are modeled by hybrid dynamical systems that incorporate continuous-time dynamics with momentum and discrete-time coordinated resets. We developed model-based algorithms, as well as algorithms suitable for games with partial information and model-free settings where players have access only to measurements of their cost. In each case, we established robust stability and accelerated convergence properties using multi-time scale techniques for hybrid dynamical systems.

\appendices
\bibliographystyle{ieeetr}
\bibliography{Bibliography}

\end{document}